\documentstyle[amssymb,11pt]{article}
\makeatletter
\renewcommand\section{\@startsection{section}{1}{\z@}%
                                     {-6.25ex\@plus -1ex \@minus -.2ex}%
                                     {1.5ex \@plus .2ex}%
                                     {\reset@font\large\bfseries}}
\renewcommand\subsection{\@startsection{subsection}{2}{\z@}%
                                     {3.25ex \@plus1ex \@minus.2ex}%
                                     {-1em}%
                                     {\reset@font\normalsize\bfseries}}
\@addtoreset{equation}{section}
\makeatother
\renewcommand{\^}{\hspace*{-3pt}{}^} 
\renewcommand{\_}{\hspace*{-3pt}{}_}
\renewcommand{\theequation}{\thesection.\arabic{equation}}
\newtheorem{theo}{Theorem}      \newtheorem{defn}{Definition}
\newtheorem{lemma}{Lemma}       \newtheorem{prop}[theo]{Proposition}   
\def\sg{{\scriptstyle {\frak G}}}   \def\End{{\mbox{\it End\/}}}
\def\mod{{\mbox{\it mod\/}}}        \def\tr{{\mbox{\it tr\/}}}
\def\:{:\,=}
\def\sd{\times\hspace*{-9.7pt}\rule{0.25pt}{5.5pt}\hspace*{10.5pt}}
\def\bo{\mbox{\,\raisebox{-0.65mm}{$\Box$} \hspace{-5.3mm}
${\scriptstyle\times}$ \/}}                
\def\be{\begin{equation}\label}      \def\ee{\end{equation}}   
\def\ba{\begin{eqnarray}}            \def\ea{\end{eqnarray}}
\def\ar{\begin{array}}               \def\er{\end{array}}
\def\nn{\nonumber}                   \def\ns{\normalsize}
\def\rar{\rightarrow}                \def\o{\otimes} 
\def\wt{\widetilde}                  \def\wh{\widehat}
\def\H{{\cal H}}        \def\cM{{\cal M}}    
\def\G{{\cal G}}        \def\F{{\cal F}}     \def\K{{\cal K}}      
\def\N{{\cal N}}        \def\S{{\cal S}}     \def\J{{\cal J}}
\def\C{{\cal C}}        \def\W{{\cal W}}     \def\Z{{\cal Z}}
\def\RR{{\cal R}}       \def\V{{\cal V}}     \def\P{{\cal P}} 
\def\R{{\cal R}}         
\def\vp{\varpi }        \def\vth{\vartheta}     \def\vs{\varsigma}
\def\ve{\varepsilon}    \def\s{\sigma}          \def\t{\tau}
\def\a{\alpha }         \def\e{\epsilon}        \def\k{\kappa}
         \def\ga{\gamma}         
\def\dl{\delta}         \def\D{\Delta }      
        \def\rf{{\sf f}}         \def\rv{{\sf v}}
\def\re{{\sf e}}                 
\def\rw{{\sf w}}        \def\rk{{\sf \kappa}}
\def\sp{{\sf p}}                 \def\sF{{\sf F}}
\newcommand{\se}[1]{\stackrel{\scriptscriptstyle #1}{\sigma}}
\newcommand{\Ph}[2]{\stackrel{\scriptscriptstyle #1}{\Phi}
            \hspace*{-4mm} \phantom{\Phi}^{#2}}
\newcommand{\up}[2]{\stackrel{\scriptscriptstyle #1}{#2}}
\def\ew{\hspace*{-2mm}}   \def\ppe{\hspace*{-2.5mm}}
\newcommand{\Fus}[6]{F_{{\scriptstyle #1}{\scriptstyle #2}}
  \hspace*{-.0mm}[ \ew \begin{array}{ll} {\scriptstyle #3 }
  \ppe & {\scriptstyle #4} \ppe \\[-2mm] {\scriptstyle #5}\ppe &
  {\scriptstyle #6}\ew \end{array} ]}
\newcommand{\CG}[6]{[ \ew \begin{array}{lll} {\scriptstyle #1} \ppe
  & {\scriptstyle #2} \ppe & {\scriptstyle #3} \ew \\[-2mm] {\scriptstyle
  #4} \ppe & {\scriptstyle #5}\ppe & {\scriptstyle #6} \ew\end{array}  ]}
\newcommand{\SJS}[6]{ \{ \ew \begin{array}{lll} {\scriptstyle #1} \ppe  &
  {\scriptstyle #2} \ppe & {\scriptstyle #3}
  \ppe \\[-2mm]{\scriptstyle #4}  \ppe & {\scriptstyle #5} \ppe &
 {\scriptstyle #6} \ew \end{array}  \} }
\newcommand{\Br}[7]{B^{#1}_{{\scriptstyle #2}{\scriptstyle #3}}
  \hspace*{-.0mm}[ \ew \begin{array}{ll} {\scriptstyle #4 }
  \ppe & {\scriptstyle #5} \ppe \\[-2mm] {\scriptstyle #6}\ppe &
  {\scriptstyle #7}\ew \end{array} ]}
\newcommand{\vvert}[3] {(\ew \!\begin{array}{ll} {\ \scriptstyle #1}\!
 & \!{\scriptstyle #3} \\ [-2mm] & \! \ppe {\scriptstyle #2}
 \ew \end{array} \!\! )}
\begin{document}
\hyphenation{ope-ra-tor ope-ra-tors ge-ne-ra-tor ge-ne-ra-tors
 mo-du-lar}
\thispagestyle{empty}
November~~1996
\hfill{\begin{tabular}{l}  \tt PDMI~-~3/1996 \\
    \tt Berlin Sfb288~No.238 \\
    \tt q-alg/9611010
       \end{tabular} }
\vspace*{1.3cm}

\begin{center}
{\Large\bf Vertex Operators -- {}From a Toy Model to
Lattice Algebras} \\ [1.1cm]
{\sc Andrei G. Bytsko}~${}^*$ \\ [1mm]
Institut f\"ur Theoretische Physik, Freie Universit\"at Berlin \\
Arnimallee 14, 14195 Berlin, Germany \\
and \\
Steklov Mathematical Institute, Fontanka 27,\\
 St.Petersburg~~191011, Russia \\ [3mm]
{\sc Volker Schomerus}~${}^{**}$ \\   [1mm]
II. Institut f\"ur Theoretische Physik, Universit\"at Hamburg \\
Luruper Chaussee 149,  22761~~Hamburg, Germany  \\ [1.5cm]
{\bf Abstract} \\ [2.9mm]

\parbox{13cm}{
Within the framework of the discrete Wess-Zumino-Novikov-Witten
theory we analyze the structure of vertex operators on a lattice.
In particular, the lattice analogues of operator product expansions
and braid relations are discussed. As the main physical application,
a rigorous construction for the discrete counterpart $g_n$ of the
group valued field $g(x)$ is provided. We study several automorphisms
of the lattice algebras including discretizations of the evolution
in the WZNW model. Our analysis is based on the theory of modular
Hopf algebras and its formulation in terms of universal elements.
Algebras of vertex operators and their structure constants are
obtained for the deformed universal enveloping algebras $U_q(\sg)$.
Throughout the whole paper, the abelian WZNW model is used as a
simple example to illustrate the steps of our construction.
}  \end{center}

\vfill{ \hspace*{-9mm}
 \begin{tabular}{l}
\rule{6 cm}{0.05 mm} \\
${}^*$ \ e-mail:\ bytsko@physik.fu-berlin.de,\ \  bytsko@pdmi.ras.ru \\
${}^{**}$  e-mail:\ vschomer@x4u.desy.de
 \end{tabular} }
\newpage
\setcounter{page}{1}
\tableofcontents

\newpage
\section{INTRODUCTION}

\noindent
{\it Quantization of the WZNW model.}\
The Wess-Zumino-Novikov-Witten (WZNW) model \cite{WeZu,Nov,Wit,KnZa} is one
of the most famous examples of a rational conformal field theory (CFT)
\cite{BPZ,KaTs,MoSe}. On the classical level it describes some time
evolution for a field $g(x)$ mapping points $x$ of the circle ${\bf S}^1$
into a compact Lie group $G$. Among the dynamical variables of the theory,
the currents $j^r(x)=g^{-1} \partial_- g$, $j^l=(\partial_+ g) g^{-1}$
are of particular interest. In contrast to the field $g$, the currents
$j^r$ and $j^l$ are chiral, so that $\partial_+ j^r = 0 $ and 
$\partial_- j^l = 0$. Moreover, their Poisson structure is well known to 
give rise to two commuting
copies of Kac-Moody (KM) algebras. Even though numerous papers have
been devoted to the quantization of the WZNW-model (e.g. \cite{Blo,Fa1,
AlSh,Gaw,BDF,CGO,CG1,FHT}), a rigorous construction of the continuum
theory (which requires field strength renormalization) is not fully
understood. This motivates the search for lattice regularizations
of the theory (i.e., the circle ${\bf S}^1$ is replaced by a
periodic lattice with lattice spacing $a$) which preserve much
of the symmetry structure of the continuum WZNW-model. One may
construct appropriate discretizations of the classical model
(i.e., $\hbar = 0$) first and then quantize the classical lattice
theory to obtain a well defined discrete quantum theory
(i.e., $\hbar \neq 0, a\neq 0$). Investigation of the latter is
expected to provide insights into the structure of the continuum
model. A final step would involve performing the limit $a \rar 0$
while keeping $\hbar \neq 0$.

The realization of this program was started in \cite{AFS,AFSV,Fa2} where
a lattice regularization of the Kac-Moody algebra has been proposed.
Classical and quantum lattice current algebras were further investigated
in \cite{FaGa,AFFS}. Our aim here is to extend the analysis of \cite{AFFS}
by introducing {\em chiral vertex operators}. In comparison with the
current algebra, the algebras of vertex operators contain
(a finite number of) additional generators. Within these larger
algebras we will be  able to prepare a discrete analogue of the group
valued field $g(x)$ by combining left- and right chiral vertex
operators.
\vspace*{4mm}

\noindent
{\it Quantum symmetry structure of the WZNW model.}\
Let us recall that solutions of the classical Yang-Baxter equation
appear already in the Poisson structure of the classical lattice
current algebras (see \cite{FaGa} and references therein).
After quantization, quantum groups and quantum universal
enveloping algebras $\G$ \footnote{For shortness, we will often
refer to $\G$ as {\em quantum algebra}.} are expected to emerge.
Throughout this paper we will meet (global and local) objects
(the monodromies $M^\a$ and a discrete field $N_n$, see below)
whose nature reflects a quantum algebraic structure as well as
an object (the discrete field $g_n$) which displays features of a
quantum group. The corresponding deformation parameter is of the
form $q=exp\{i\gamma \hbar \}$ with $\gamma=\pi/(k+\nu)$ (where
$k$ is the level of the KM algebra and $\nu$ is the dual Coxeter
number of $G$) and does not depend on the lattice spacing $a$.
Therefore, the quantum group structures of the continuum and
lattice WZNW model coincide. It is also worth mentioning that
some aspects of the quantum symmetry structure survive reductions
to other theories so that part of what we describe below may be
compared with studies of the quantum Liouville and Toda models
\cite{GeNe,FaTa,Bab,CrGe}.
\vspace*{4mm}

\noindent
{\it Remarks on lattice current algebras.}\
Before we summarize our results, let us briefly  review
the discretization used in \cite{AFFS} for the chiral currents
$j^l(x), j^r(x)$. Recall that the latter are Lie-algebra valued
fields which depend periodically on the variable $x$. Instead of
working with these standard variables, we prefer to pass to the
fields $j^r(x)$ and $ \eta(x) = j^r(x)-j^l(x)$ and describe their
lattice counterparts. Our lattice divides the circle into $N$
links of length $a = 2\pi/N$. So there are $N$ vertices at the
points $x = a n$ which are numbered by $n = 0, \dots, N-1,
N\equiv 0$ and the $n^{th}$ link runs from the $(n-1)^{st}$ vertex to
the $n^{th}$. We may discretize the field $\eta(x)$ by the simple
prescription $\eta_n := \int_{(n-\frac 12)a}^{(n+ \frac 12)a} \eta(x) dx 
\,=\, a \eta(na) + O(a^2)$ so that the lattice field
$\eta_n$ has values in a tensor product of $N$ copies of the Lie
algebra which are assigned to the $N$ vertices on the lattice. For the
right chiral current $j^r(x)$ the discretization scheme is different.
In this case we encode the information about the field in the
holonomies along links, i.e., we define the lattice field $j^r_n$ by
$$  j^r_n \ :=\ P exp(\int_n j^r(x) dx) \ \ . $$
Here $\int_n$ denotes integration along the $n^{th}$ link.
By construction,
this classical lattice field $j^r_n$ has values in the Lie group.
The rather different treatment of the fields $\eta(x)$ and $j^r(x)$
may be understood from the Poisson structure of the classical
theory, which is ultralocal for $\eta(x)$ but contains terms
proportional to $\dl'(x-y)$ if the field $j^r(x)$ is involved
(see \cite{AFS,FaGa,AFFS}).

When we pass to the quantum theory, the functions on the
space of field configurations become operators and generate
some non-commutative algebra $\K_N$. More concretely, the
algebra $\K_N$ is generated from the quantum lattice
fields $J^r_n, N_n$ which correspond
to the classical fields $j^r_n, \eta_n$ described above. We
review the explicit definition of {\em lattice current
algebras} in Section~4. Let us only mention here that a very
elegant formulation for commutation relations of the quantum
operators can be given in the $R$-matrix language.

In mathematical terms, one has to regard the quantum fields
$N_n$ and $J_n = J^r_n$ as objects in the tensor product $\G_a
\o \K_N$ of the deformed universal enveloping algebra $\G_a
= U_q(\sg)$ with the lattice current algebra $\K_N$.
We can understand this by looking at the classical lattice
field $j^r_n$, for instance. It was constructed as the
holonomy of the Lie-algebra valued field $j^r(x)$ and may be
evaluated with irreducible representations of the Lie algebra.
Let us denote such representations by $\t^I$ and introduce the
symbols $V^I, \dl_I$ for their carrier spaces and dimensions,
respectively. Then
we see that $j_n = j^r_n$ gives rise to $\dl_I\times\dl_I$-matrices
$j^I_n$ of dynamical variables. Accordingly, the corresponding
quantum operators $J^I_n$ are matrices of generators for $\K_N$
which is to say that $J^I_n \in \End(V^I) \o \K_N$. All these
objects $J^I_n$ may be assembled back into one {\em universal
element} $J_n \in \G_a \o \K_N$. More details will be presented
later; we anticipated this heuristic discussion of universal
elements only to prepare for some formulae below.

One of the main aims in \cite{AFFS} was to develop a complete
representation theory for lattice current algebra $\K_N$. It
turned out that $\K_N$ possesses a family of irreducible
$*$-representations on vector spaces $W^{IJ}_N$ with
labels $I,J$ running
through classes of finite-dimensional, irreducible
representations of $U_q(\sg)$. {\em Two} such labels are needed
because of the {\em two} chiralities in the current algebra.
Furthermore, the algebra $\K_N$ was found to admit two families
of local co-actions $\Lambda_n^r, \Lambda^{l}_n: \K_N \mapsto
\G_a \o \K_N$ of the Hopf algebra $\G_a$. They may be considered
as a special case of the more general lattice fusion products
in \cite{AFFS} and give rise to a notion of tensor products for
representations of $\K_N$ (see also \cite{NiSz} for related
results).
\vspace*{4mm}

\noindent
{\it Vertex operators on a lattice.}\
Product structures in the representation theory are precisely
what is needed to initiate a theory of vertex operators. More
technically, we employ the homomorphisms $\Lambda_n^r, \Lambda_n^{l}$
in extending the lattice current algebra $\K_N$ by chiral vertex
operators $\Phi_n^{r}, \Phi_n^{l}$ so that the following {\em
intertwining relations} hold for both chiralities $\a = r,l$,
\be{vint}
      A \ \Phi^{\a}_n \,=\, \Phi^{\a}_n \ \Lambda^{\a}_n (A) \ \ \
     \ \ \mbox{ for all } \ \ \ \ A \ \in \ \K_N \ \ .
\ee
The elements $\Phi^\a_n$ generate an extension $\W_N$ of the
lattice current algebra $\K_N \subset \W_N$. Since
$\Lambda^{\a}_n(A)$ is an element of $\G_a \o \K_N$
and hence also of the extension $\G_a \o \W_N$, the
product on the r.h.s. of (\ref{vint}) is well defined for
$\Phi^{\a}_n \in \G_a \o \W_N$.
On the l.h.s., $A = e \o A \in \G_a \o \K_N$
with $e \in \G_a$ being the unit element.

Our vertex operators $\Phi^{\a}_n$ on the lattice possess
a number of properties which are all closely related to properties
of vertex operators in the continuum theory. Let us highlight
some of them without going into a detailed discussion.\footnote{
A construction of (non-chiral) vertex operators for infinite open
lattices has been suggested in \cite{NiSz}. Some properties of these 
vertex operators are similar to what we shall consider here. However, 
these are different structures, in particular, because for a finite 
lattice the current algebra $\K_N$ has a non-trivial center $\C$. 
An action of our vertex 
operators on elements from $\C$ will play a crucial role in the theory.}

\begin{enumerate}
\item Lattice vertex operators $\Phi^{\a}_n$ at a fixed lattice
      site obey {\em operator product expansions} of the form
     \be{latOPE}
      \Ph{2}{r}_n \ \Ph{1}{r}_n \, = \, F_r \ \D_a (\Phi^{r}_n)
       \ \ \ \ \mbox{ and } \ \ \ \
      \Ph{1}{l}_n \ \Ph{2}{l}_n \,=\,  F_l \ \D_a (\Phi^{l}_n) \ .
     \ee
      As usual, the notation $\Ph{1}{\a}_n$ means that we regard
      the vertex operator $\Phi^{\a}_n$ as an element of $\G_a \o
      \G_a \o \W_N$ with trivial entry in the second tensor factor
      etc. We have also used the shorthand $\D_a (\Phi^{\a}_n) =
      (\D \o id)(\Phi^{\a}_n) \in \G_a \o \G_a \o \W_N$ for the
      action of the co-product on the first tensor factor of
      $\Phi^\a_n$. The objects $F_\a$ are analogues of the {\em
      fusion matrix} in the continuum theory. We describe their
      general properties and, in particular, their relation with
      $6j$-symbols in Section 2.
\item Lattice vertex operators $\Phi^{\a}_n$ assigned to different
      lattice sites obey {\em braid relations}
     \be{lbraid}
      \Ph{1}{r}_n  \ \Ph{2}{r}_m \, =\, \R^r_- \ \Ph{2}{r}_m\ \Ph{1}{r}_n
      \ \ \ \ \mbox{ and } \ \ \ \
      \Ph{2}{l}_n \  \Ph{1}{l}_m \, =\,
      \R^l_+\  \Ph{1}{l}_m \  \Ph{2}{l}_n
     \ee
      for all $0 \leq n < m < N$. Here the objects $\R_\pm$
      play the role of the {\em braiding matrix} in the continuum
      theory. Let us add that lattice vertex operators of different
      chirality commute for all $n,m$. Furthermore, $\Phi^\a_n$ commute
      with $N_m$ for $m \neq n$ and with $J_m$ for $m \neq n,n+1$,
      that is, the vertex operators have local exchange relations
      with elements of the current algebra.
\item Lattice vertex operators $\Phi^{\a}_n, \a = r,l,$
      satisfy the following {\em difference equation }
     \be{ldiff}
       \Phi^{\a}_{n+1}\, =\, \Phi^{\a}_n\  J^{\a}_{n+1} \ \ .
     \ee
      In the naive continuum limit, we have $J^{\a}_n = e\o e - a
      J^{\a}(x)+O(a^2)$ with $x = a n $ and the  difference equation
      becomes a differential equation which expresses $\partial_x
      \Phi^{\a}(x)$ as a (normal ordered) product of  $\Phi^{\a}(x)$
      and $J^{\a}(x)$. Such an equation is well known for the
      quantized continuum theory.
\end{enumerate}
As one may infer from the third property in this short list, lattice
vertex operators (much like their continuum counterparts) cannot be
periodic. Indeed, starting from $\Phi^{\a}_0$ an iterated application
of eq.~(\ref{ldiff}) gives
$$
      \Phi^{\a}_{N}\, \equiv\,  \Phi^{\a}_0\  M^{\a} \ \ \ \mbox{ with }
       \ \ \ \ M^{\a}\, =\,  J_1^{\a} \dots J_N^\a \ \ .
$$
The objects $M^{\a}, \a = r,l,$ are called {\em chiral monodromies}.
Actually, the lattice rotation $n \mapsto n+N$ gives rise to an
inner automorphism of the algebra of vertex operators which acts
trivially on the lattice fields $J^\a_n$ and  $N_n$. We show in
Section 6 that this automorphism can  be generated by conjugation
with a unitary element ${\sf v}$. The latter is constant on the
irreducible representation spaces $W^{IJ}_N$ of the lattice current
algebra $\K_N$ and its value
$$      {\sf v}^{IJ} = e^{2 \pi i( h_J - h_I)}\ \ \ \  $$
can be expressed in terms of the conformal dimensions $h_I$ of the
WZNW model. This leads us to identify ${\sf v}$ with the operator
$\exp\{2 \pi (L_0 - \bar{L}_0)\}$ which generates rotations by
$2\pi$ in the continuum theory. In the lattice theory ${\sf v}$ is
obtained from quantum traces of chiral monodromies $M^\a$ and
is related to the ribbon element of $U_q(\sg)$.

It will be shown in Sections 3 and 5 that the field $\S_a(\Phi^l_n)
\Phi^r_n$ \footnote{We use the notation $\S_a(\Phi^l_n) =
(\S \o id) (\Phi^l_n)$ with $\S$ being the antipode of $\G$.}
can be restricted to the {\em diagonal subspace} $\bigoplus_K
W_{\,N}^{\bar{K} K}$. Let us denote this restriction by $g_n$ which
suggests that it is a quantum lattice analogue of the group
valued field $g(x)$ in the WZNW model. In fact, our analysis
will reveal that $g_n$ is a local and quantum group valued
field, i.e.
\be{gfield}
 \up{1}{g}_n \ \up{2}{g}_m \,=\   \up{2}{g}_m \ \up{1}{g}_n \ 
 \ \ ( n \neq m ) \ \ \
 \mbox{ and }  \ \ \ \ R \ \up{2}{g}_n \ \up{1}{g}_n  \,=\
 \up{1}{g}_n\  \up{2}{g}_n\,  R \ \ .
\ee
Moreover, $g_n$ turns out to be periodic. In contrast to the
chiral currents $j^\a(x)$, the time evolution of the group
valued field $g(x)$ is described by a nontrivial second order
differential equation. Its discrete analogue is discussed in Section 6.

Before we address the full lattice theory we explain some basic
constructions in a simple toy model (cf. Section 3). Here one
studies the algebra generated by the  monodromies $M^r,M^l$
instead of the whole (lattice) current algebra and universal
(deformed) tensor operators for the quantum algebra $\G$ as
simple examples of vertex operators \cite{AF1,CG2,AF2,ByFa}.
This finite dimensional toy model may be regarded as a special
case of the discrete WZNW theory where $N = 1$ and it describes
the non-local degrees of freedom for an arbitrary number $N$ of
lattice sites. Objects and relations of the toy model admit for a
nice pictorial presentation which, in particular, brings new
light into the {\em shadow world} \cite{KR1}.

\section{HOPF ALGEBRAS AND VERTEX OPERATORS}
\setcounter{equation}{0}

\subsection{Semi-simple modular Hopf algebras.} \
By definition, a Hopf algebra is a quadruple $(\G,\e,\D,\S)$ of
an associative algebra $\G$ (the ``{\it symmetry algebra}\,'')
with unit $e \in \G$,
a one-dimensional representation $\e: \G \mapsto {\Bbb C}$
(the ``{\it co-unit}\,''), a homomorphism $\D: \G \mapsto \G \o \G$
(the ``{\it co-product}\,'') and an anti-automorphism $\S: \G \mapsto\G$
(the ``{\it antipode}\,''). These objects obey a set of basic axioms
which can be found, e.g., in \cite{Abe,Swe}. The Hopf algebra
$(\G,\e,\D,\S)$ is called quasi-triangular if there is an invertible
element $R \in  \G \o \G$ such that \vspace*{-2mm}
\ba
     R\   \D (\xi) & = & \D'(\xi) \ R \
      \ \mbox{ for all } \ \ \xi \in \G\ , \nn \\[1mm]
     (id \o \D)(R) = R_{13}R_{12} \ \ & , & \ \
     (\D \o id)(R) = R_{13}R_{23} \ \ .  \nn
\ea
Here $\D'(\xi) = P \D(\xi)P $, with $P$ being the permutation, i.e.,
$P (\xi \o \eta) P = \eta \o \xi$ for all $\xi,\eta \in \G$,
and we are using the standard notation for the elements $R_{ij}
\in \G \o \G \o \G$.

For a ribbon Hopf-algebra one postulates, in addition, the existence of
a certain invertible central element
$v \in \G$ (the ``{\it ribbon element}\,'')
which factorizes $R'R \in \G \o \G$ ( here $R' = PRP$), in the sense that
\be{rib}
 R'R \,=\, (v \o v)\, \D(v^{-1})\ , \ \ \ \ \S(v)\,=\,v \ ,
 \ \ \ \ \e(v) \,=\, 1
\ee
(see \cite{ReTu,KR2} for details).
We want this structure to be consistent with a $*$-operation
on $\G$. To be more precise, we require that
\footnote{We fix $*$ on $\G \o \G$ by $(\xi \o \eta)^* =\xi^* \o\eta^*$.
Following \cite{MaSc}, we could define an alternative involution
$\dagger$ on $\G \o \G$ which involves a permutation of
components, i.e., $(\xi \o \eta)^\dagger = \eta^\dagger \o
\xi^\dagger$ and $\xi^\dagger = \xi^* $ for all $\xi, \eta \in \G$.
With respect to $\dagger$, $\D$ becomes an ordinary $*$
-homomorphism and $R$ is unitary.}
\be{stop}
 R^* = (R^{-1})' = PR^{-1}P \ ,\ \ \ \D(\xi)^* = \D'(\xi^*)\ , \ \ \
 v^* \,=\, v^{-1}\ .
\ee
This structure is of particular interest, since it appears in the theory
of the quantized universal enveloping algebras $U_q(\sg)$ when the
complex parameter $q$ has values on the unit circle \cite{MaSc}.

At this point we assume that {\em $\G$ is
semi-simple}, so that every finite dimensional representation of $\G$ can
be decomposed into a  direct sum  of finite dimensional, irreducible
representations. {}From each equivalence class $[I]$ of irreducible
representations of $\G$, we may pick a representative $\t^I$, i.e.,
an irreducible representation of $\G$ on  a $\dl_I$-dimensional Hilbert
space $V^I$. The {\em quantum trace} $\tr_q^I$ is a  linear functional
acting on elements $X \in \End (V^I)$ by
$$     \tr^I_q (X) = \tr^I( X \t^I(w))\ \ . $$
Here $\tr^I$ denotes  the standard trace on $\End (V^I)$ with
$\tr^I (e^I) = \dl_I$ and $w\in\G$
is a distinguished group-like element constructed from
the ribbon element $v$ and the element $R$ by the formula
$ w^{-1} = v^{-1} \sum \S(r^2_\vs) r^1_\vs$,
where the elements $r^i_\vs$ come from the expansion
$R = \sum r^1_\vs \o r^2_\vs$.

Evaluation of the unit element $e^I \in \End (V^I)$ with $\tr^I_q$ gives
the {\em quantum dimension}, $d_I := \tr^I_q(e^I)$, of the representation 
$\t^I$. Furthermore, we  assign a  number
${\sf S}_{IJ}$ to every pair of representations $\t^I,\t^J$:
$$    {\sf S}_{IJ} :=  {\vartheta} (\tr^I_q \o \tr^J_q) (R'R)^{IJ}
  \ \ \ \mbox{ with } \ \ (R'R)^{IJ} =(\t^I \o \t^J)(R'R) \ \ , $$
with a suitable, real normalization factor ${\vartheta}$.
The numbers ${\sf S}_{IJ}$ form the so-called {\em $S$-matrix}
${\sf S}$. {\it Modular} Hopf algebras are ribbon
Hopf algebras with an  invertible $S$-matrix.\footnote{If
a diagonal matrix ${\sf T}$ is introduced according to
${\sf T}_{IJ} = \varpi \dl_{I,J} d_I^2 \t^I(v)$ (with an appropriate
complex factor $\varpi$), then ${\sf S}$ and ${\sf T}$ furnish a
projective representation of the modular group $SL(2,{\Bbb Z})$.}

Let us finally recall that the  tensor
product, $\t \bo \t'$, of two representations $\t,\t'$ of a Hopf
algebra is defined by
$$ (\t \bo  \t')(\xi) = (\t \o \t') \D (\xi) \ \ \mbox { for all }
 \ \ \xi \in \G\ \ . $$
In particular, one may construct the tensor product $\t^I \bo \t^J$
of  two irreducible representations. According to  our assumption
that $\G$ be semi-simple, such tensor products of representations can
be decomposed into a direct sum  of irreducible representations.

Among all our assumptions on the structure of the Hopf-algebra $(\G,
\e,\D,\S)$ (quasi-triangularity, existence  of a ribbon  element $v$,
semi-simplicity of $\G$ and invertibility of $S$-matrix $\sf S$),
semi-simplicity of
$\G$ is the most problematic one. In fact it is violated by the
algebras $U_q(\sg)$ when $q$ is  a root of unity. It is
sketched in \cite{AGS} how ``truncation'' can cure this problem,
once the theory has been extended to weak quasi-Hopf algebras
\cite{MaSc}.
\vspace*{5mm}

\noindent
{\sc Example:} ({\em Modular Hopf-algebra $\Z_q$} \cite{AFFS})
We wish to give one fairly trivial example for the
algebraic structure discussed so far which comes from
the  group ${\Bbb Z}_p$. To be more precise, we consider
the  associative algebra $\Z_q$ generated by one element $h$ subject
to the relation $h^p = e$. Co-product, co-unit and antipode
for this algebra can be defined by
$$ \D(h) = h \o h \ \ , \ \ \S (h) =  h^{-1} \ \ , \ \
 \e(h) = 1 \ \  . $$
We observe that  $\Z_q$ is a commutative semi-simple algebra. It has
$p$ one-dimensional representations $\t^t(h) = q^t, t = 0, \dots, p-1,$
where $q$ is a root of unity, $q = e^{2 \pi i/p}$. We may
construct characteristic projectors $P^t \in \Z_q$ for these
representations according to
\be{zP}
 P^t = \frac{1}{p} \sum_{m=0}^{p-1}  q^{-tm}  h^m\ \ \mbox{ for }
\ \ t= 0, \dots, p-1\ \  .
\ee
One can easily check that $\t^t(P^s) = \delta_{t,s} $.
The elements $P^t$ are employed to obtain a nontrivial $R$-matrix:
\be{RZ}
 R = \sum_{t,s=0}^{p-1} q^{ts} P^t \o P^s \  .
\ee
It is easy to see that
$(\t^t \o \t^s)R = q^{ ts}$. The $R$-matrix satisfies all
the  axioms stated above and thus turns $\Z_q$ into a quasi-triangular
Hopf algebra. Moreover, a ribbon element is given by  $v =
\sum q^{-t^2} P^t$.

It is natural to introduce a $*$-operation on $\Z_q$ such that
$h^*= h^{-1}$. The relations (\ref{stop}) hold due to
the co-commutativity of $\D$, i.e., $\D' = \D$, and
the property $R = R'$. A direct computation shows that the $S$-matrix
$\sf S$ is invertible only for  odd integer $p$.
Summarizing all this, the algebra $\Z_q, q= exp(2\pi i /p) $ is a
semi-simple ribbon Hopf-$*$-algebra. It is a modular Hopf algebra for 
all odd integer $p$. The reader is invited to check that for $\Z_q$ 
the quantum trace $tr^I_q$  coincides with the standard one.

\subsection{Universal elements and $R$-matrix formalism.}
Modular Hopf algebras admit a very elegant $R$-matrix description.
For its presentation, let us introduce another (auxiliary)
copy, $\G_a$, of $\G$ and let us consider the $R$-matrix as an
object in $\G_a \o \G$. To distinguish the latter clearly from
the usual $R$, we denote it by
$$ N_+ \equiv R \in \G_a \o \G \ , \ \ \ \ \
 N_- \equiv  (R')^{-1} \in \G_a \o \G \ . $$
At the same time let us introduce the standard symbols $R_+ =
R$ and $R_- = (R')^{-1} \in \G_a \o \G_a$.  Quasi-triangularity
of the $R$-matrix furnishes the relations
\ba \label{Npmeq}
 \D_a(N_{\pm}) \,=\, \up{1}{N}_{\pm}\, \up{2}{N}_{\pm} \ \ & , & \ \
 R_+ \,\up{1}{N}_+ \,\up{2}{N}_- \,=\,\up{2}{N}_- \,\up{1}{N}_+ \,R_+\ \ ,
 \\[1mm]  R_+ \,\up{1}{N}_\pm \, \up{2}{N}_\pm  & = &
   \up{2}{N}_\pm \, \up{1}{N}_\pm \, R_+ \ \ . \nn
\ea
Here we use the same notations as in the introduction, and
$\D_a (N_\pm) = (\D \o id) (N_\pm) \in \G_a \o \G_a \o \G$.
The subscript $\ _a$ reminds us that $\D_a$ acts on the
auxiliary (i.e., first) component of $N_\pm$. To be perfectly
consistent, the objects $R_\pm$ in the preceding equations should
all be equipped with a lower index $a$ to show that $R_\pm \in \G_a
\o \G_a$ etc. We hope that no confusion will arise from omitting
this subscript on $R_\pm$. The equations (\ref{Npmeq}) are somewhat
redundant: in fact, the exchange relations in the second line
follow from the first equation in the first line. This
underlines that the formula for $\D_a(N_\pm)$ encodes information
about the product in $\G$ rather than the co-product.%
\footnote{The co-product $\D$ of $\G$ acts on $N_\pm$ according
to $\D(N_\pm)= (id \o \D)(N_\pm) = N_\pm'  N_\pm'' \in
\G_a \o \G \o \G$. Here $N_\pm'$ and $N_\pm''$
have the unit element $e \in \G$ in
the third and second tensor factor, respectively.}

Next, we combine $N_+$ and $N_-$ into one element
$$   N \,:=\, N_+ (N_-)^{-1}   \ \ \in \ \G_a \o \G\ \ .$$
{}From the properties of $N_\pm$ we obtain an expression for the
action of $\D_a$ on $N$,
\be{DN}
 \ar{c}
  R_+\, \D_a(N) \, = \, R_+\,\up{1}{N}_+ \up{2}{N}_+ \up{2}{N}\^{-1}_-
             \up{1}{N}\^{-1}_- \nn \, = \,
 \up{2}{N}_+ \up{1}{N}_+ R_+ \, \up{2}{N}\^{-1}_-
             \up{1}{N}\^{-1}_- \,=\, \\ [2mm]
 \, = \, \up{2}{N}_+ \up{2}{N}\^{-1}_- R_+ \up{1}{N}_+
             \up{1}{N}\^{-1}_-
         \, = \, \up{2}{N}  R_+  \up{1}{N}  \ . \er
\ee
As seen above, the formula for $\D_a(N)$ encodes relations in the
algebra $\G$ and  implies, in particular, the following exchange
relations for $N$:
\be{NN}
    R_-^{-1} \up{2}{N} R_+ \up{1}{N}
    \,=\,  R_-^{-1} R_+ \D_a (N)  \,=\,
    R_-^{-1} \D'_a(N) R_+
    \,=\, \up{1}{N} R_-^{-1} \up{2}{N} R_+ \  .
\ee
This kind of relations appeared first in \cite{ReSe} to describe
relations in $U_q(\sg)$. One may in fact also go in the other
direction, which means to reconstruct a modular Hopf algebra
$\G $ from an object $N $ satisfying the above exchange relations.
To begin with, one has to choose linear
maps $\pi: \G_a \mapsto {\Bbb C}$ in the dual $\G'_a$ of $\G_a$. When
such linear forms $\pi \in \G'_a$ act on the first tensor factor
of $N \in \G_a \o \G$ they produce elements in $\G$:
$$
          \pi(N)  \equiv (\pi\o id) (N) \ \ \in \ \G \ \
          \mbox{ for all } \ \ \pi \in \G'_a\ \ .
$$
$\pi(N) \in \G$ will be called the $\pi$-component of $N$ or just
{\em component of $N$}.
It has been shown in \cite{AlSc} that  the components of $N$ generate
the algebra $\G$, that is, one can reconstruct the modular
Hopf algebra $\G$ from the object $N$. A more precise formulation
is given by the following lemma.

\begin{lemma} {\em \cite{AlSc}} 
Let $\G_a$ be a finite-dimensional,
semi-simple modular Hopf algebra and $\N$ be the algebra generated
by components of $N \in \G_a \o \N$ subject to the relations
\be{Nprop}  
   \up{1}{N} R_+ \up{2}{N} = R_+ \D_a(N) \ \ ,  
\ee  \label{lemma}
where we use the same notations as above. Then $N$ can be
decomposed into a product of elements $N_\pm \in \G_a \o \N$, \
$N=N_+ N_-^{-1}$, such that
\ba
 \D(N)\, \equiv \, N_+'\,  N''\,  (N_-')^{-1} & &
 \ \in \ \G_a \o \N \o \N   \nn \\ [2mm]
 \e(N_\pm)\, \equiv \, e\  \in \ \G_a \ \ \ \  , \ \ \  \  \S(N_\pm)
 &\equiv&  N_\pm^{-1} \ \in \ \G_a \o \N \ \ \ \ ,\ \ \ \ N_\pm^*
 \,\equiv \, N_\mp \nn
\ea
define a Hopf-algebra structure on $\N$. Here, the action of
$\D,\e,\S$ on the second tensor component of $N,N_\pm$ is
understood.  $N_+'$, $N_-'$ and $N''$ are regarded as elements of
$\G_a\o\G\o\G$ with trivial entry in the third and in the second
tensor factors, respectively.
As a Hopf algebra, {\em $\N$ is isomorphic to $\G_a$}.
\end{lemma}

There is another object, similar to $N$, that is equally natural
to consider and that will appear later in the text,
$$
 \wt{N} \,:=\, N_+^{-1} N_-  \ \ \in
   \ \G_a \o \G\ .
$$
Its properties are derived in complete analogy with our treatment
of $N$,
\be{N'}
  R_- \, \D_a(\wt{N})  \,=\, \up{1}{\wt{N}} R_- \up{2}{\wt{N}}
  \ , \ \ \ \  R_+^{-1} \up{1}{\wt{N}} R_- \up{2}{\wt{N}} \,=\,
  \up{2}{\wt{N}} R_+^{-1} \up{1}{\wt{N}} R_- \ .
\ee
An appropriate version of Lemma 1 establishes an isomorphism between
the algebra $\wt{\cal N}$ generated by components of $\wt{N}$ and the
algebra $\G_{op}$. The latter stands for the quantum algebra $\G$ with
opposite multiplication, i.e., elements $\xi,\eta \in \G_{op}$ are
multiplied according to $\xi \cdot \eta : = \eta \xi$.

Observe that property (\ref{stop}) implies that $N^*=\wt{N}$.
We wish to rewrite this simple formula for the action of $*$ on $N$
in a more sophisticated way which proves to be useful in the sequel.
For this purpose, let us introduce an element $S\in\G_a \o \G$ as
follows
\be{SS}
  S \,:=\, N_+ \, \D(\k)\,(\k\o \k)^{-1} \,=\,
 N_- \, \D(\k^{-1})\,(\k\o \k) \  ,
\ee
where $\k$ is some central square root of the ribbon element $v\in\G$,
i.e., $\k^2=v$ and $\k$ commutes with all $\xi\in\G$. The two expressions
for $S$ given in (\ref{SS}) are equivalent due to (\ref{rib}).
It is easy to check that
\be{propS}
  S^* \,=\, S \ , \ \ \ \ S' \,\equiv\, PSP \ = \ S^{-1}\ \ .
\ee
Now we are able to rewrite the $*$-operation on $N$ and $\wt{N}$ with 
the help of $S$:
\be{*N}
  N^* \,=\, S^{-1} \, N^{-1} \, S \ , \ \ \ \
 \wt{N}\^{{}\,\,*} \,=\,  S \, \wt{N}\^{{}\,\,-1} \, S^{-1} \ .
\ee

\noindent
{\sc Example:} ({\em The universal elements for $\Z_q$})
The notion of universal elements can be illustrated with the example
of $\Z_q$. The elements $N_\pm\in\G_a\o\G$ are constructed from
the $R$-matrix (\ref{RZ}):
$N_\pm=\sum_{t,s} q^{\pm ts}\,P^t\o P^s=\sum_s P^s\o h^{\pm s}$
and hence $N=\sum_s P^s \o h^{2s}$. The functorial properties
(\ref{Npmeq}),(\ref{Nprop}) can be verified by using the obvious
identity $\D(P^s)=\sum_k\,P^k\o P^{s-k}$. In order to make these
properties more transparent, we introduce an Hermitian operator
$\wh{p}$ such that $h=q^{\wh p}$. It follows from the definition
of $\t^s$ that $\tau^s(\wh{p})=s$ and that the co-product, antipode
and co-unit act on $\wh{p}$ according to
$$ \D(\wh{p})\ =\ \wh{p}\o e\ +\ e\o\wh{p} \ \ \ , \ \ \
\S(\wh{p})\ =\ -\wh{p} \ \ \ , \ \ \ \e(\wh{p})\ =\ 0\ \ . $$
In these notations, the characteristic projector (\ref{zP}) acquires
the form $P^s=\frac 1p \sum_m q^{m\,(\wh{p}-s)}$ and the universal
elements $N_\pm, N, \wt{N}$ are given by
\be{Nz}
 N_\pm\,=\,q^{\pm\,\wh{p}\,\o\,\wh{p}} \ , \ \ \ \
 N\,=\,q^{2\,\wh{p}\,\o\,\wh{p}}\ , \ \ \ \
 \wt{N}\,=\,q^{-2\,\wh{p}\,\o\,\wh{p}}\ .
\ee
These expressions simplify the task of checking the functoriality
relations in (\ref{Npmeq}),
$$
 \D_a(N_\pm)\,=\,q^{\pm(\D\o id)\,\wh{p}\,\o\,\wh{p}} \,=\,
 q^{\pm\,(\wh{p}\,\o\, e + e\,\o\,\wh{p})\,\o\,\wh{p}} \,=\,
 \up{1}{N}_\pm\,\up{2}{N}_\pm \ .
$$
Observe that the ribbon element $v=\sum_s q^{-s^2}P^s$ can be
written as $v=q^{-\wh{p}\^{{}\;\;2}}$ and hence we may choose
$\k=q^{-\frac 12 \wh{p}\^{\;\;2}}$. A simple calculation gives
$S = e \o e$ for the element $S$ defined in (\ref{SS}). Thus,
formulae (\ref{*N}) simplify for $\Z_q$ and become $N^*=N^{-1}$
and $\wt{N}\^{{}\,\,*}=\wt{N}\^{{}\,\,-1}$.

\subsection{Vertex operators and their structure data.}
Our next aim is to recall the theory of tensor operators for a semi-simple
modular Hopf algebra $\G$. To this end, we combine the carrier spaces
$V^I$ of its finite dimensional irreducible $*$-representations $\t^I$ into
the {\em model space} $\cM = \oplus_I V^I$. Each subspace $V^I\subset\cM$
appears with multiplicity one. The model space $\cM$ comes
equipped with a canonical action of our modular Hopf algebra so that
we can think of $\G$ as being contained in the associative algebra
$\V = \End(\cM)$ of endomorphisms on $\cM$.
Let us also introduce $\C \subset \V$ to denote the center of $\G\subset
\V$ and $\re$ for the unit element of $\V$.

\begin{defn} {\em (Vertex operator)} \label{vertex}
An invertible element $\Phi
\in \G_a \o \V$ is called a {\em vertex operator} for $\G$, if
\begin{enumerate}
\item $\Phi$ intertwines the action of $\G$ on the model space
      $\cM$ in the sense that
      \be{cov} \xi\  \Phi = \Phi \ \D'(\xi) \ \ \mbox{ for all } \
            \xi \in  \G \ \ . \ee
      Here $\xi = e \o  \xi$ on the l.h.s. and $\D'(\xi)= P \D (\xi )P$
      on the r.h.s. are both regarded as elements in $\G_a \o \V$.
\item $\Phi$ obeys the following generalized unitarity relation
      \be{star} \Phi^* =  S^{-1} \Phi^{-1} = \k_a \rk N_+^{-1} \Phi^{-1}
      \rk^{-1}\ \ , \ee
      where $S\in\G_a\o\G$ was defined in (\ref{SS})-(\ref{propS}).
      On the r.h.s., $\rk^{\pm 1} = (e \o \rk^{\pm 1})$ and
      $\k_a = (\k \o {\re})$, so that all these factors are elements
      of $\G_a \o \V$.
\end{enumerate}
Invertibility of $\Phi$ means that there exists an element $\Phi^{-1}
\in \G_a \o \V$ such that $\Phi\ \Phi^{-1} = e \o {\re} = \Phi^{-1} \Phi$.
\end{defn}

Since Definition \ref{vertex} is fundamental to what follows below, let
us discuss it in more detail. In eqs.~(\ref{cov})-(\ref{star})
it would be possible to replace $\D'$ by $\D$ and at the same time
$S$ by $S^{-1}$. We shall meet elements $\Phi$ with such properties
later and call them vertex operators as well.

The relation (\ref{cov}) describes the {\it covariance property} of
$\Phi$. It means that $\Phi$ is a {\em universal tensor operator} for
$\G$ (see, e.g., \cite{MaSc}). More precisely, we may evaluate the element
$\Phi \in \G_a \o \V$ with representations $\t^I$ of $\G$ to
obtain matrices $\Phi^I = (\t^I \o id) (\Phi) \in \End(V^I) \o \V$.
The rows of these matrices form tensor operators which transform
covariantly according to the representation $\t^I$ of $\G$. The
relation (\ref{cov}) may be rewritten in the $R$-matrix formalism
of Subsection 2.2 (see \cite{AFFS}, where a similar calculation was
discussed) as follows:
\be{Ncov}
   \up{1}{N}_\pm \Ph{2}{} \, = \, \Ph{2}{}\, R_\pm \up{1}{N}_\pm  \ \
   \mbox{ or } \ \
   \up{1}{N} \Ph{2}{}\  R_-\, =\,  \Ph{2}{}\  R_+ \up{1}{N} \ \ .
\ee
These relations are equivalent \cite{By1} to the definition of
deformed tensor operators in terms of generalized adjoint actions
of $\G$ which is often used in the theory of ($q$-deformed) tensor 
operators (see, e.g., \cite{Bie}).

Our formula (\ref{star}) for the $*$-operation on $\Phi$ certainly
deserves a more detailed explana\-tion.\footnote{$*$-operations
of a similar form have appeared in \cite{MaSc,AF2,AGS}.}
Both expressions we have provided describe
$\Phi^*$ in terms of $\Phi^{-1}$. Using the intertwining relation
(\ref{cov}) one concludes that the conjugated vertex operator
obeys a transformation law which differs from the covariance
properties of the inverse $\Phi^{-1}$:
$$    \Phi^{-1}  \xi \,=\, \D'(\xi) \, \Phi \ \ \mbox{ while } \ \ \
      \Phi^* \, \xi \,=\, \D(\xi) \,  \Phi^*\ \ . $$
The second relation follows from our assumption (\ref{stop}) on
the behaviour of the co-product under conjugation. Comparison
of the two transformation laws motivates to multiply $\Phi^{-1}$
with a factor $N_+^{-1}$ so that we obtain two objects with
identical covariance properties, namely $\Phi^*$ and $N_+^{-1}
\Phi^{-1}$. In addition, the operation $*$ is supposed to be an
involution, i.e., $(\Phi^*)^* = \Phi$. This requires to dress the operator
$N_+^{-1}\Phi^{-1}$ with factors of $\k$ as we did in the second
expression for $\Phi^*$ in (\ref{star}). All these factors can be moved
to the left of $\Phi^{-1}$ with the help of eq.~(\ref{cov}), so that
$\Phi^* = S^{-1} \Phi^{-1}$. The identity $(\Phi^*)^* = \Phi$
holds then as a consequence of (\ref{propS}).

Suppose for the moment that we are given a vertex operator $\Phi$
in the sense of our Definition \ref{vertex}. Then we can use it
to construct the following {\em structure data} of $\Phi$,
\ba
         F & :=&  \Ph{2}{} \ \Ph{1}{} \ \D_a(\Phi^{-1})
         \  \in \ \G_a \o \G_a \o \V \label{str1}\ , \\ [1mm]
       \s(\rf)& :=&  \Phi \ (e \o\rf) \ \Phi^{-1}\ \ \mbox{ for all }
             \ \   \rf \in \C \subset \V\ , \label{str2}\\[1mm]
         D & := & \Phi \ N \ \Phi^{-1} \in \G_a \o \V \label{str3} \ .
\ea
As they are defined, the last tensor components of $F,D$ and
$\sigma(\rf)$ belong to the algebra $\V$. However, with the help
of relation (\ref{cov}) and standard axioms of Hopf algebra
it is easy to see that $F,D$
and $\s(\rf)$ commute with all elements
$\xi \in \G \subset \V$ and hence that $F \in \G_a \o \G_a \o \C$
while $\sigma(\rf), D \in \G_a \o \C$. Before we give a
comprehensive list of properties of the structure data, we
introduce some more notations,
\ba
        \RR_\pm & \equiv & F'\, R_\pm \ F^{-1} \ \in \
        \G_a \o \G_a \o \C
        \ \ \mbox{ and } \label{RRelem} \\[1mm]
        \D_F(\xi) & \equiv &  F \ (\D(\xi) \o \re) \ F^{-1}
        \ \in \ \G_a \o \G_a \o \C \ \ \label{DelF} .
\ea
Here $F'=(P\o\re)F(P\o\re)$. As a consequence of 
eqs.~(\ref{str1})-(\ref{str2})
and our definition (\ref{RRelem}) we obtain the following
exchange relations for vertex operators:
\be{PhPh}
 \RR_\pm \, \up{2}{\Phi} \, \up{1}{\Phi} \ = \ \up{1}{\Phi} \,
 \up{2}{\Phi} \, R_\pm  \ , \ \ \
 \RR_\pm \, \up{2}{\s} \, \up{1}{\s} (\rf) \,=\, \up{1}{\s} \,
 \up{2}{\s} (\rf) \, \RR_\pm \ . \label{ssig}
\ee
It is also worth noticing that one may think of $\RR$ and $\D_F$
as being obtained from $R$ and $\D$ through a twist with $F$ in
the sense of Drinfeld \cite{Dr2}.

\begin{prop}\label{strucprop}  {\em (Properties of the structure data)}
Let the structure data be defined as in eqs.~(\ref{str1})-(\ref{str3}).
Then it follows from Definition \ref{vertex} that
\begin{enumerate}
\item the element $D \in \G_a \o \C$ may be expressed in terms
      of $\s$ and the ribbon element $v$ so that
      \be{Dop} 
          D = v_a \rv^{-1} \s({\rv}) \ .  
      \ee
      Here $v_a  = (v \o {\re}) \in \G_a \o \C$ and
      $\rv=(e \o {\rv}) \in \G_a \o \C$, that is, we denote
      the ribbon element by $v_a$ and ${\rv}$ when it is regarded as
      an element of $\G_a$ or $\C$, respectively.
\item the elements $F,\R_\pm \in \G_a \o \G_a \o \C$ and $D\in
      \G_a \o \C$ together with the homomorphism $\s:\C \rightarrow
      \G_a \o \C$ obey the following set of relations:
   \ba \label{ax1}
       (e \o F)\; \Bigl( (id \o \D_a )(F)\Bigr) \,  &= &
       \up{3}{\s}(F)\;  \,\Bigl( (\D_a \o id)(F)\Bigr) \ \   , \\[1mm]
         \label{ax3}
      \up{1}{D}\,{\RR}_- \, =\, \RR_+ \,\up{2}{\s}(D)\ \
       & \mbox{ , } & \ \
       \RR_- \up{2}{D}\, =\, \up{1}{\s}(D)\,\RR_+ \ \ ,\\[1mm]
      \label{ax2}
      \up{2}{\s}\,\up{1}{\s} (\rf)  =  \D_F (\s(\rf))
      \ \ &  & \mbox{ for all }  \ \ \rf \in \C \ \ ,  \\[1mm]
      \RR_{\pm,12}\ \up{2}{\s}(\RR_{\pm,13}) \ \RR_{\pm,23}  & =  &
       \up{1}{\s}(\RR_{\pm,23})\ \RR_{\pm,13} \
        \up{3}{\s}(\RR_{\pm,12}) \ \ .   \label{qYB}
    \ea
      The symbol \hbox{$\up{2}{\s}\!(D)$} denotes $(id\o\s)(D)\in\G_a\o
      \G_a\o\C$ and \hbox{$\up{1}{\s}\!(D) = (P\o\re) \up{2}{\s}\!(D)
      (P \o \re)$} with $P$ being the permutation. Similar
      conventions apply to eqs.~(\ref{ax1}),(\ref{qYB}).
 \item the behaviour of the structure data with respect to the
      $*$-operation is given by
     \ba    \label{*F}
     F^* = S_a F^{-1}    \ \ & \mbox{ with }&  \ \
     S_a = (R_+ \D(\k)\,(\k\o\k)^{-1}) \o \re \  \in
    \G_a  \o \G_a \o \C \ , \nn \\[1mm]
      \label{*RR}  \label{*D}
    \RR_\pm^* \, = \,\RR^{-1}_\pm  \ \ & , & \ \ D^*\;=\;D^{-1}\ ,\\[1mm]
    \label{*DF} \label{*s}  \s(\rf )^* = \s(\rf^*)
    \  &, & \ \ \left(\D_F(\xi)\right)^*\;=\; \D_F(\xi^*) \ ,\nn
     \ea
     for all $\xi \in \G$ and $\rf \in \C$. It means, in particular,
     that $D, \RR_\pm$ are unitary while $\s, \D_F$ act as 
     $*$-homomorphisms.
\end{enumerate}
\end{prop}

A proof of the main statements can be found in Appendix A.1. It should
be mentioned that some of the relations given in Proposition
\ref{strucprop} have appeared in the literature before. The
equation (\ref{qYB}) is probably the most characteristic in our
list as it generalizes the usual Yang-Baxter equation. It
appeared first in connection with the quantum Liouville model
\cite{GeNe}; later some universal solution for equation (\ref{qYB})
in the case of $\G=U_q(sl(2))$ has been found \cite{Bab}.
More recently in \cite{BBB}, the elements $F$ and $\RR$ and their
relations were reinterpreted
in the language of quasi-Hopf algebras \cite{Dr2}. As we remarked
already, $F$ may be regarded as a twist and it follows from eq.
(\ref{ax1}) that the twisted co-product $\D_F$ is quasi-coassociative
with co-associator \hbox{$\phi =\, \up{3}{\s}\!(F_{12})\,F_{12}^{-1}$}. 
The latter can be used to rewrite eq.~(\ref{qYB}) as a quasi Yang-Baxter
equation (more details are discussed, e.g., in \cite{BBB,By2}).

Relations (\ref{str3}), (\ref{ax3}) and the first equation in (\ref{PhPh})
have been introduced in \cite{AF1} in a description of deformed
cotangent bundles $T_q^* G$. There, an object $N$ was defined
in terms of $\Phi$ and $D$ through eq.~(\ref{str3}). The relation
(\ref{ax3}) allowed to derive exchange relations for $N$ which
guaranteed that coordinate functions for the fibers of $T_q^* G$
could be obtained from $N$.

We shall see later that the equations in Proposition \ref{strucprop}
have a number of important implications for the lattice theories. 
Reversing this logic, many of the relations in Proposition \ref{strucprop} 
were conjectured as natural properties of a coordinate dependent braiding 
matrix in the continuum WZNW-model \cite{Fa1,CGO,CG1}.

\subsection{Gauge transformations of vertex operators.}

There exists a large gauge freedom in the choice of vertex
operators $\Phi$. In fact, one may replace $\Phi \mapsto \Lambda
\Phi$ with $\Lambda \in \G_a \o \C$ being invertible and unitary.
This transformation does not change the general properties
(\ref{cov})-(\ref{star}) of vertex operators but certainly effects
their structure data. Namely, after the action of $\Lambda$ on
$\Phi$ the initial structure data transform into following ones:
\ba
  F \ \mapsto \ \up{2}{\Lambda} \; \up{2}{\s}\!(\Lambda) \, F \,
 \D_a(\Lambda^{-1}) \ \ \ \ &,& \ \ \ \ D \ \mapsto \ \Lambda\,D\,
 \Lambda^{-1}\ , \nn\\ [1mm]
 \s(\rf) \ \mapsto \ \Lambda \, \s(\rf) \, \Lambda^{-1} &&
 {\rm for\ all}\ \rf\in\C \ , \nn
\ea
where \hbox{$\up{2}{\s}(\Lambda)=(id\o\s)(\Lambda)\in\G_a\o\G_a\o\C$},
as before. One may reduce such a gauge freedom by additional requirements
on the structure data or on the vertex operators. For instance,
the gauge freedom allows to normalize the vertex operators in the
following sense. Consider the element $\rw := \e_a(\Phi)\equiv(\e\o id)
\Phi\in\V$, where $\e:\G\mapsto {\Bbb C}$ stands for the co-unit of $\G$.
An application of the Hopf algebra axiom $(\e\o id)\D\, = id\ $ to
(\ref{cov}) furnishes the identity $\xi \rw = \rw \xi$ and hence
$ \rw \in \C$. {}From this and eqs.~(\ref{str2})-(\ref{str3}) we
conclude that
$$
 (\e\o id)D \,=\, \re \ \ \ \ , \ \ \ (\e\o id)\s(\rf) \,=\, \rf \ \
 {\rm for\ all}\ \rf\in\C \ .
$$
Moreover, (\ref{star}) implies unitarity of $\rw$ (observe that $(\e
\o id)(S) = e$). Therefore, we can perform the gauge transformation
$\Phi \mapsto (e\o\rw^{-1}) \Phi$, which does not change $\s$ and $D$
but normalizes $F$ and $\Phi$ so that, without loss of generality, we
may assume
$$
 \e_a(\Phi) \,=\, \re \in\C\ , \ \ \ \  (id\o\e\o id)F \,=\,
 (\e\o id\o id)F  \,=\, e\o\re\in\G_a\o\C \ .
$$
The normalization of $F$ follows from the normalization and
operator product expansion of $\Phi$ with the help of
$(\e \o id)\D \,=\, id\, =\, (id \o \e)\D$. It also
leads to the identities $(\e\o id)\RR_\pm= (id\o\e)\RR_\pm=
e\o\re$.

{}Finally, let us notice that multiplication of vertex operators $\Phi$ by
element $\sF \in \G_a \o \G$ from the right, i.e., $\Phi \mapsto \Phi
\sF$, corresponds to twisting the co-product of $\G$.
\footnote{The object $\sF$ should not be confused with
our $F \in \G_a \o \G_a \o \C$. We use similar letters
mainly for historical reasons.} Transformations of this kind
relate vertex operators $\Phi_q = \Phi_1 \sF_q$ for the deformed
universal enveloping algebras $U_q(\sg)$ with unitary vertex
operators $\Phi_1$ of the undeformed algebras $U(\sg)$ \cite{Dr2}.

\subsection{On the construction of vertex operators.}
So far, we have considered the vertex operators as given objects.
In the spirit of Lemma \ref{lemma}, however, we can reverse our approach
and think of them as being defined through eqs.~(\ref{str1})-%
(\ref{str3}) with an appropriate set of structure data.
This is made more precise in the following proposition.

\begin{prop} {\em (Reconstruction of $\Phi$ from structure data)}
\label{Rec} For a modular Hopf algebra $\G \cong \G_a$ with
center $\C$, let $F \in \G_a \o \G_a \o \C$ and a homomorphism
$\s : \C \rar \G_a \o \C$ be given. Define the elements
$D \in \G_a \o \C$, $\RR_\pm \in \G_a \o \G_a \o \C$ through
equations (\ref{Dop}),(\ref{RRelem}), respectively, and suppose
that $F, \s$ (together with $D,\RR_\pm$) satisfy the relations
(\ref{ax1})-(\ref{*D}). Then there exists a vertex operator
$\Phi \in \G_a \o \V$ for $\G$ such that
\be{Phdef}     \Ph{2}{}\,  \Ph{1}{} \, = \,  F \, \D_a(\Phi)
       \ \ \ \  ,\ \   \ \ \Phi\ \rf \, =  \, \s(\rf )\  \Phi \ \ .
\ee
In particular, the invertible element $\Phi \in \G_a \o \V$
has the  properties (\ref{cov})-(\ref{star}) and the algebra $\V$
generated by its components is associative. $\V$ may be identified 
with the algebra of operators on the model space $\cM = \bigoplus_I
V^I $, as before.
\end{prop}

{\sc Proof:} Let us only sketch the proof since it is based on the
same computations that are involved in the proof of Proposition
\ref{strucprop}. The construction
of $\Phi$ starts from eqs.~(\ref{Phdef}). In fact, one can
use them to build an abstract algebra $\wt{\V}$ which is
generated by components of an object $\Phi \in \G_a \o \wt{\V}$
and elements in $\C$ such that the two relations (\ref{Phdef})
hold. The properties (\ref{ax1}), (\ref{ax2}) ensure this algebra
to be well defined and associative. Due to eqs.~(\ref{*D}), $\wt{\V}$
admits a consistent $*$-operation which makes $\Phi$
unitary in the sense of eq.~(\ref{star}). In the next step,
an element $N \in \G_a \o \wt{\V}$ is defined by eq.~(\ref{str3}).
With the help
of eqs.~(\ref{ax3}) one proves that $N$ obeys the  relations
(\ref{Nprop}), (\ref{Ncov}) and hence that $\wt{\V}$ contains $\G$
as a subalgebra. This subalgebra is finally used to analyze a concrete
representation of $\wt{\V}$ and to show that $\wt{\V} \cong \V
= \End(\cM)$; hence, components of $\Phi$ become operators on
the model space $\cM$.

Let us apply Proposition \ref{Rec} to the example of $\G \cong U_q(\sg)$.
To this end we need to define appropriate candidates for $F$ and
$\s$ which is achieved with the help of the Clebsch-Gordan maps
$C[TL|S]: V^T \o V^L \rar V^S$ and the $6j$-symbols ${\footnotesize
\SJS{.}{.}{.}{.}{.}{.}}$ of $U_q(\sg)$. Within the space $V^L$ of
highest weight $L$, we fix a basis of eigenvectors $e^L_\lambda$
for the Cartan subalgebra with eigenvalues $\lambda$ and denote
the associated Clebsch-Gordan coefficients by
\raisebox{-.5ex}{$\CG{T}{L}{S} {\vth}{\lambda}{\varsigma}$}.
Now define $F,\s $ such that
$$ F^{TL} = (\t^T \o \t^L)(F)\ \ \  \mbox{ and }\ \ \  \s^L(\wh\sp)
   = (\t^L \o id)(\s(\wh\sp)) \ \ \ \mbox{ have matrix elements } $$
\ba
     F^{TL}_{\vth\lambda,\vth'\lambda'}  & = &    \sum_{S,\vs}
      \SJS{T}{L}{S}{\hat\sp}{\hat\sp+\vth+\lambda}
     {\hat\sp+\lambda}^*\      \CG{T}{L}{S}{\vth'}{\lambda'}{\vs} \ ,
     \label{Fqdef}  \\[1mm]
     \s^L(\wh\sp)_{\lambda,\lambda'} & =&  (\wh\sp +\lambda)
     \ \dl_{\lambda,\lambda'} \ . \label{sigdef}
\ea
Here ${\wh\sp}$ is a {\it rank\/}$(\sg)$-dimensional vector of
elements in $\C$ with $\t^K(\wh\sp) = K$. Other notations
and conventions are explained in Appendix A.2.

\begin{prop} {\em (Vertex operators for $U_q(\sg)$)}
\label{Uqvert}
There exist vertex operators $\Phi_q$ for the
deformed universal enveloping algebras $U_q(\sg)$ such that
$$
    \Ph{2}{}_q\,  \Ph{1}{}_q \, = \,  F \, \D_a(\Phi_q)
   \ \ \ \  ,\ \   \ \ \Phi_q\ \rf \, =  \, \s(\rf )\  \Phi_q\ \ .
$$
Here $F$ is built up from the $6j$-symbols and the Clebsch-Gordan
maps of $U_q(\sg)$ as in eq.~(\ref{Fqdef}) and $\s$ is given by
(\ref{sigdef}).
\end{prop}

The statement follows directly from the Proposition \ref{Rec}
once the relations (\ref{ax1})-(\ref{*D}) have been checked to
hold for $F, \s$. The latter is done in Appendix A.2. Let us mention
that formulae similar to (\ref{Fqdef}) were considered in
\cite{BBB,By2}.

\subsection{Vertex operators for $\Z_q$.}

To conclude our discussion of vertex operators, let us
provide an explicit formula for $\Phi$ in our
standard example $\G = \Z_q$. Let us fix a set of normalized
basis vectors $ | s \rangle, s = 0, \dots , p-1,$ for
the one-dimensional
carrier spaces $V^s$ of the representations $\t^s$. They
span the $p$-dimensional model space $\cM = \bigoplus_s V^s$.
On this space one can introduce a unitary operator
$\wh{Q}\in End(\cM)$ by
$$
  \wh{Q} \ |{p-1} \rangle \, = \, |0 \rangle \ \ \ \
  \mbox{ and } \ \ \ \ \wh{Q} \ | s \rangle \, = \, |{s+1}  \rangle
$$
for all $s = 0, \dots, p-2$. This operator obeys Weyl commutation
relations with the generator $h \in \Z_q$, i.e., \hbox{$q\,\wh{Q}\,h=
h\,\wh{Q}$}. With the help of $\wh{Q}$ and the characteristic projectors
$P^s$ introduced in Subsection 2.1 we are able to define $\Phi$:
$$ \Phi \,:=\,  \sum_s \ P^s \o \wh{Q}  ^s \, = \,  {1 \over p}
   \sum_{s,t} q^{-st} h^t \o \wh{Q} ^s \ \ \in \ \G_a \o End(\cM)\ \ .
$$
It follows from the unitarity of $\wh{Q} $ and the Weyl relations
of $\wh{Q}$ and $h$ that $\Phi$ obeys all the defining properties
of a vertex operator (as we explained in Subsection 2.2, the element
$S$ in eq.~(\ref{star}) becomes trivial for $\G=\Z_q$). One may then
compute the structure data. To this end it is convenient to employ
the operator $\wh{p}$ introduced in Subsection 2.2  such that
$h=q^{\wh{p}}$. Since the commutative algebra $\Z_q$ is isomorphic
to its center $\C$, all elements in $\Z_q$ can be regarded as
elements of $\C$ and we use our standard notational conventions
whenever we do so, in particular we shall use ${\sf h}=q^{\wh\sp}$
for ${\sf h}, {\wh\sp} \in \C$. We also introduce an anti-Hermitian
operator $\wh{\vs}$ by $\wh{Q}=e^{\wh{\vs}}$, so that the Weyl
relations for $\wh{Q}$ and $h$ imply $\,[\,\wh\sp\,,\wh\vs\,]=\re$.
Within these notations our basic objects look as follows:
$$
 h\,=\,q^{\wh{p}}\,\in\G\ , \ \ \ \
 \rv=q^{-\wh{\sp}\^{{}\;\; 2}}\,\in\C \ , \ \ \
 \Phi \,=\,e^{\wh{p}\,\o\,\wh\vs}\,\in\G_a\o End(\cM)\ .
$$
Now expressions for the structure data may be obtained by short
computations,
$$ F \,=\, e\o e \o  \re  \ , \ \ \ \
    \s (q^{\wh\sp}) \, = \, h^{-1} \o q^{\wh\sp} \,=\,
   q^{-\wh{p}\,\o\,\re+e\,\o\,\wh\sp}  \ , \ \ \ \
   D \,=\, q^{-2\,(\wh{p}^2 \,\o\, \re -\wh{p}\,\o\,\wh{\sp})}
$$
and \hbox{$\RR_\pm = q^{\wh{p}\,\o\,\wh{p}\,\o\,\re}$}.

Let us remark that, although the example of vertex operators for $\Z_q$
is fairly trivial, it nevertheless shares some features with the case
of $\G = U_q(\sg)$. Indeed, the ribbon element of $U_q(\sg)$ is given by
$\rv = q^{-\wh{\sp}\,(\wh{\sp}+\rho)}$ \cite{Dr2} where $\wh\sp\in\C^{\o_r}$
is a $r = rank(\sg)$-dimensional vector such that $\tau^K(\wh\sp)=K$ and
$\rho$ is the sum of the positive roots. Our above formula (\ref{sigdef})
means that
$$
\s({\wh{\sp}})\,=\, {{\wh H}\o\re+e\o \wh{\sp}} \ , \ \ \
D \, = \, (\chi\, \o\, \re) \cdot \,q^{-2 \wh{H}\o\wh{\sp}}  \ ,
$$
where $\wh H$ is a vector of elements in the Cartan subalgebra
such that $\wh H e_\lambda^L  = \lambda e^L_\lambda$ and the
element $\chi \in \G = U_q(\sg)$ can be worked out easily
with the help of eq.~(\ref{Dop}). Such expressions, or special
cases thereof, may be found in \cite{CrGe,CGO,AF1,CG1,By2}).
The element $F$ and the vertex operators $\Phi$ are certainly
quite non-trivial for $U_q(\sg)$ (for some explicit examples
see \cite{FaGa,CG2,ByFa,By2}).

\section{A TOY MODEL FOR THE DISCRETE WZNW THEORY}

In the rest of this paper we shall apply the theory of modular Hopf
algebras and their vertex operators to construct and investigate
the lattice  WZNW-model. We start with a simple {\em toy model}
for which the lattice consists of only one site and one edge
(see Figure 1). When we discuss the general notion of lattice
current algebras in Section 4, we shall understand that they
contain chiral observables $M$ (the {\em chiral monodromies})
being assigned to the edge.

\hbox{\begin{picture}(300,50)
    \put(200,20){\circle{30}}
    \put(216,20){\circle*{3}}
    \put(171,17){$\large M$}
    \put(220,17){$\large \Phi$}
\end{picture} }
\begin{center}
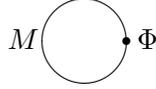
\parbox{13cm}{ \small {\bf Figure 1:} Single-vertex lattice.
Chiral observables $M$ are assigned to the edge while chiral
vertex operators $\Phi$ sit on the vertex.}
\end{center}

\subsection{Properties of chiral vertex operators.}
Later in the text we shall find that the global chiral observable
$M$ in the lattice current algebra obeys the following relation
\be{rM}
\up{2}{M} R_+ \up{1}{M} \ = \ R_- \D_a (M) \ ,
\ee
where $R_\pm$, $\D_a$ are attributes of the modular Hopf algebra $\G$
as before. Components of $M$ generate an algebra $\J$ with center
denoted by $\C$.

Eq.~(\ref{rM}) reminds us of the defining relation (\ref{DN})
for the universal element $N$, which contains all the information
about the structure of $\G$. Indeed, the only difference is that
the $R_+$ on the l.h.s.
of eq.~(\ref{DN}) has been replaced by $R_-$. A short computation
reveals that we can pass from eq.~(\ref{rM}) to (\ref{DN}) by
rescaling $M$ with the ribbon element $v_a =(v \o e) \in \G_a \o \G$.
This implies that $N \mapsto v_a M$ provides an isomorphism of
the algebras $\G$ and $\J$. In particular, the commutation relations
for $M$,
\be{MMr}
  R_\pm^{-1}\, \up{2}{M} R_+ \up{1}{M} \ =\
  \up{1}{M} R_-^{-1} \up{2}{M} \, R_\mp
\ee
coincide with eqs.~(\ref{NN}) for the element $N$. The isomorphism
of $\J$ and $\G$ certainly implies that there is a $*$-operation on
$\J$ given by the formula (\ref{*N}) with $N$ replaced by $M$ (notice
that the factor $v_a$ is unitary). The lattice theories, however,
choose a different conjugation which we discuss in Subsection 3.3
below.

Now let us introduce a vertex operator $\Phi$ for $\J \cong \G$.
It will be called {\em chiral vertex operator} of the toy model
and its properties can be copied from the relations (\ref{cov})
-(\ref{str3}) when we keep in mind to replace $N$ by $v_a M$,
\ba
       \eta\ \Phi \, = \, \Phi \ \D'(\eta) \ \  & , & \  \
       \up{1}{M} \Ph{2}{}\  R_-\, =\,  \Ph{2}{}\  R_+ \up{1}{M} \ ,
       \label{rcov}  \\[1mm]
       \Ph{2}{}\,  \Ph{1}{} \, = \,  F \, \D_a(\Phi) \ \ & , & \ \
       \R_\pm \Ph{2}{}\  \Ph{1}{} \, = \, \Ph{1}{}\ \Ph{2}{}\ R_\pm\ ,
       \label{rOPE}\\[1mm]
       D\,\Phi \, = \, v_a \Phi \ M  \ \ &,& \ \ \label{rmon}
       \Phi\, \rf  =  \, \s(\rf )\,  \Phi  \ \
       {\mbox{for  all}}\ \rf\in\C \ . \label{rs}
\ea
Here $\eta\in\J$, $\C$ stands for the center of $\J$,
and we used the same notations as in the previous section. 
The components of $\Phi\in\G_a\o\V$ give rise to the 
{\em algebra $\V$ of chiral vertex operators}. 
Together with components of $M$, they act on the model 
space  $\cM \,=\,\bigoplus_I  V^I$.

We refer to the first equation in (\ref{rOPE}) as
{\em operator product expansions} (OPE) for $\Phi$ and call $F$ the
{\em universal fusion matrix}. The second formula in (\ref{rOPE})
follows from the operator product expansions; it describes {\em braid
relations} for the chiral vertex operators and hence leads to
interpret $\RR_\pm$ as the {\em braiding matrix} of our model.~\footnote{
This will become clearer in the full lattice theory where braid
relations of vertex operators assigned to different sites contain
only $\RR_\pm$ and the factor $R_\pm$ is absent. Observe also that in
the quantum non-deformed limit, i.e., $\ga\rar 0$, $\hbar\neq 0$,
$q=e^{i\hbar\ga}\rar 1$, the $R$-matrix $R_\pm$ approaches $e\o e$
whereas the limit of $\RR_\pm$ is non-trivial (cf. also
\cite{AF1,ByFa}).}

There exists a nice pictorial presentation for the described 
algebraic structure. Definitions for the basic objects -- 
except from $D, M$ -- are given in Figure 2. Pictures for $M$
and $D$ may, in principle,  be constructed with the help
of eq.~(\ref{Dop}), eqs.~(\ref{rmon})
and an appropriate presentation of the ribbon element.
{}From the basic blocks we can built up the equations
(\ref{rcov})-(\ref{rs}) as in Figure 3.
All these pictures are separated by a thick solid line into
left and  right halves with dotted lines appearing on the left
side while thin solid lines exist only on the right side. Our
graphical rules are the same as in \cite{KR1}, and, in their
terminology, the dotted lines may be said to live in the
{\em shadow world}.

\begin{figure}
\setlength{\unitlength}{.5pt}
\begin{picture}(380,120)(0,0)
\put(50,70){\line(1,-1){25}}
\multiput(49,71)(-2.8,2.8){9}{\circle*{1}}
\put(100,70){$\:\Phi$}
\put(240,70){\line(1,-1){45}}
\put(262,48){\line(-1,-1){22}}
\multiput(240,70)(-2.8,2.8){16}{\circle*{1}}
\put(312,70){$\:\Delta_a(\Phi)$}
\thicklines
\multiput(25,45)(0.25,0){6}{\line(1,1){50}}
\multiput(195,25)(0.25,0){6}{\line(1,1){90}}
\end{picture} \hspace*{.5cm}
\begin{picture}(400,100)(0,0)
\multiput(62,85)(-4,0){12}{\circle*{1}}
\multiput(39,60)(-2.8,2.8){10}{\circle*{1}}
\put(37,69){\footnotesize$\rf$}
\put(92,70){$\:{\rf}$}
\multiput(266,93)(-4,0){17}{\circle*{1}}
\multiput(247,72)(-4,0){7}{\circle*{1}}
\multiput(232,59)(-2.8,2.8){18}{\circle*{1}}
\put(231,77.5){\footnotesize$\rf$}
\put(300,70){ $\:\sigma({\rf})$}
\thicklines
\multiput(23,45)(0.25,0){6}{\line(1,1){54}}
\multiput(215,40)(0.25,0){6}{\line(1,1){70}}
\end{picture}\\
\begin{picture}(380,140)(0,0)
\thicklines
\multiput(75,45)(-2.8,2.8){18}{\circle*{1}}
\thinlines
\multiput(25,45)(2.8,2.8){8}{\circle*{1}}
\multiput(54,74)(2.8,2.8){8}{\circle*{1}}
\put(80,70){$\:\RR_+$}
\thicklines
\multiput(225,45)(-2.8,2.8){8}{\circle*{1}}
\multiput(196,74)(-2.8,2.8){8}{\circle*{1}}
\thinlines
\multiput(175,45)(2.8,2.8){18}{\circle*{1}}
\put(230,70){$\:\RR_-$}
\multiput(355,45)(-2.8,2.8){18}{\circle*{1}}
\multiput(333,67)(2.8,2.8){10}{\circle*{1}}
\put(362,70){$\:F$}
\end{picture}\hspace*{.5cm}
\begin{picture}(400,140)(0,0)
\thicklines
\multiput(25,45)(0.25,0){6}{\line(1,1){50}}
\thinlines
\put(40,60){\line(1,0){42}}
\put(61,81){\line(1,-1){21}}
\put(100,65){ $\:\eta$}
\put(58,65){\footnotesize$\eta$}
\thicklines
\multiput(205,45)(0.25,0){6}{\line(1,1){50}}
\thinlines
\put(245,85){\line(1,-1){40}}
\put(233,73){\line(1,0){24}}
\put(216,56){\line(1,0){58}}
\put(290,63){ $\:\Delta'(\eta)$}
\put(243,61){\footnotesize$\eta$}
\end{picture}
\begin{center}
\parbox{13cm}{ \small {\bf Figure 2:}
Graphical presentation of our basic objects. Pictures for $D$ and
$M$ exist as well, but they are more complicated (cf. remarks
in the text).}
\end{center}
\setlength{\unitlength}{1pt}
\end{figure}

\begin{figure}
\setlength{\unitlength}{.5pt}
\begin{picture}(400,200)(0,0)
\put(70,70){\line(1,-1){45}}
\multiput(70,70)(-2.8,2.8){16}{\circle*{1}}
\put (70,25){\line(-1,1){23}}
\multiput (47,47)(-2.8,2.8){8}{\circle*{1}}
\put(165,70){$\stackrel{(\ref{rOPE})}{=}$}
\thinlines
\put(280,70){\line(1,-1){45}}
\put(302,48){\line(-1,-1){22}}
\multiput(280,70)(-2.8,2.8){16}{\circle*{1}}
\multiput(258,92)(2.8,2.8){8}{\circle*{1}}
\thicklines
\multiput(235,25)(0.25,0){6}{\line(1,1){90}}
\multiput(25,25)(0.25,0){6}{\line(1,1){90}}
\end{picture}
\begin{picture}(400,240)(0,-20)
\put(50,70){\line(1,1){25}}
\multiput(50,70)(-2.8,-2.8){9}{\circle*{1}}
\multiput(50,20)(-2.8,2.8){10}{\circle*{1}} 
\put(50,20){\line(1,-1){25}}
\put (150,95) {\line(0,-1){100}}
\multiput(240,95)(0,-4){25}{\circle*{1}}
\put(350,20){\line(1,1){25}}
\multiput(350,20)(-2.8,-2.8){9}{\circle*{1}}
\multiput(350,70)(-2.8,2.8){9}{\circle*{1}}
\put(350,70){\line(1,-1){25}}
\put(95,45){=}   \put(35,-50){$\Phi^{-1} \Phi = e \o \re $}
\put(290,45){=}  \put(240,-50){$e \o \re = \Phi\  \Phi^{-1}$}
\thicklines
\multiput(125,95)(0.25,0){6}{\line(0,-1){100}}
\multiput(265,95)(0.25,0){6}{\line(0,-1){100}}
\multiput(350,70)(0.25,0){6}{\line(1,1){25}}
\multiput(350,70)(0.25,0){6}{\line(-1,-1){25}}
\multiput(350,20)(0.25,0){6}{\line(-1,1){25}}
\multiput(350,20)(0.25,0){6}{\line(1,-1){25}}
\multiput(50,20)(0.25,0){6}{\line(1,1){25}}
\multiput(50,20)(0.25,0){6}{\line(-1,-1){25}}
\multiput(50,70)(0.25,0){6}{\line(-1,1){25}}
\multiput(50,70)(0.25,0){6}{\line(1,-1){25}}
\end{picture}\\[-.3cm] \setlength{\unitlength}{.6pt}\hspace*{-.5cm}
\begin{picture}(400,140)(-20,0)
\put(50,70){\line(1,-1){25}}
\multiput(49,71)(-2.8,2.8){9}{\circle*{1}}
\multiput(48,68)(-4,0){9}{\circle*{1}}
\multiput(31.5,51.5)(-2.8,2.8){6}{\circle*{1}}
\put(30,56.5){\footnotesize $\rf$}
\put(250,70){\line(1,-1){25}}
\multiput(249,71)(-2.8,2.8){11}{\circle*{1}}
\multiput(227,93.1)(4,0){12}{\circle*{1}}
\multiput(243.6,79.7)(4,0){5}{\circle*{1}}
\put(248,82.3){\footnotesize $\rf$}
\put(150,68){$\stackrel{(\ref{rs})}{=}$}
\thicklines
\multiput(25,45)(0.25,0){6}{\line(1,1){50}}
\multiput(225,45)(0.25,0){6}{\line(1,1){55}}
\end{picture}\\[-1cm] \hspace*{-4mm}
\begin{picture}(400,140)(-20,0)
\put(50,70){\line(1,-1){30}}
\multiput(49,71)(-2.8,2.8){10}{\circle*{1}}
\thinlines
\put(54,74){\line(1,0){38}}
\put(73,93){\line(1,-1){19}}
\thicklines
\multiput(20,40)(0.25,0){6}{\line(1,1){60}}
\put(70,79){\footnotesize $\eta$}
\multiput(249,71)(-2.8,2.8){10}{\circle*{1}}
\thicklines
\multiput(220,40)(0.25,0){6}{\line(1,1){60}}
\thinlines
\put(250,70){\line(1,-1){30}}
\put(240,60){\line(1,0){20}}
\put(225,45){\line(1,0){50}}
\put(249,49){\footnotesize$\eta$}
\put(150,68){$\stackrel{(\ref{rcov})}{=}$}
\end{picture}\hspace*{-1.0cm}
\begin{picture}(300,140)(-20,0)
\put(53,73){\line(1,0){36}}
\put(71,91){\line(1,-1){18}}
\thicklines
\multiput(22,42)(0.25,0){6}{\line(1,1){56}}
\thinlines
\multiput(48,68)(-4,0){10}{\circle*{1}}
\multiput(30.6,49.4)(-2.8,2.8){7}{\circle*{1}}
\put(27.5,56){\footnotesize$\rf$}
\put(68,77){\footnotesize$\eta$}
\thicklines
\multiput(222,42)(0.25,0){6}{\line(1,1){56}}
\thinlines
\put(247,67){\line(1,-1){18}}
\put(229,49){\line(1,0){36}}
\multiput(252.4,72.4)(-2.8,2.8){7}{\circle*{1}}
\multiput(272,92)(-4,0){11}{\circle*{1}}
\put(250,79.5){\footnotesize$\rf$}
\put(244,54){\footnotesize$\eta$}
\put(160,68){=}
\end{picture}
\begin{center}
\parbox{13cm}{ \small {\bf Figure 3:}
Pictorial presentation of some basic relations. Only the left
equations in (\ref{rcov}), (\ref{rOPE}) and the right equation in
 (\ref{rs}) are depicted.
The figure in the lower right corner means that $\rf \in \C$ is
central in $\G$. More rules are explained in the text.}
\end{center}
\setlength{\unitlength}{1pt}
\end{figure}

\subsection{Second chirality.} What we have discussed so far
will be relevant for right chiral objects in the discrete WZNW model.
Now we have to describe an analogous construction for the left chiral
sector of the theory. To distinguish the two chiralities, we mark the
objects of the previous subsection by an extra index $r$ so that
$M^r = M,\Phi^r = \Phi, F_r = F, \s_r = \s  \dots $ etc. Their left
chiral counterparts will have an index $l$.

To introduce left chiral vertex operators $\Phi^l$ we follow the
same strategy as in the previous subsection. Namely, we postulate
algebraic relations for an object $M^l$ (which will be justified
in Section 4) and use them as the basic input for our left
chiral theory. So let us assume that we are given some object $M^l$
such that
\be{lM}
   \up{1}{M}\^l \ R_- \up{2}{M}\^l\   = \
   R_+ \ \D_a (M^l)\ \ .
\ee
The algebra generated by components of $M^l$ will be denoted by
$\J^l$ and we use the symbol $\C^l$ for its center.

It is easy to see that the properties of $v_a^{-1}M^l$ coincide with
those of the element $\wt{N}$ introduced in Subsection 2.2, eq.
(\ref{N'}). This holds, in particular, for the commutation relations,
\be{MMl}
  R_\pm^{-1}\up{1}{M}\^l \ R_- \up{2}{M}\^l \,  = \,
 \up{2}{M}\^l \ R_+^{-1} \up{1}{M}\^l\ R_\mp\,  \ \ .
\ee
Thus, the algebra $\J^l$ is isomorphic to the algebra generated by
components of $\wt{N}$, i.e., to $\G_{op}$ ( ${}_{op}$ means the
opposite multiplication, cf.~Subsection 2.2).

Since eqs.~(\ref{lM})-(\ref{MMl}) differ from the properties of $M^r$,
the relations for the left chiral vertex operators will differ from
those we had in the right chiral sector. The consistent definition
of the left vertex operators is provided by the following list of
fundamental relations:
\ba
   \eta\ \Phi^l \, = \, \Phi^l \ \D(\eta) \ \ & , & \ \
  \up{1}{\Phi}\^l \, R_- \, \up{2}{M}\^l \,  =  \,
  \up{2}{M}\^l \, \up{1}{\Phi}\^l \, R_+ \label{lcov}\ \ , \\[1mm]
  \up{1}{\Phi}\^l\,\up{2}{\Phi}\^l\,=\, \,F_l \,\D_a(\Phi^l)\ \ &,& \ \
   \RR^l_\pm \, \up{1}{\Phi}\^l\,\up{2}{\Phi}\^l\,  =  \,
  \up{2}{\Phi}\^l\,\up{1}{\Phi}\^l \, R_\pm  \label{lOPE} \ \ ,\\[1mm]
  D_l \ \Phi^l  \, = \,   v_a^{-1}  \Phi^l \  M^l \label{lD} \ \ &,& \ \
    \Phi^l\,\rf \,=\, \s_l(\rf)\,\Phi^l \ \
  \mbox{ for all }\ \rf\in\C^l  \ \ . \label{ls}
\ea
Components of $\Phi^l\in\G_a\o\V^l$ generate the algebra $\V^l$ of left
chiral vertex operators and act on the left model space $\cM^l
\cong \bigoplus_I V^I$. Starting
from the defining equation (\ref{lM}) for $M^l$ one may check that
the exchange relations (\ref{lcov}) describe a consistent
transformation law of the vertex operators $\Phi^l$. It is
then clear that the left vertex operators obey eqs.
(\ref{lOPE})-(\ref{ls}) with some appropriate structure data $F_l,
\s_l, D_l, \R_\pm ^l$. The consistency relations for the
left structure data can be worked out in analogy to our
discussion of Proposition \ref{strucprop}. For more detailed
explanations see Appendix A.3.

Let us now combine the two chiral theories by constructing
their tensor product so that all operators act on the space
$\cM^l \o \cM^r$ with trivial action of the right chiral objects
on the first tensor factor and vice versa. In terms of exchange
relation this corresponds to
\ba
 \up{1}{\Phi}\^r\, \up{2}{\Phi}\^l \,=\,  \up{2}{\Phi}\^l
 \up{1}{\Phi}\^r\  &\ \ \  ,\ \ \  &
 \up{1}{M}\^r\, \up{2}{M}\^l \,=\,  \up{2}{M}\^l
 \up{1}{M}\^r \ ,   \label{comm1} \\ [1mm]
 \up{1}{\Phi}\^r\, \up{2}{M}\^l \,= \,  \up{2}{M}\^l
 \up{1}{\Phi}\^r\  &  ,  &
 \up{1}{\Phi}\^l\, \up{2}{M}\^r \,= \, \up{2}{M}\^r \
 \up{1}{\Phi}\^l\ .  \label{comm2}
\ea
Components of the chiral vertex operators $\Phi^l, \Phi^r$
generate an algebra $\W = \V^l \o \V^r$. Although this combination
of chiral theories appears to be quite trivial, it sets the stage
for the construction of the quantum group valued field $g$ that
we are about to discuss in Subsection 3.4.

Before we get there, let us explain how to incorporate our new
left chiral
objects into the graphical presentation discussed at the end of
the previous subsection. The pictures for the left chiral theory
are simply mirror images of those in Figures 2,3, that is, left chiral
objects have their dotted lines on the right side and thin solid lines
on the left side of the thick solid line. To present the tensor
product of the left- and right theory, we draw all objects
into the same pictures. Now there are dotted and solid lines
on both sides. If we add the rule that these lines of different
style do not interfere, we obtain commutativity of the two
chiralities as expressed in eqs.~(\ref{comm1})-(\ref{comm2}).

\subsection{$*$-operation for chiral vertex operators.}
In principle, a $*$-operation for $M^l, M^r$ and the associated
vertex operators could be introduced along the lines of Section 2.
But as we indicated the lattice models choose a slightly
different conjugation. Its description requires to introduce a
new object.

By definition, the models spaces $\cM^l, \cM^r$ carry an action of
the modular Hopf-algebra $\G$. With the help of the co-product
$\D$ this gives rise to a canonical action of $\G$ on the tensor
product $\cM^l \o \cM^r$ and hence to an embedding $\iota$ of the
quantum algebra $\G$ into the algebra $\W = \V^l \o \V^r$ of
chiral vertex operators. For the exchange relations of $\iota(\xi)$
and chiral vertex operators, our construction implies:
$$  \iota (\xi) \ \Phi^r \ =\  \Phi^r \  \D'_\iota(\xi) \ \ \ ,
    \ \ \ \iota (\xi) \ \Phi^l \ = \ \Phi^l \ \D_\iota(\xi)
    $$
for all $\xi \in \G$; we used $\D_\iota(\xi) = (id \o \iota)\D(\xi)$
and similarly for $\D'_\iota$.
These relations imply that $\Phi^l, \Phi^r$ transform covariantly
with respect to our new action $\iota$ of $\G$ on $\W$. They can
be rewritten in the $R$-matrix formulation,
$$
 \up{1}{N}_{\pm} \, \up{2}{\Phi}\^r \,=\,
 \up{2}{\Phi}\^r \, R_\pm \, \up{1}{N}_{\pm} \ \ , \ \
 \up{1}{N}_{\pm} \, \up{2}{\Phi}\^l \,=\,
 \up{2}{\Phi}\^l \, \up{1}{N}_{\pm} \, R_\pm \ ,
$$
where $N_\pm = (id \o \iota)(R_\pm) \in \G_a \o \W$. In our
pictorial presentation the objects $N_\pm$ would appear as
over-/under- crossings of thin and thick solid lines. Hence,
they have thin solid lines on both sides of the boundary between
the left and the right world. This corresponds to the fact that
components of $N_\pm$ act nontrivially on both factors in 
$\cM^l\o \cM^r$, that is, they are {\em not chiral}. The same holds 
true for the product $N = N_+ (N_-)^{-1}$.

Now we are prepared to describe the $*$-operation which is relevant
for the toy model. To this end, we build an object $S_\iota$ with
the help of $\iota$ by $S_\iota = (id \o \iota)(S) \in \G_a
\o \W$ and $S \in \G_a \o \G$ is defined as in Subsection 2.2. It
is used to extend the $*$-operation on $\G \cong \iota(\G) \subset
\W$ to the algebra of chiral vertex operators:
$$   (\Phi^r)^* \ = \ S_\iota^{-1} (\Phi^r)^{-1} \ \ \ , \ \ \
     (\Phi^l)^* \ = \ S_\iota \ (\Phi^l)^{-1} \ \ . $$
The first formula looks familiar already and since the exchange
relations of $\iota (\xi)$ with $\Phi^r$ coincide with eq.
(\ref{cov}), consistency need not to be checked again. The second
formula is a variant of eq.~(\ref{star}) which is adapted to
the algebraic properties of the left chiral theory. To prove
that it is consistent one has to modify our discussion in
Subsection 2.3 slightly. We leave this to the reader. It
remains to show that the adjoints of $\Phi^l$ and $\Phi^r$
commute; this is not obvious at all, since $S_\iota$ is not
a chiral object. Commutativity of the adjoints may be seen
most easily if we rewrite the adjoints in the form (\ref{star})
which involves conjugation with $\k$ (which is $\iota(\k)$ in
our case). Then the desired consistency follows from the
transformation law of vertex operators under the action of
$\iota(\xi)$ and the Yang-Baxter equation (see also \cite{AFFS}).

It follows from eqs.~(\ref{rmon}),(\ref{lD}) that the conjugation
acts on the chiral monodromies $M^r, M^l$ according to
$$ (M^r)^* = S_\iota^{-1} (M^r)^{-1} S_\iota \ \ \ , \ \ \
   (M^l)^* = S_\iota (M^l)^{-1} S_\iota^{-1} \ \ .
$$
We shall rediscover such a behaviour for the chiral monodromies
of the lattice theory in Subsection 4.3.

\subsection{Quantum group valued field $g$.}
So far we have reached a good level of understanding for our
right- and left chiral theories which act naturally on
the tensor product $\cM^l \o \cM^r$ of chiral model spaces.
In this subsection we would like to have a closer look at
the {\em diagonal subspace}
$$ \H = \bigoplus_K V^{\bar K} \o V^K \
    \subset \cM^l \o \cM^r\ \ .$$
While components of $M^l, M^r$ leave $\H$ invariant, this is
certainly not the case for the vertex operators $\Phi^l,
\Phi^r$. Nevertheless, the vertex operators can be combined
into a new object $g$ which admits restriction to the diagonal
subspace $\H$.

The construction of $g$ requires a careful preparation.
Let us begin this with some remarks on the center $\C$ of
$\G$ (recall that $\C^r\cong\C^l\cong\C$). First, observe that  $\C$
is spanned by the characteristic projectors $P^J$ of irreducible
representations $\t^J$ of $\G$, i.e., by  projectors $P^J\in \C$
which obey $\t^K(P^J) = \dl_{K,J}$. Notice also that
the antipode $\S$ maps the element $P^K\in \C$ to the characteristic
projector $P^{\bar K} \in \C$ of the conjugate
representation $\t^{\bar K}$, i.e., $\S (P^K) = P^{\bar K}$.\footnote{
Strictly speaking, the conjugate
of $\t^K$ is obtained with the help of a transpose $\ ^t$
as $\ ^t \t^K \circ \S$. The latter is isomorphic to
$\t^{\bar{K}}$ (this property defines the label $\bar{K}$).}

Returning to our toy model, we combine the canonical isomorphism
$\nu: \C^r \rar \C^l$ and action of the antipode $\S$ into
a map $\S_{lr}: \C^r \rar \C^l$, $\S_{lr}(\rf)=\S(\nu(\rf))$.
With the help of this map we can characterize the diagonal
subspace $\H$ as a subspace generated by all vectors
$\phi \in \cM^l \o \cM^r$ such that $\rf \phi = \S_{lr}(\rf) \phi$
holds for all $\rf \in \C^r$. In this language, the restriction
to $\H$ means to impose the constraint $ \rf = \S_{lr}(\rf)$ for
all $\rf \in \C^r$. This constraint couples the two chiralities
and it seems natural to restrict the choice of the left- and right
structure data $F_\a, \s_\a, D_\a, \RR^\a_\pm$ at the
same time. Notice that they were completely independent until
now, as long as they solved the appropriate consistency relations.
So let us agree to adjust the choice of the
structure data for the  left chirality to whatever we use in the
right chiral part such that
\ba
    F_l \,=\, \S^{(2)}_{lr}({F'_r}^{-1}) \hspace*{1.5cm}
    \ \ \ \ & , & \ \ \ \ D_l\, = \, \S^{(1)}_{lr}(D^{-1}_r) \ ,
    \label{lstruct}\\[1mm]
    \RR_\pm^l \, = \, \S^{(2)}_{lr} ({\RR^r_\pm}') \hspace*{1.5cm}
    \ \ \ \ & , &\ \ \ \ \s_l(\rf)\, =\,
    \S^{(1)}_{lr} ( \s_r \circ \S^{-1}_{lr} (\rf))\ , \label{rltrafo}\\[2mm]
   \mbox{ with }\ \  \S^{(n)}_{lr}\, :=\, (\S^{-1} \o \S^{(n-1)}_{lr})
    & &\mbox{ and }\ \ \ \
    \S^{(0)}_{lr} \,:=\, \S_{lr}  \label{Sdef}
\ea
and the prime on $F_r$ and $\RR_\pm^r$ denotes permutation
of the first two tensor factors in $\G_a \o \G_a \o \C$.
It is not difficult to show that these formulae give
consistent structure data for the left chiral theory
(cf. also Appendix A.3). The motivation for eqs.~(\ref{lstruct}),
(\ref{rltrafo}) comes from the construction of the quantum group valued
field $g$. So let us define
\be{defg}
       g  \,:=\,  \S_a(\Phi^l)\  \Phi^r \ \ \in\G_a \o \W \ ,
\ee
where $\S_a(\Phi^l) \equiv (\S \o id)(\Phi^l)$. The element $g$
indeed preserves the constraint which we discussed above, that is
\be{fSf}
  {\rm if}\ \ \ \rf \ g \,= \, \S_{lr}(\rf)\ g \ , \ \ \
 {\rm then}\ \    g \ \rf \,= g \ \S_{lr}(\rf) \ \ \mbox{ for all }
       \ \ \rf \in \C^r \ .
\ee
Therefore, components of $g$ map the diagonal space $\H$ into
itself.  This remarkable property is established by a straightforward
computation (see Appendix A.4).

To study properties of $g$ it is helpful to have some knowledge
about the object $\S_a(\Phi^l)$. Simple applications of the
standard Hopf algebra axioms allow to deduce
\ba   
   \S_a(\Phi^l)\ \xi \ = \ \D_\iota(\xi)\ \S_a(\Phi^l) \ \ & , & \ \
   R_+ \ \S_a(\up{2}{\Phi}\^l)\ \up{1}{M}\^l \ = \
   \up{1}{M}\^l \ R_- \ \S_a( \up{2}{\Phi}\^l ) \ \ ,
       \label{SPcov}  \\[2mm]     \label{theta}
   \S_a(\Phi^l)  \,=\, (\Phi^l)^{-1} \theta_l \ \ \ & &
   \mbox{ with } \ \ \ \theta_l \in \G_a\o\C^l \ .
\ea
Here $\xi \in \J^l \cong \G$ in the first equation, 
$\S_a( \up{2}{\Phi}\^l )$ is a shorthand for 
$(id \o \S_a)( \up{2}{\Phi}\^l )$, and the relation (\ref{theta})
may be regarded as a definition of $\theta_l$.
The transformation laws of vertex operators show that $\theta_l$
commutes with $e \o \xi \in \G_a \o \J^l$ and hence $\theta_l \in
\G_a \o \C^l$. We can actually give an explicit formula for
$\theta_l$ in terms of $F_l$. If we assume for simplicity that
$\e_a(\Phi^l) = \re$ (cf. Subsection 2.4), then $\theta_l=\sum_\vs
f^1_\vs \S(f^2_\vs) \o f^3_\vs $, where $f^i_\vs$ come from the
expansion $F_l=\sum_\vs f^1_\vs\o f^2_\vs \o f^3_\vs$.

\begin{prop} {\em (Properties of $g$)} Let $g$ denote the object
defined in eq.~(\ref{defg}) and restricted to the subspace $\H$.
This element $g \in \G_a \o \End(\H)$ obeys the following
relations:  \label{gprop} \vspace*{-3mm}
\ba
     \up{2}{g}\  \up{1}{g} \, = \, \D_a(g)\ ,  && \ \
    R_\pm\ \up{2}{g} \ \up{1}{g} \ =\ \up{1}{g}\  \up{2}{g}\  R_\pm
    \label{qg} \ , \\[1mm]
    \up{1}{M}\^r \up{2}{g}\  R_-\, =\, \up{2}{g}\  R_+ \up{1}{M}\^r \ , &&
  \ \   \up{1}{M}\^l\  R_- \up{2}{g}\, = \,
   R_+ \up{2}{g}\   \up{1}{M}\^l \ , \label{qs}\\[1mm]
     M^l\  g \,  =  \,    g\  M^r \ ,  && \ \
      {\sf v} \   g \  {\sf v}^{-1}\, =\, g \ , \label{qr}
\ea
where in the last line ${\sf v} = \rv_r \rv_l^{-1}$ is a combination of
the ribbon elements $\rv_r \in \C^r$ and $\rv_l \in \C^l$.
Moreover, $g$ is normalized, $\e_a(g)=\re$, and invertible with
inverse $g^{-1} = \S_a(g)$.
\end{prop}
Proofs of all these relations are given in Appendix A.4.

Eqs.~(\ref{qg}) mean that $g$ obeys the defining relation of a
quantum group $\F \, = \, \mbox{\it Fun\/}_q(G)$. More precisely,
components of $g$ generate
the dual of the quantum algebra $\G$.  The elements $M^r$, $M^l$
furnish algebras of left- and right-invariant vector fields for $\F$
and they are related to each other by means of eq.~(\ref{qr}).
All these equations are well known in the theory of quantum
groups. In more geometric terms, they describe the deformed
co-tangent bundle $T_q^*G$ \cite{AF1}. \footnote{Similarly, the
algebra generated by components of $\Phi^r$, $M^r$ only, is a
deformation of the co-tangent bundle $T^* B$ of the Borel
subgroup of $G$ \cite{ByFa}.}

Let us now explain the pictorial presentation of the object $g$
(see Figure 4). First, recall that so far left and right
chiral objects lived on the same plane but on different
sides of the thick solid line and there was no interaction
between them. But if we want to consider the restriction from
$\cM^l \o \cM^r$ to the diagonal subspace $\H$, we have to
modify the rules. Namely, the restriction $\rf=\S_{lr}(\rf)$
enforces us to change the topology by gluing the plane into a
cylinder. Then we can join ends of dotted lines from both
sides of the thick solid line and thus combine objects of
different chirality. This is demonstrated by the graphical
presentation of $g$ in Figure 4 (the dashed line continues
the dotted line around the back side of cylinder and, hence,
cannot interfere with any line on the front side). Figure~4
also sketches the proof of the operator product expansions
of $g$ in eq.~(\ref{qg}).

Before concluding this subsection we would like to compare our
construction of $g$ with the one discussed in \cite{AF1}. There,
two decompositions of $g$ into triple products of elements, 
\hbox{$g \,=\, u\, Q^{-1} v \,=\, u_0 \, Q v_0$},  have been provided. 
All operators which appear in these relations act on the diagonal
subspace $\H$. The variables $v_0, u_0$ are chiral observables,
i.e., $u_0 \in \G_a \o\J^l, v_0 \in \G_a \o \J^r$, and hence
they commute with each other: $\up{1}{u}_0 \up{2}{v}_0 \, =
\, \up{2}{v}_0 \up{1}{u}_0$. Notice that components of $u_0,v_0$
leave the subspaces $V^{\bar I} \o V^I \subset \H$ invariant and
hence their actions are, in principle, expressible through 
the chiral objects $M^\a$ (in particular, $u_0,v_0$
are not to be confused with our vertex operators). The exchange
relations of
$u_0$ (respectively $v_0$) can be controlled only after
multiplication with $Q \in \G_a \o \End(\H)$. In fact, the
elements $u \,=\, u_0\, Q $ and $ v \, = \, Q\, v_0$ possess
the same exchange relations as our chiral vertex operators.
On the other hand, they are certainly not chiral any more
(because chiral vertex operators cannot act on $\H$).
In particular, $u$ does not commute with $v$. One may think
of $u$ (and similarly of $v$) as a left chiral vertex
operator dressed with a right chiral factor which leaves
the quadratic relations unchanged and, at the same time,
produces an operator acting on $\H$. Our construction in
terms of chiral vertex operators and the restriction from
$\cM$ to the diagonal subspace $\H$ is similar to
\cite{Fa1,FHT}.

\setlength{\unitlength}{.5pt}
\begin{figure}
\begin{picture}(400,200)(0,0)
\thicklines
\multiput(25,75)(0.25,0){6}{\line(1,1){50}}
\thinlines
\put(50,100){\line(-1,1){25}}
\multiput(75,75)(-2.8,2.8){9}{\circle*{1}}
\put(110,96){$\:\S_a(\Phi^l)$}
\end{picture} \\[-3.9cm]
\begin{picture}(400,180)(0,0)
\thicklines
\multiput(25,45)(0.25,0){6}{\line(1,1){100}}
\thinlines
\multiput(50,70)(-2.8,2.8){9}{\circle*{1}}
\put(75,45){\line(-1,1){25}}
\multiput(125,95)(-2.8,2.8){9}{\circle*{1}}
\put(100,120){\line(-1,1){25}}
\multiput(25,95)(8,0){13}{\line(1,0){5}}
\put(164,92){ $\: g$}
\end{picture}\hspace*{-1cm}
\begin{picture}(400,400)(0,0)
\thicklines
\multiput(80,220)(0.25,0){6}{\line(1,1){114}}
\thinlines
\multiput(110,250)(-2.8,2.8){9}{\circle*{1}}
\put(140,220){\line(-1,1){32}}
\put(122,238){\line(-1,-1){20}}
\multiput(185,275)(-2.8,2.8){9}{\circle*{1}}
\put(160,300){\line(-1,1){32}}
\put(148,312){\line(1,1){20}}
\multiput(85,275)(8,0){13}{\line(1,0){5}}
\put(220,275){=}
\put(-45,274){ $\Delta_a(g)\ =$}
\thicklines
\multiput(280,220)(0.25,0){6}{\line(1,1){114}}
\thinlines
\multiput(285,275)(2.8,-2.8){6}{\circle*{1}}
\put(300,240){\line(1,-1){25}}
\put(320,260){\line(1,-1){25}}
\multiput(300,240)(0,3.6){6}{\circle*{1}}
\multiput(320,260)(-3.6,0){6}{\circle*{1}}
\multiput(385,275)(-2.8,2.8){6}{\circle*{1}}
\put(350,290){\line(-1,1){25}}
\put(370,310){\line(-1,1){25}}
\multiput(350,290)(3.6,0){6}{\circle*{1}}
\multiput(370,310)(0,-3.6){6}{\circle*{1}}
\multiput(285,275)(8,0){13}{\line(1,0){5}}
\put(415,275){=}
\thicklines
\multiput(60,45)(0.25,0){6}{\line(1,1){100}}
\thinlines
\multiput(85,70)(-2.8,2.8){9}{\circle*{1}}
\put(110,45){\line(-1,1){25}}
\multiput(160,95)(-2.8,2.8){9}{\circle*{1}}
\put(135,120){\line(-1,1){25}}
\multiput(60,95)(8,0){13}{\line(1,0){5}}
\thicklines
\multiput(75,60)(0.25,0){6}{\line(1,1){100}}
\thinlines
\multiput(100,85)(-2.8,2.8){9}{\circle*{1}}
\put(125,60){\line(-1,1){25}}
\multiput(175,110)(-2.8,2.8){9}{\circle*{1}}
\put(150,135){\line(-1,1){25}}
\multiput(75,110)(8,0){13}{\line(1,0){5}}
\put(220,90){=}
\thicklines
\multiput(245,45)(0.25,0){6}{\line(1,1){66}}
\thinlines
\multiput(270,70)(-2.8,2.8){3}{\circle*{1}}
\put(295,45){\line(-1,1){25}}
\multiput(286,86)(2.8,-2.8){3}{\circle*{1}}
\put(286,86){\line(-1,1){25}}
\multiput(262,78)(8,0){4}{\line(1,0){5}}
\thicklines
\multiput(280,80)(0.25,0){6}{\line(1,1){66}}
\thinlines
\multiput(305,105)(-2.8,2.8){3}{\circle*{1}}
\put(330,80){\line(-1,1){25}}
\multiput(321,121)(2.8,-2.8){3}{\circle*{1}}
\put(321,121){\line(-1,1){25}}
\multiput(297,113)(8,0){4}{\line(1,0){5}}
\put(370,90){=}
\put(395,90){ $\stackrel{2}{g}\,\stackrel{1}{g}$}
\end{picture}
\begin{center}
\parbox{13cm}{ \small {\bf Figure 4:}
The definition of $g$ in terms of $\Phi^r$ and $\S_a(\Phi^l)$ is
shown on the left side. The right side of the figure is a
pictorial proof of the operator product expansion for $g$
(first eq.~in (\ref{qg})).}
\end{center}
 \end{figure}
\setlength{\unitlength}{1pt}

\subsection{Toy model for $\Z_q$.}

It is instructive to realize the constructions of the toy model
in the case of $\G=\Z_q$. Now we have two commuting copies, $h_\a$,
$\a=r,l$, of the element $h$ (see Subsection 2.1) generating the
chiral algebras $\J^\a$. We can also introduce Hermitian
operators $\wh{p}_\a$ such that $h_\a=q^{\wh{p}_\a}$. To introduce the
chiral monodromies $M^r$ and $M^l$ we use the expressions (\ref{Nz})
for the elements $N$ and $\wt{N}$. Since $M^r$ and $M^l$ differ from
them only by factors $v_a^{-1}$ and $v_a$, we get
\be{Mz}
 M^r \,=\,q^{\wh{p}\^{\;\; 2}\,\o\, e + 2\,\wh{p}\,\o\,\wh{p}_r}\ , \ \
 \ \ M^l \,=\,q^{-\wh{p}\^{\;\; 2}\,\o\, e - 2\,\wh{p}\,\o\,\wh{p}_l} \ .
\ee
The reader is invited to verify the functorial properties (\ref{rM}),
(\ref{lM}) for these objects (in fact, the check repeats the computations
performed in Subsection 2.2). As was explained in Subsection 2.2, the
element $S_\iota$ is trivial in the case of $\Z_q$, therefore
the chiral monodromies are unitary. The components of $M^\a$ act on
the model spaces $\cM^\a$.

Next we need to  construct the diagonal subspace
$\H=\oplus V^{\bar K} \o V^K$. It can be seen from the explicit formula
for the characteristic projectors (\ref{zP}) that $\S_a(P^s)=P^{-s}$,
i.e., the representation conjugate to $\tau^s$ is $\tau^{-s}$
(where $s$ is taken modulo $p$, $q^p=1$). Therefore, \hbox{$\H=\oplus
|\!-s\rangle \o |s\rangle$} is a $p$-dimensional subspace in
$p^2$-dimensional space $\cM^l\o\cM^r$. Using the operators
$\wh{p}_\a$, we can  characterize the subspace $\H$
as follows: $\wh{p}_r \phi=-\wh{p}_l\phi$ for all $\phi\in\H$.

Now we employ the construction for vertex operators which we
provided in Subsection 2.5. Let $\wh{Q}_\a$, $\a=r,l$  be
unitary operators acting on $\cM^l\o\cM^r$ such that
$\wh{Q}_r |s'\rangle\o|s''\rangle=|s'\rangle\o|s''+1\rangle$ and
$\wh{Q}_l |s'\rangle\o|s''\rangle=|s'+1\rangle\o|s''\rangle$. It is
convenient to introduce also two operators $\wh{\vs}_\a$ by
$\wh{Q}_\a=e^{\wh{\vs}_\a}$. With these notations it is easy to
verify that
$$
  \Phi^\a \,=\, \sum_{s=0}^{p-1} \, P^s \o \wh{Q}_\a\^{\!s} \,=\,
 e^{\wh{p}\,\o\,\wh{\vs}_\a } \, \in\G_a\o\End(\cM^\a) \ ,
  \ \ \a=r,l
$$
are right- and left chiral operators obeying all the properties spelled
out in Subsections 3.1 and 3.2, respectively. In particular, the
$R$-matrix commutation relations in (\ref{rcov}), (\ref{lcov}) boil
down to Weyl relations: $e^{e\,\o\,\wh{p}\,\o\,\wh{\vs}_\a }$
$q^{\pm 2\,\wh{p}\,\o\,e\,\o\,\wh{p}_r}$ =
$q^{\pm 2\wh{p}\,\o\,\wh{p}\,\o\,e}$
$q^{\pm2\,\wh{p}\,\o\,e\,\o\,\wh{p}_r}$
$e^{e\,\o\,\wh{p}\,\o\,\wh{\vs}_\a }$.
Recall that the universal $R$-matrices in (\ref{rcov}), (\ref{lcov})
are regarded as elements in $\G_a\o\G_a\o End(\H )$ with trivial entry in
the third tensor factor; hence, the factor $q^{\pm2\wh{p}\,\o\,\wh{p}
\,\o\,e}$ converts $R_\mp$ into $R_\pm$ (cf. Subsection 2.2).

Now, applying (\ref{defg}), we get an explicit expression for $g$:
\be{zg}
  g \,=\, \sum_{s=0}^{p-1} \, P^s \o \wh{Q}_r\^{\!s} \, \wh{Q}_l\^{-s}
 \,=\,  e^{\wh{p}\,\o\,(\wh{\vs}_r-\wh{\vs}_l) }  \ .
\ee
This object manifestly maps the diagonal subspace into itself and hence
we may regard $g$ as an element in $\G_a\o End(\H)$. The operator
product expansion (\ref{qg}) is obvious (see the analogous
computation for $N_\pm$ in Subsection 2.2). Moreover, the first
equation in (\ref{qr}) is again of Weyl-type (notice that here
the factors $v_a^{\mp 1}$ of $M^\a$ are essential):
$$
 \ar{c}  M^l\,g \,=\, q^{-\wh{p}\^{\;\; 2}\o\,e}\,
 q^{-2\,\wh{p}\,\o\,\wh{p}_l}\, e^{\wh{p}\,\o\,(\wh{\vs}_r-\wh{\vs}_l)}
 \,=\, q^{-\wh{p}\^{\;\; 2}\o\,e}\,e^{\wh{p}\,\o\,(\wh{\vs}_r-\wh{\vs}_l)}
 \,q^{-2\,\wh{p}\,\o\,\wh{p}_l}\,q^{2\,\wh{p}\^{\;\; 2}\o\,e}\,= \\ [1.5mm]
 =\, e^{\wh{p}\,\o\,(\wh{\vs}_r-\wh{\vs}_l) } \,
 q^{\wh{p}\^{\;\; 2}\o\,e}\, q^{-2\,\wh{p}\,\o\,\wh{p}_l}\,=\,
 \, e^{\wh{p}\,\o\,(\wh{\vs}_r-\wh{\vs}_l) } \,
 q^{\wh{p}\^{\;\; 2}\o\,e}\, q^{2\,\wh{p}\,\o\,\wh{p}_r}\,=\,
 g\,M^r \ . \er
$$
In the last line we used the constraint $\wh{p}_r=-\wh{p}_l$
valid on the diagonal subspace. To conclude, we notice that
in the $\Z_q$ case the vertex operators and the field $g$
are  unitary.

\section{REVIEW ON LATTICE CURRENT ALGEBRAS}
\setcounter{equation}{0}

In the previous section we have considered the toy model for
the WZNW theory which certainly did not go much beyond the
theory of vertex operators for quasi-triangular modular Hopf
algebras (except that we had two commuting copies of this
theory). Vertex operators for the infinite dimensional
current algebras of the WZNW-model depend in addition on
a spatial coordinate $x$. This brings new locality features
into the theory. Our aim is to describe them for a lattice
regularization of the WZNW-model developed in
\cite{AFS,AFSV,Fa1,FaGa,AFFS},
where the spatial coordinate assumes the discrete values,
$x=2\pi n/N$, $n=0,..,N-1$. We begin this discussion with a
brief review on lattice current algebras $\K_N$. Our notations
are close to those adopted in \cite{AFFS}.

\subsection{Definition of lattice current algebras.}

Let us consider a one-dimensional periodic lattice which consists of $N$
vertices. It is convenient to enumerate the vertices from 0 to $N-1$ and
the corresponding edges from 1 to N as shown below. 

\hbox{\begin{picture}(300,50)
    \put(90,20){\line(1,0){230}}
    \put(90,20){\circle*{3}} \put(130,20){\circle*{3}}
    \put(170,20){\circle*{3}} \put(280,20){\circle*{3}}
    \put(320,20){\circle*{3}}
    \put(87,25){$\large \G_0$} \put(127,25){$\large \G_1$}
    \put(167,25){$\large \G_2$} \put(273,25){$\large \G_{N-1}$}
    \put(311,25){$\large \G_N\equiv\G_0$}
    \put(107,8){$\large J_1$} \put(147,8){$\large J_2$}
    \put(297,8){$\large J_{N}$}
\end{picture} }
\begin{center}
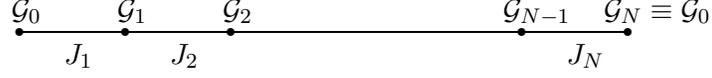
\parbox{13cm}{ \small {\bf Figure 5:} $N$-vertex periodic lattice.
 Each site is equipped with a copy of the symmetry algebra $\G$.
 The discrete currents $J_n$ are assigned to edges.  }
\end{center}

{}According to the ideology of \cite{AFFS}, the definition of
the algebra $\K_N$ involves two kinds of objects -- those associated
with the sites and those associated with the edges. The $n^{th}$ site
of the lattice is equipped with a copy $\G_n$ of the algebra $\G$ and
copies for different sites commute. In other words, $\K_N$ contains
$\G_n$ and the whole tensor product $\G^{\otimes N}$ as subalgebras.
The canonical isomorphism of $\G$ and $\G_n\subset \G^{\otimes N}$
furnishes the embeddings $\iota_n:\G\mapsto\G^{\otimes N}$ for
$n=0,..,N-1$:
$$ \iota_n(\xi) = e\o \dots\o\xi\o\dots\o e \in\G^{\otimes N} \ \
{\mbox{ for all}}\ \xi\in \G\ \ ,  $$
where the only nontrivial entry of the tensor product on r.h.s. 
appears in the $n^{th}$ position. The definition of $\K_N$ also 
involves generators $J_n^r$, $n\!=\! 1,..,N,$ (the right currents) 
which are discrete analogues of the continuum holonomies along 
the edges (cf. Introduction).

\begin{defn} {\em \cite{AFFS}}  \label{CAD}
The lattice current algebra $\K_N$ is generated by components of
invertible elements $J_n^r \in\G_a\o\K_N$, $n=1,..,N$ along with
elements in $\G^{\otimes N}$. These generators are subject
to the following relations :
\ba
 \up{2}{J}\^r_n  \up{1}{J}\^r_n \ = \  R_- \, \D_a( J^r) \ , &&
 (J^r_n)^*\, =\, S_n^{-1} (J^r_n)^{-1} S_{n-1} \ , \label{K2} \\ [1mm]
 \up{1}{J}\^r_{n+1}\, \up{2}{J}\^r_n \,=\,
 \up{2}{J}\^r_n \, R_+ \, \up{1}{J}\^r_{n+1} \label{K4}\ , &&
 \up{1}{J}\^r_n  \up{2}{J}\^r_m \,=\,  \up{2}{J}\^r_m \up{1}{J}\^r_n
 \ \ for\ n \neq m,m \pm 1\  (\mod N) \label{Jbr} \ , \\ [1.5mm]
 \iota_n(\xi)\, J_n^r \,=\, J_n^r \, \D_n' (\xi) \ ,  &&
  \D_{n-1}' \, (\xi) \, J_n^r \,=\, J_n^r \, \iota_{n-1}(\xi) \ \
  \mbox{ for all } \ \xi \in \G , \label{K1} \\ [1.5mm]
  \iota_m(\xi) \, J_n^r \,=\, J_n^r \, \iota_m(\xi) &&  \mbox{ for all }
  \ \xi \in \G,\; m \neq n,n-1\  (\mod N) \ . \nn
\ea
Here $R_\pm$ denote the elements $R_\pm\o \re\in\G_a\o\G_a\o\K_N$,
$S_n=(id\o\iota_n)(S) \in \G_a \o \G_n \subset \G_a \o \K_N$ with $S$
defined as in (\ref{SS}), and $\D'_n(\xi)=(id\o\iota_n)(\D'(\xi))$ where
$\D'(\xi)= P\D(\xi)P$ as usual. Invertibility of $J^r_n$ means that
there exists an element $(J^r_n)^{-1}\in\G_a\o\K_N$ such that
$J^r_n (J^r_n)^{-1} =  e\o \re = (J^r_n)^{-1}  J_n^r$.
\end{defn}

The lattice current algebra $\K_N$ contains a subalgebra $\J_N^r$
generated by components of the currents $J_n^r$ only. They are subject
to relations (\ref{K2})-(\ref{Jbr}).
The full lattice current algebra $\K_N$ can be regarded as a semi-direct
product  \hbox{$\G^{\o N}\sd \J_N^r$}, where the action of $\G^{\o N}$ 
on $\J_N^r$ is given by the covariance relations (\ref{K1}).

Taking into account the quasi-triangularity of $R_\pm$, we obtain the
following consequence of (\ref{K2})
\be{JJ}
 R_\pm\,\up{1}{J}\^r_n\,\up{2}{J}\^r_n\,R_\mp \,=\,
 \up{2}{J}\^r_n\,\up{1}{J}\^r_n \ \  .
\ee
These $R$-matrix relations for the description of the lattice Kac-Moody
algebras have been introduced first in \cite{AFS}. Following
our discussion in Subsection 2.2, one can introduce the objects
$N_{n,\pm}=(id\o \iota_n)(R_\pm) \in \G_a \o \G_n \subset \G_a \o \K_N$,
which obey the standard relations (\ref{Npmeq}). They are used to
rewrite the relations (\ref{K1}) in the following $R$-matrix form
\be {JN}
 \up{1}{N}_{n,\pm} \, \up{2}{J}\^r_n \,=\, \up{2}{J}\^r_n \,R_{\pm}\,
 \up{1}{N}_{n,\pm} \,, \;\;\;\;
 \up{2}{J}\^r_n \,\up{1}{N}_{n-1,\pm} \,=\, R_{\pm}\,
 \up{1}{N}_{n-1,\pm} \, \up{2}{J}\^r_n \ \ .
\ee

\subsection{Left currents.}
The continuum WZNW model possesses two chiral subalgebras, that is,
along with the (right) current $j^r(x)$ it involves the left current
$j^l(x)$ such that left and right currents commute. A nice feature
of the lattice current algebra $\K_N$ is that it already contains the
second chirality in an encoded form. Indeed,  one may introduce the
following new variables $J^l_n \in \G_a \o \K_N$:
\be{nchir}
  J^l_n \ := \ v_a^2\, N_{n-1,+}^{-1} \, J^r_n \, N_{n,-}  \ .
\ee
In the notations of Definition \ref{CAD} they obey (see \cite{AFFS}
for details)
\ba
 \up{1}{J}\^l_n\, \up{2}{J}\^l_n \, = \, R_+\, \D_a (J^l_n)  &,&
 (J^l_n)^*\, =\, S_n\, (J^l_n)^{-1} S_{n-1}^{-1} \ , \label{DI} \\ [1mm]
  \up{1}{J}\^l_{n} \, R_- \, \up{2}{J}\^l_{n+1} \,=\,
 \up{2}{J}\^l_{n+1} \,  \up{1}{J}\^l_{n} &,&
 \up{1}{J}\^l_n \up{2}{J}\^l_m \,=\, \up{2}{J}\^l_m \up{1}{J}\^l_n \ \
 \mbox{ for } \ \ n \neq m,m \pm 1\  (\mod N) \ , \label{II} \\ [1mm]
  \up{1}{J}\^l_n \up{2}{J}\^r_m \,=\, \up{2}{J}\^r_m \up{1}{J}\^l_n
 &&  \mbox{ for all}\ \ \  m,\,n ,  \label{JI} \\ [1.5mm]
  \iota_n(\xi)\, J_n^l \,=\, J_n^l \, \D_n (\xi)\ \   &,& \ \
  \D_{n-1} \, (\xi) \, J_n^l \,=\, J_n^l \, \iota_{n-1}(\xi) \ \
  \mbox{ for all } \ \xi \in \G , \label{FI} \\ [1.5mm]
  \iota_m(\xi) \, J_n^l \,=\, J_n^l \, \iota_m(\xi) \ \ & & \ \
   \mbox{ for all }\ \xi \in \G,\; m \neq n,n-1\  (\mod N) \ . \nn
\ea
Due to these properties, the objects $J^l_n$ may be interpreted as
left counterparts of the right currents $J^r_n$. Notice that there is
a manifest symmetry between the defining relations for the right currents
and the properties of left currents. It underlines the fact that left and
right chiralities in the WZNW model appear on an equal footing.
In fact, (\ref{DI})-(\ref{FI}) could be regarded as an alternative
definition of the lattice current algebra $\K_N$.

It also follows that $\J^l_N$ and $\J^r_N$, i.e., the algebras
generated by components of $J^l_n$ and $J^r_n$, respectively,
are commuting subalgebras in $\K_N$ and $\J^r_N$ is isomorphic to
$(\J^l_N)_{op}$. Here the subscript $\ _{op}$ means opposite
multiplication as before.

\subsection{Holonomies and monodromies.}
The currents $J_n^\a, \a = r,l$ were defined as discrete
analogues of holonomies along the $n^{th}$ edge. Similarly,
one may introduce the holonomies along the link connecting
the $0^{th}$ and the $n^{th}$ sites :
\be{4.6}
 U_n^{\a}\,:=\,J_1^{\a} \dots J_n^{\a}, \;\;\;\; n=1,\dots,N-1\ \ .
\ee
As one might expect, the properties of such holonomies are similar
to those of chiral currents.\footnote{
Let us mention here some subtle point in the definition of the lattice
current algebra. Notice that relations (\ref{JJ}) would not change if we
replaced $R_-$ by $R_+$ in the definition (\ref{K2}). However,
this ambiguity disappears if we demand that $U^r_n$ and $J^r_n$ have
the same functoriality relation (compare (\ref{K2}) and first eq.~in
(\ref{DU})). A similar subtlety appears once more in the construction
of the left currents. Indeed, we could replace factor $v^2_a$ by $v_a$
in the definition (\ref{nchir}); then we would obtain the relation
(\ref{DI}) with $R_-$ instead of $R_+$. But in this case functorial
properties of $J^l_n$ and $U^l_n$ would be different.}
Namely, it is easy to verify that
\ba
 \up{2}{U}\^r_n\ \up{1}{U}\^r_n \,=\, R_- \ \D_a (U_n^r) \ \ &,&\ \
 \up{1}{U}\^l_n\ \up{2}{U}\^l_n \,=\, R_+ \ \D_a (U_n^l)
 \label {DU} \ , \\ [1.5mm]
 (U^r_n)^*\, =\, S_n^{-1} (U^r_n)^{-1} S_0 \ \  &,& \  \
 (U^l_n)^*\, =\, S_n \, (U^l_n)^{-1} S_0^{-1} \ , \label{cU} \\ [1mm]
 R_\pm\  \up{1}{U}\^r_n\ \up{2}{U}\^r_n\ R_\mp \,=
 \, \up{2}{U}\^r_n\ \up{1}{U}\^r_n \ \  &,&\ \
 R_\pm \,\up{2}{U}\^l_n\,\up{1}{U}\^l_n\,R_\mp \,=
 \,\up{1}{U}\^l_n\,\up{2}{U}\^l_n \ , \label{uu} \\ [1.5mm]
 \D_0'(\xi) \ U_n^r\,  = \, U_n^r \ \iota_0(\xi) \ \ &,&\ \
 \iota_n(\xi) \ U_n^r \, = \, U_n^r \ \D_n' (\xi) \ \
 \mbox{ for all }\ \  \xi\ \in\ \G, \label{u'} \\ [1.5mm]
 \D_0(\xi) \ U_n^l  \, = \, U_n^l \ \iota_0(\xi)\ \  &,&\ \
 \iota_n(\xi) \ U_n^l \, =\,  U_n^l \ \D_n (\xi) \ \
 \mbox{ for all }\ \  \xi\ \in\ \G \label{u''}
\ea
and $U_n^{\a}$ commute with $\iota_m(\xi)$ for all $m\neq 0, n$.
However, there is an important difference between currents and
holonomies: since the latter are localized on the chain of edges 
that runs from the $0^{th}$ vertex to the $n^{th}$, the localization 
domains of all holonomies overlap. This is reflected in their
mutual exchange relations for $1 \leq n < m \leq N\!-\!1$:
\be{UU}
 \up{2}{U}\^r_n\,\up{1}{U}\^r_m \,=\,
 R_-\,\up{1}{U}\^r_m\,\up{2}{U}\^r_n \ \ \ , \ \ \
 \up{1}{U}\^l_n\,\up{2}{U}\^l_m \,=\,
 R_+\, \up{2}{U}\^l_m \,\up{1}{U}\^l_n \ .
\ee

As we have argued in the introduction, holonomies of chiral fields
along the whole circle (i.e., the chiral monodromies) are of
particular interest. In the continuum case they are given by
$ m^\a = {\cal P}\, exp\{\oint j^\a (x) dx\}$.
Monodromies for the quantum lattice theory are defined by a
natural discrete analogue of this formula,
\be{mon}
     M^{\a}\,=\,J^\a_1 \,J^\a_2 \dots \,J^\a_N \ .
\ee
Simple calculations allow to derive the following properties of
the monodromies $M^\a$:
\ba
 \up{2}{M}\^r R_+ \up{1}{M}\^r \,=\, R_- \D_a (M^r)\ \  & , & \ \
 \up{1}{M}\^l R_- \up{2}{M}\^l \,=\,
 R_+ \D_a (M^l) \ , \label{DM} \\ [1mm]
 (M^r)^*\, =\, S_0^{-1} (M^r)^{-1} S_0 \ \  &,& \ \
 (M^l)^*\, \,= S_0 \, (M^l)^{-1} S_0^{-1} \ ,  \label{cM} \\ [1mm]
 R_+\,\up{1}{U}\^r_n\,\up{2}{M}\^r \,=\, \up{2}{M}\^r \,
 R_+\,\up{1}{U}\^r_n  \ \ & , & \ \
 R_-\,\up{2}{U}\^l_n\,\up{1}{M}\^l \,= \,\up{1}{M}\^l \,R_-\,
 \up{2}{U}\^l_n \  ,  \label{UM} \\ [1mm]
 \D_0'(\xi)\,M^r \, =\, M^r \, \D_0'(\xi) \ \ & , & \  \
 R_\pm \, \up{1}{N}_{0,\pm} \, \up{2}{M}\^r \,=\,
 \up{2}{M}\^r \, R_\pm \, \up{1}{N}_{0,\pm} \ , \label{Mrcov} \\ [1mm]
 \D_0(\xi)\, M^l\, =\, M^l \, \D_0 (\xi) \  \  &,& \ \
 \up{1}{N}_{0,\pm} \,R_\pm \, \up{2}{M}\^l \,=\,
 \up{2}{M}\^l \,  \up{1}{N}_{0,\pm} \,R_\pm \ \label{Mlcov}
\ea
for all  $\xi\in \G$ and $M^\a$ commute with $\iota_m(\xi)$ for all
$m\neq 0$.

Now we see that the relations (\ref{rM}) and (\ref{lM}) which we postulated
in the toy model construction indeed describe properties of the chiral
monodromies. Our next aim is to extend the toy model to the full
lattice theory. Recall that the structure data of the toy model were
built from elements in the center $\C^\a$ of the algebra $\J^\a$
spanned by components of $M^\a$. Elements in these algebras $\C^\a$
are still central in the full lattice theory. 
In fact, it follows from (\ref{MMr}), (\ref{MMl}) that the
algebras $\C^\a$ are spanned by the elements $c^I_\a= \tr_q^I \t^I(M^\a)$
where $\tau^I$ runs through irreducible representations of $\G$
and $\t^I(M^\a) = ( \t^I \o id)( M^\a)$ \cite{AGS}. Equipped with this
explicit description of $\C^\a$ one concludes from 
eqs.~(\ref{UM})-(\ref{Mlcov}) that the elements $c^I_\a$ commute 
with $U^\a_n$, $N_{n,\pm}$ for all $n$ and hence they are central 
elements in $\K_N$.
Actually, the following stronger statement holds \cite{AFFS}: the
elements $c^I_\a \in \K_N $ with $I$ running through the classes of
irreducible representations of $\G$ generate the full center of the
lattice current algebra $\K_N$. This explains why the structure data
for vertex operators on the lattice will be built from the commuting
subalgebras $\C^\a$ exactly as in the toy model.

\subsection{Current algebra for $\Z_q$.}
Let us consider the current algebra in the case of  $\Z_q$. The
algebras $\G_n$ assigned to the sites are generated by the elements
$h_n=\iota_n(h)$. As usual, we can introduce $\wh{p}_n$ such that
$h_n=q^{\wh{p}_n}$. For $N_{n,\pm}\in\G_a\o\G_n$ we have
$N_{n,\pm}=\sum_s P^s\o h_n^{\pm s} = q^{\pm\wh{p}\,\o\,\wh{p}_n}$
(cf. Subsection 2.2). Next we build the chiral currents
\be{cz}
   J^r_n \,=\, \sum_{s=0}^{p-1} \, q^{\frac 12 s^2} \,
 P^s\o (\wh{W}_n^r)^s \ \ \ \ , \ \ \ \ \
  J^l_n \,=\, \sum_{s=0}^{p-1} \, q^{-\frac 12 s^2} \,
 P^s\o (\wh{W}_n^l)^s \
\ee
{}from a family of unitary elements $\wh{W}_n^r$ and
$\wh{W}_n^l := h_{n-1}^{-1} \wh{W}_n^r h_n^{-1}$ which obey
the following Weyl-type relations
\ba
 h_n \, \wh{W}_n^\a \,=\, q \, \wh{W}_n^\a \, h_n \ ,&& h_{n-1}\,
 \wh{W}_n^\a \,=\,q^{-1}\, \wh{W}_n^\a \, h_{n-1} \ ,\label{hW}\\ [1mm]
 h_m \, \wh{W}_n^\a \,=\, \wh{W}_n^\a \, h_m \ &&
   {\rm for}\ m\neq n,n-1 \ , \nn \\ [1mm]
 \wh{W}_{n+1}^r \, \wh{W}_n^r \, =\, q\, \wh{W}_n^r \,
  \wh{W}_{n+1}^r \ , &&
 \wh{W}_{n+1}^l \, \wh{W}_n^l \, =\, q^{-1}\, \wh{W}_n^l \,
  \wh{W}_{n+1}^l \ , \nn \\ [1mm]
 \wh{W}_{n}^\a \, \wh{W}_m^\a \, =\, \wh{W}_m^\a\,
  \wh{W}_{n}^\a \ &&  {\rm for}\ m\neq n\pm 1 \ \ . \label{WW}
\ea
Since the element $(\wh{W}_n^\a)^p$ is obviously central for this
algebra, we additionally impose the condition: $(\wh{W}_n^\a)^p=e$
for all $n$ (which is, in fact, a choice of a normalization).

The algebra generated by $\wh{W}_n^\a$ is known as lattice $U(1)$-current
algebra \cite{FV1,Fa3}. The relation we have used to obtain the
elements $\wh{W}_n^l$ from the $h_m$ and $\wh{W}_n^r$ is a special
case of formula (\ref{nchir}) and it implies that $\wh{W}_n^l
\wh{W}_m^r = \wh{W}_m^r \wh{W}_n^l$ for all pairs $n,m$.

The functorial properties (\ref{K2}) and (\ref{DI}) of currents
(\ref{cz}) can be checked in the same way as we did this for the
elements $N_\pm$ in Subsection 2.2. The exchange relations
(\ref{K4}), (\ref{JN}) and (\ref{II}) are again reduced to
Weyl relations. Since $S_n=e\o e$, the chiral currents are
unitary, $(J^\a_n)^*\,=\, (J^\a_n)^{-1}$, which is in agreement
with the unitarity of $\wh{W}_n^\a$.

To proceed, we introduce anti-Hermitian  operators $\wh\vp_n^\a$
such that $\wh{W}_n^\a = e^{\wh\vp_n^\a}$.\footnote{Strictly speaking,
the algebra generated by $\wh\vp_n^\a$ and $\wh{p}_n$ is larger than
one generated by $\wh{W}_{n}^\a$ and $h_n$ (see, e.g., \cite{Fa3}).
The latter is called compactified form of the former. }
In these notations the commutation relations given above acquire the form:
\be{zv}
 \ar{ll} [\, \wh\vp_{m}^r \,,\,\wh\vp_n^r \,] \,=\,
   \ln q \, (\dl_{m,n+1} - \dl_{m,n-1})  \ \ \ , & \
   [\, \wh\vp_{m}^l \,,\,\wh\vp_n^l \,] \,=\,
   \ln q \,( \dl_{m,n-1} - \dl_{m,n+1} )  \ , \\ [1.5mm]
   [\, \wh{p}_m \,,\,\wh\vp_n^\a \,] \,=\, \dl_{m,n} - \dl_{m,n-1}\ , & \
   [\, \wh\vp_{m}^l \,,\,\wh\vp_{n}^r \,] \,=\, 0 \ .\er
\ee
The chiral currents now can be rewritten in the following form:
\be{cz'}
  J^r_n \,=\, \k_a^{-1}\, e^{\wh{p}\,\o\,\wh\vp_n^r} \,=\,
  e^{ (\frac 12\ln{q})\,\wh{p}\^{\;\;2} \o\,e+\wh{p}\,\o
   \,\wh\vp_n^r} \ , \ \ \ \ \
  J^l_n  \,=\, \k_a \, e^{\wh{p}\,\o\,\wh\vp_n^l} \,=\,
  e^{-(\frac 12\ln{q})\, \wh{p}\^{\;\;2} \o\,e + \wh{p}\,\o
   \,\wh\vp_n^l} \ .
\ee

Next we can construct the chiral holonomies. For this purpose the
variables $\wh\vp_n^\a$ are more convenient. Indeed, applying the 
special case of the Campbell-Hausdorff formula,
$e^a e^b = e^{a+b} e^{\frac 12 [a,b]}$ valid if
$[a,[a,b]]=[b,[a,b]]=0$, we easily obtain:
$$
  U^r_n \,=\, q^{\frac 12 \wh{p}\^{\;\;2}\,\o\,e}\,
  e^{\wh{p}\,\o\,\sum_{k=1}^n  \wh\vp_k^r} \ , \ \ \ \ \
  U^l_n \,=\, q^{-\frac 12 \wh{p}\^{\;\;2}\,\o\,e}\,
  e^{\wh{p}\,\o\,\sum_{k=1}^n \wh\vp_k^l} \ .
$$
It is obvious now why relations (\ref{DU})-(\ref{u''}) for the holonomies
copy those for the currents. The exchange relations (\ref{UU}) are again
reduced to Weyl-type relations.

{}Proceeding in the same way, we can construct the chiral monodromies
as $M^\a=U^\a_{N-1}J^\a_N$, which needs again an application of
the Campbell-Hausdorff formula. The result reads
\be{zM}
 M^r \,=\, q^{ \wh{p}\^{\;\;2}\,\o\,e}\,
  e^{\wh{p}\,\o\,\sum_{k=1}^N  \wh\vp_k^r} \ , \ \ \ \ \
  M^l \,=\, q^{- \wh{p}\^{\;\;2}\,\o\,e}\,
  e^{\wh{p}\,\o\,\sum_{k=1}^N \wh\vp_k^l} \ .
\ee
Bearing in mind that for $\Z_q$ the quantum trace coincides with the
standard one (see Subsection 2.1), we conclude from (\ref{zM}) that
the algebras $\C^\a$ are generated by exponentials of the elements
${\sf\wh{p}}_\a=\sum_{k=1}^N  \wh\vp_k^\a$. Indeed, using the commutation
relations given above, it is easy to verify that these operators
commute with all elements of the currents algebra. Performing a formal
replacement $\sum_{k=1}^N  \wh\vp_k^\a \rar \pm (2 \ln q)\,\wh{p}_\a$
in(\ref{zM}) (the sign depends on the chirality), we recover the
formulae (\ref{Mz}) of the toy model.

\section{VERTEX OPERATORS ON A LATTICE}
\setcounter{equation}{0}
\subsection{Definition of $\W_N$.}

In Section 3 we have considered algebras $\J^\a$ and $\V^\a$, $\a=r,l$
generated by components of the chiral monodromies $M^\a$ and the chiral
vertex operators $\Phi^\a$, respectively. Both chiralities together
were used to generate the algebra $\W = \V^l \o \V^r$ of our toy
model. Below we shall define an algebra $\W_N$ of vertex operators
on a lattice. For this purpose, we shall replace the algebras $\J^\a$ in
the definition of $\W$ by their lattice counterparts $\J^\a_N$. So
we assume that we are given the lattice current algebra $\K_N$ with
center $\C^l \o \C^r$ (recall that $\C^\a \cong \C \cong$ center of
$\G$) and two sets of structure data $F_\a, \s_\a, \RR^\a_\pm, D_\a$, 
$\a = l,r,$ which obey the standard relations. The last
tensor components of the structure data are regarded as elements
in the center of the lattice current algebra $\K_N$, i.e., we
have $F_\a \in \G_a \o \G_a \o \C^\a \subset \G_a \o \G_a \o
\K_N$ etc.

\begin{defn}\label{W_N}{\em (Algebra of vertex operators on a lattice)}
The algebra $\W_N$ is generated by elements in $\K_N$ and components
of the vertex operators $\Phi^\a_0 \in \G_a \o \W_N$. Generators
$N_n, J^\a_n\in\G_a\o\K_N$ obey the defining
relations (\ref{K2}),(\ref{K4}) and (\ref{JN}) for lattice
current algebras. The elements $\Phi^\a_0 \in\G_a\o\V^\a \subset
\G_a \o \W_N$, $\a=r,l$ are subject to the following conditions:
\begin{enumerate}
\item
 They satisfy operator product expansions and exchange relations
 with elements in the center $\C^l \o \C^r \subset \K_N$ given by
 \ba   \label{OPE0r}
  \up{2}{\Phi}\^r_0 \ \up{1}{\Phi}\^r_0   =  \ F_r \D_a (\Phi^r_0)
  \ \ \ \ &, &\ \ \ \ \Phi^r_0\ \rf_r\ =\s_r(\rf_r)\ \Phi^r_0 \ , \\[1mm]
  \up{1}{\Phi}\^l_0 \ \up{2}{\Phi}\^l_0   =  \ F_l \D_a (\Phi^l_0)
  \ \ \ \ &, & \ \ \ \ \Phi^l_0 \ \rf_l \ =\s_l(\rf_l)\ \Phi^l_0 \ .
   \label{OPE0l}
 \ea
 Here $F_\a \in \G_a \o \G_a \o \C^\a$; $\s_\a$ are homomorphisms
 from $\C^\a$ to $\G_a \o \C^\a$ and $\rf_\a \in \C^\a$. Moreover,
 $\Phi^\a_0$ are invertible and vertex operators of different
 chirality commute:
 \be{W'}
  (\Phi_0^\a)^{-1} \Phi_0^\a \,=\, e\o \re \,=\, \Phi_0^\a \,
  (\Phi_0^\a)^{-1} \ ,\ \a=l,r\ ,\ \ \ \ \up{1}{\Phi}\^r_0 \,
  \up{2}{\Phi}\^l_0 \,=\, \up{2}{\Phi}\^l_0 \, \up{1}{\Phi}\^r_0 \ .
 \ee
\item
 $\Phi^r_0$ and $\Phi^l_0$ are chiral vertex operators for the
 algebra $\K_N$ in the sense that the following exchange relations
 with $J^r_n, N_{n,\pm} \in \G_a \o \K_N$ hold:
 \ba
  \up{1}{J}\^r_1\,\up{2}{\Phi}\^r_0 \,=\,
  \up{2}{\Phi}\^r_0 \, R_+ \, \up{1}{J}\^r_1 \ &,& \
  \up{2}{\Phi}\^r_0 \, \up{1}{J}\^r_N \,=\,
  \up{1}{J}\^r_N \, \up{2}{\Phi}\^r_0 \, R_- \ ,\label{FJp}\\[1mm]
  \up{1}{\Phi}\^r_0\,\up{2}{J}\^r_n \,=\,
  \up{2}{J}\^r_n\,\up{1}{\Phi}\^r_0 \ \mbox{ for }
  \ n \neq  1, N \ & , & \
  \up{2}{J}\^r_n \ \up{1}{\Phi}\^l_0\, = \, \up{1}{\Phi}\^l_0\
  \up{2}{J}\^r_n \ \mbox{ for all } \ n \label{Fj}\ ,   \\[1mm]
  \up{1}{N}_{0,\pm} \, \up{2}{\Phi}\^r_0 \,=\,
  \up{2}{\Phi}\^r_0 \, R_\pm \, \up{1}{N}_{0,\pm} \ \ &,& \ \
  \up{1}{N}_{0,\pm} \, \up{2}{\Phi}\^l_0 \,=\,
  \up{2}{\Phi}\^l_0 \, \up{1}{N}_{0,\pm} \, R_\pm \ ,
  \label{NP0}
 \ea
  and components of $N_{m,\pm}$ commute with components of the
  vertex operators for \hbox{$m\neq 0$}.
\item
 The $*$-operation on $\K_N$ can be extended to $\W_N$ by the
 following prescription
 \be{staronphi}
  (\Phi^r_0)^* \,=\, (S_0)^{-1}\,(\Phi^r_0)^{-1}\ \  \ ,\ \ \
  (\Phi^l_0)^* \,=\, S_0 \,(\Phi^l_0)^{-1}\ , \label{W2}
 \ee
 where $S_0=(id\o \iota_0)(S) \in \G_a \o \G_0 \subset \G_a\o\K_N$
 with $S \in \G_a \o \G$ being constructed by formula (\ref{SS}).
\end{enumerate}
Further relations involving left currents $J^l_n$ and the monodromies
$M^\a$ follow and will be spelled out below.
\end{defn}

Let us underline once more that the structure data for $\Phi^\a_0$
are constructed from elements in the center $\C^l \o \C^r$ of $\K_N$.
This is possible because both algebras $\C^\a$ are isomorphic to
the center of $\G$ \cite{AFFS} (see our short discussion at the end
of Subsection 4.3).

Next we would like to supplement our definition of $\W_N$ by a list
of consequences which follow from the stated relations. They concern
exchange relations of chiral vertex operators with left currents
$J^l_n$ and elements $\xi \in \G_n \subset \K_N$,
\ba
 \up{2}{J}\^l_1\,\up{1}{\Phi}\^l_0 \, = \,
 \up{1}{\Phi}\^l_0 \, R_- \, \up{2}{J}\^l_1 &,&
 \up{1}{\Phi}\^l_0 \, \up{2}{J}\^l_N \,=\,
 \up{2}{J}\^l_N \, \up{1}{\Phi}\^l_0 \, R_+ \ , \label{FJm} \\ [1mm]
 \up{1}{\Phi}\^l_0\,\up{2}{J}\^l_n \,=\,
 \up{2}{J}\^l_n\,\up{1}{\Phi}\^l_0 \ \mbox{ for } \
 n \neq 1,N \ & , & \  \up{2}{J}\^l_n \, \up{1}{\Phi}\^r_0 \ = \
 \up{1}{\Phi}\^r_0 \, \up{2}{J}\^l_n \ \mbox{ for all } \ n  ,
 \label{Fj2} \\[2mm]
 \iota_0(\xi)\, \Phi^r_0 \, =\, \Phi^r_0 \, \D_0' (\xi) \ \  &,& \ \
 \iota_0(\xi)\, \Phi^l_0 \, =\, \Phi^l_0 \, \D_0 (\xi) \ ,
 \label{W0} \\[1mm]
 \iota_m(\xi) \, \Phi^r_0 \,=\, \Phi^r_0 \, \iota_m(\xi) \ \ &,& \ \
 \iota_m(\xi) \, \Phi^l_0 \,=\, \Phi^l_0 \, \iota_m(\xi) \ \
 \mbox{ for }\ m \neq 0 \ ,
\ea
for all $\xi\in\G$ and we used the same notations as in Subsection
4.1. The first set of relations, i.e. eqs.~(\ref{FJm}),(\ref{Fj2}),
are obtained with the help of eq.~(\ref{nchir}). {}From our earlier
discussion it is clear that the relations (\ref{W0}) are equivalent
to eqs.~(\ref{NP0}). All exchange relations of the chiral vertex
operators with elements in $\K_N$ are local in the sense that
objects assigned to sites $n \neq 0$ or links $m \neq 1,N$
commute with $\Phi^\a_0$. This means that we can think of
$\Phi^\a_0$ as being assigned to the vertex $n = 0$ and hence
explains the subscripts $\ _0$.
The precise form of the nontrivial exchange relations involving
$\Phi^\a_0$ may be understood in terms of co-actions of $\G$
on $\K_N$ (see remarks in the introduction and \cite{AFS,AFFS}).

Now let us compare the definition of algebras $\W_N$ of vertex
operators on the lattice with our toy model. To
this end we derive exchange relations between the chiral vertex
operators $\Phi^\a_0$ and the chiral monodromies $M^\a$ (see eqs.
(\ref{mon})),
\be{MP0}
 \up{1}{M}\^r\ \up{2}{\Phi}\^r_0 \ R_- \ = \
 \up{2}{\Phi}\^r_0 \  R_+ \  \up{1}{M}\^r \ \ \ \ , \ \ \  \
 \up{2}{M}\^l\ \up{1}{\Phi}\^l_0 \  R_+ \ = \
 \up{1}{\Phi}\^l_0 \  R_- \  \up{2}{M}\^l\ \ .
\ee
The answer is to be compared with the relations (\ref{rcov}),
(\ref{lcov}) in the toy model and shows that the objects
$(\Phi^\a_0, M^\a, N_0)$ of the algebra $\W_N$ obey the
same exchange algebra as $(\Phi^\a, M^\a, N)$  in the
toy model. Thus, the toy model may not only be considered
as a special case of a lattice theory with $N=1$ but also
it is embedded as a subalgebra in all lattice algebras
$\W_N$ for arbitrary $N$. We can use this insight to rewrite some
of the relations we discussed for the toy model in terms of the
corresponding lattice objects. In particular, one has
\ba \label{RPP0}
 \RR^r_\pm\,\up{2}{\Phi}\^r_0\,\up{1}{\Phi}\^r_0 \,=\,
 \up{1}{\Phi}\^r_0\,\up{2}{\Phi}\^r_0 \,R_\pm \ \ &,& \ \
 \RR^l_\pm\,\up{1}{\Phi}\^l_0\,\up{2}{\Phi}\^l_0 \,=\,
 \up{2}{\Phi}\^l_0\,\up{1}{\Phi}\^l_0 \,R_\pm \ , \\[1mm]
  v_a^{-1}\, (\Phi^r_0)^{-1}\, D_r \, \Phi^r_0 \,=\,M^r \ \ &,& \ \ 
  v_a \, (\Phi^l_0)^{-1}\, D_l \, \Phi^l_0 \,=\,M^l \ \ \mbox{with}
  \label{MDM} \\[1mm]      \label{DD}
  D_r\,=\, (v_a  \rv_r^{-1}) \, \s_r(\rv_r) \ \ &,&  \  \
  D_l\,=\, (v_a^{-1} \rv_l) \ \s_l(\rv_l^{-1})  \ \ .
\ea
Here, $\rv_\a\in\C^\a$ are images of the ribbon element $v \in \G$
under the canonical isomorphisms from the center of $\G$ into the
subalgebras $\C^\a \subset \K_N$ (cf.~eq.~(\ref{Dop})). 
The latter are generated by
quantum traces of monodromies, i.e., by elements of the form
$\tr^I_q \t^I(M^\a)$ (for notations see Subsection 4.3). The
elements $\RR^\a_\pm$ in eqs.~(\ref{RPP0}) are given through
the standard formula $\RR^\a_\pm= F_\a' R_\pm (F_\a)^{-1}\in
\G_a\o\G_a\o\C^\a$.

\subsection{Vertex operators at different sites.}

Definition 3 involves only vertex operators assigned to the $0^{th}$
site of the lattice. We may now try to construct vertex operators
$\Phi^\a_n \in \G_a \o \W_N$ from elements in the algebra $\W_N$
which are assigned to other sites $n \neq 0$. In particular, they
are required to satisfy the characteristic fusion and braid
relations of vertex operators and, moreover, we want them to
commute with all elements in $\K_N$ which are assigned to sites
$m \neq n$ or edges $m \neq n, n+1$. The solution to this problem
is certainly not unique. In the following, we shall describe just one
possible construction. The idea is to introduce the vertex
operators $\Phi^\a_n$ with the help of the holonomies $U^\a_n\in
\J^\a_N$ by the simple formulae:
\be{Pn}
 \Phi^r_n \ :=\ \Phi^r_0\ U^r_n\ \ \ \ ,\;\;\;\;
 \Phi^l_n\ :=\ \Phi^l_0\ U^l_n\ \ \mbox{ for } \;\;\; n=1,
 \dots,N-1 \  .
\ee
Using the relations (\ref{DU})-(\ref{u''}) for chiral
holonomies, it is easy to verify the following properties of
$\Phi^\a_n$:
\ba
 \up{1}{N}_{n,\pm} \, \up{2}{\Phi}\^r_n \,=\,
 \up{2}{\Phi}\^r_n \, R_\pm \, \up{1}{N}_{n,\pm} \ \ &,& \ \
 \up{1}{N}_{n,\pm} \, \up{2}{\Phi}\^l_n \,=\,
 \up{2}{\Phi}\^l_n \, \up{1}{N}_{n,\pm} \, R_\pm \ , \label{NPn} \\[1mm]
 (\Phi^r_n)^* \,=\, (S_n)^{-1} (\Phi^r_n)^{-1}  \ \ &,& \ \
 (\Phi^l_n)^* \,=\, S_n \, (\Phi^l_n)^{-1} \ , \label{*Pn} \\ [1mm]
 \iota_n(\xi)\,\Phi^r_n \,=\, \Phi^r_n \,\D'_n(\xi) \ \ &,& \ \
 \iota_n(\xi)\,\Phi^l_n \,=\, \Phi^l_n \,\D_n(\xi) \ \ \
 {\rm for\ all}\ \xi\in\G  \ , \label{xPn}
\ea
and $\iota_m(\xi)$ commute with $\Phi^\a_n$ for any $m\neq n$.
Here $S_n=(id\o \iota_n)(S)$ with  $S \in\G_a\o \G$ as before.
Next, one has to investigate fusion and braiding properties of
$\Phi^\a_n$. The computation (see Appendix A.5) reveals that
the elements $\Phi^\a_n$ obey the same relations as our
vertex operators $\Phi^\a_0$ at the $0^{th}$ site, i.e.
\ba
 \up{2}{\Phi}\^r_n\,\up{1}{\Phi}\^r_n\, = \; F_r \,\D_a(\Phi^r_n) \ \
 &,& \ \  \up{1}{\Phi}\^l_n\,\up{2}{\Phi}\^l_n\, = \;
 F_l \,\D_a(\Phi^l_n)\ ,  \label{DFnp}  \\ [1mm]
 \RR^r_\pm\,\up{2}{\Phi}\^r_n\,\up{1}{\Phi}\^r_n \,=\,
 \up{1}{\Phi}\^r_n\,\up{2}{\Phi}\^r_n \,R_\pm \ \ &,& \ \
 \RR^l_\pm\,\up{1}{\Phi}\^l_n\,\up{2}{\Phi}\^l_n \,=\,
 \up{2}{\Phi}\^l_n\,\up{1}{\Phi}\^l_n \,R_\pm \label{PPn} \\ [1mm]
 \Phi^\a_n  \, \rf_\a \, =\s_\a(\rf_\a) \,\Phi^\a_n \ \ && \ \
 \mbox{ for all }\ \ \rf_\a\in\C^\a \ ,  \a=l,r \  \label{sFn}
\ea
hold with structure data $F_\a, \RR^\a_\pm, \s_\a$ being
identical to the structure data of $\Phi^\a_0$ in eqs.~(\ref{OPE0r}),
(\ref{OPE0l}) and (\ref{RPP0}). In order to get an analogue of eqs.
(\ref{MDM}), we introduce the monodromies
\be{Mn}
 M^\a_n \,:=\,J^\a_{n+1}\dots\,J^\a_{N}\,J^\a_{1}\dots\,J^\a_{n}
 \,=\,(U_n^\a)^{-1}M^\a U_n^\a \ \ \mbox{ for }  \ \ \a=r,l \ \ .
\ee
They are holonomies along the whole circle which begin and end
at the $n^{th}$ site. It is now obvious that
\ba
  v_a^{-1}\, (\Phi^r_n)^{-1}\, D_r \, \Phi^r_n \,=\,M^r_n &,&
  v_a \, (\Phi^l_n)^{-1}\, D_l \, \Phi^l_n \,=\,M^l_n \ \label{MDMn}
\ea
hold for all $0\leq n <N$ and the elements $D_\a\in\G_a\o\C^\a$
are the same as in eqs.~(\ref{MDM}),(\ref{DD}). Let us remark that
the quantum traces $\tr^I_q \t^I(M^\a_n)$
are elements of the algebras $\C^\a \subset \K_N$ from which we
constructed our structure data. Moreover, they do not depend on
the index $n$, i.e., one can prove that $\tr^I_q \t^I(M^\a_n) =
\tr^I_q \t^I(M^\a_m)$ for all pairs $n,m$ \cite{AGS}.

It still remains to investigate the exchange relations of the
vertex operators $\Phi^\a_n$ with currents $J^\a_n \in \G_a \o
\K_N$. Details are explained in Appendix A.5; here we only
state the results:
\ba
 \up{1}{J}\^r_{n+1}\,\up{2}{\Phi}\^r_n \,=\,
 \up{2}{\Phi}\^r_n \, R_+ \, \up{1}{J}\^r_{n+1} \  &,& \
 \up{2}{\Phi}\^r_n \, \up{1}{J}\^r_n \,=\,
 \up{1}{J}\^r_n \, \up{2}{\Phi}\^r_n \, R_- \ , \label{PJn} \\ [1mm]
 \up{2}{J}\^l_{n+1}\,\up{1}{\Phi}\^l_n \, = \,
 \up{1}{\Phi}\^l_n \, R_- \, \up{2}{J}\^l_{n+1} \ &,& \
 \up{1}{\Phi}\^l_n \, \up{2}{J}\^l_{n} \,=\,
 \up{2}{J}\^l_n \, \up{1}{\Phi}\^l_n \, R_+ \ , \label{PJn'} \\ [1mm]
 \up{1}{\Phi}\^r_n\,\up{2}{J}\^r_m \,=\,
 \up{2}{J}\^r_m\,\up{1}{\Phi}\^r_n \ &,& \
 \up{1}{\Phi}\^l_n\,\up{2}{J}\^l_m \,=\,
 \up{2}{J}\^l_m\,\up{1}{\Phi}\^l_n \ \
 \mbox{ for }  m\neq n\, ,n+1 \ , \\ [1mm]
   \up{1}{\Phi}\^r_n \, \up{2}{J}\^l_m \,=\,
 \up{2}{J}\^l_m \,  \up{1}{\Phi}\^r_n \  &,&  \
 \up{1}{\Phi}\^l_n \, \up{2}{J}\^r_m \,=\, \up{2}{J}\^r_m \,
 \up{1}{\Phi}\^l_n \ \ \mbox{ for all }\ n,m \ . \label{RL}
\ea
Finally, as a consequence of these relations and (\ref{W'}) we derive
that
\be{rl}
 \up{1}{\Phi}\^r_n\,\up{2}{\Phi}\^l_m \ =\
  \up{2}{\Phi}\^l_m \, \up{1}{\Phi}\^r_n \ \ \ \ \mbox{ for all }
  \ \ n,m \ \ .
\ee
To summarize, we established that the construction (\ref{Pn})
provides us with chiral vertex operators $\Phi^r_n$ and $\Phi^l_n$
which are naturally assigned to the $n^{th}$ site of the lattice.
These vertex operators share the same structure data $F_r, \RR_\pm^r,
\dots$ and $F_l, \RR_\pm^l, \dots$. Their exchange relations with
elements of the current algebra $\K_N$ are local in the sense
discussed above.

Although the vertex operators have local relations with the
observables, one should expect that they themselves are non-local.
Indeed, it is easy to derive the following exchange relations
(see Appendix A.5):
\ba
  \up{1}{\Phi}\^r_n\,\up{2}{\Phi}\^r_m \,=\,
  \RR^r_- \, \up{2}{\Phi}\^r_m\,\up{1}{\Phi}\^r_n \ &,& \
 \up{2}{\Phi}\^l_n\,\up{1}{\Phi}\^l_m \,=
 \RR^l_+ \, \up{1}{\Phi}\^l_m\,\up{2}{\Phi}\^l_n \ \ \
 {\rm for}\  0\leq n<m<N  \label{braid} \ ,  \\ [1mm]
 \up{1}{\Phi}\^r_n\,\up{2}{\Phi}\^r_m \,=\,
 \RR^r_+ \, \up{2}{\Phi}\^r_m\,\up{1}{\Phi}\^r_n \ &,& \
 \up{2}{\Phi}\^l_n\,\up{1}{\Phi}\^l_m \,=
 \RR^l_- \, \up{1}{\Phi}\^l_m\,\up{2}{\Phi}\^l_n \ \ \
 {\rm for}\  0\leq m<n<N  \label{braid'} \ .
\ea
So, elements $\Phi^\a_n$ and $\Phi^\a_m$ do not commute even if
the $n^{th}$ and $m^{th}$ site at which they are localized are far
apart. The relations (\ref{braid}),(\ref{braid'})
demonstrate clearly that $\RR^\a_\pm$
play the role of braiding matrices in local quantum field theory.

\subsection{Extension on a covering of the circle.}

In Subsection 5.2 we have listed properties of the vertex operators
$\Phi_n^\a$ which are valid for $0\leq n,m<N$. However, unlike the
generators of $\K_N$, the vertex operators live on a covering of
the circle, i.e., if we want to make sense of objects $\Phi^\a_n$
with $n \in {\Bbb Z}$, the operator $\Phi^\a_{n+N}$ necessarily
differs from $\Phi^\a_{n}$. Indeed, $\Phi^\a_n$ may be defined
for $n \in {\Bbb Z}$ by the following difference equation which
is encoded in eqs.~(\ref{Pn}):
\be{dP}
   \Phi^\a_{n+1} \,=\, \Phi^\a_{n} \, J^\a_{n+1} \ .
\ee
Here we assume that $J^\a_n$ has been extended periodically to
$n \in {\Bbb Z}$. Periodicity properties of the objects $\Phi^\a_n$
can be expressed through the monodromies $M^\a_n$ introduced in
(\ref{Mn}),
\be{nN}
 \Phi^\a_{n+kN}=\Phi^\a_n (M^\a_n)^k\ , \ \ \ 0\leq n<N\ ,
 \ \ \ k\in{\Bbb Z}\ .
\ee
To proceed, we observe that properties of $M^\a_n$ are similar to
those of $M^\a\equiv M^\a_0$. Using relations spelled out in
Section 4, we easily find that $M^\a_n$ obey the functorial
relations
\be{DMn}
 \up{2}{M}\^r_n R_+ \up{1}{M}\^r_n \,=\, R_- \D_a (M^r_n)\ , \ \ \
 \up{1}{M}\^l_n R_- \up{2}{M}\^l_n \,=\,
 R_+ \D_a (M^l_n)
\ee
which coincide with (\ref{DM}). Therefore, $M^r_n$ and $M^l_n$ obey
the exchange relations (\ref{MMr}) and (\ref{MMl}). Bearing this in
mind, we employ (\ref{UM}) to derive
\be{PMn}
 \up{2}{\Phi}\^r_{m} \, R_+ \, \up{1}{M}\^r_n \,=\,
 \up{1}{M}\^r_n \, \up{2}{\Phi}\^r_{m} \, R_- \ , \ \ \ \
 \up{1}{\Phi}\^l_{m} \, R_- \, \up{2}{M}\^l_n \,=\,
 \up{2}{M}\^l_n \, \up{1}{\Phi}\^l_{m} \, R_+
\ee
for $0\leq n<N$ and $m=n\,(mod\,N)$, i.e., $m=n+kN$, $ k\in{\Bbb Z}$.

Using the properties of the monodromies $M^\a_n$, we can establish
(see Appendix A.5) that relations (\ref{*Pn})-(\ref{sFn}),
(\ref{MDMn})-(\ref{RL}) are valid for $\Phi^\a_n$ with the coordinate
$n$ being replaced by $n'=n+kN$. Thus, the {\em local} properties of
vertex operators $\Phi^\a_{n+kN}$ living outside of the interval
$0\leq n<N$ coincide with those of $\Phi^\a_{n}$ living inside this
interval.

The extension of the exchange relations between vertex operators
to the covering of our discrete circle is slightly more subtle. For
instance, the braid relation of the vertex operator $\Phi^\a_n$
and its counterpart $\Phi^\a_{n+N}$ does not coincide with (\ref{PPn}).
Instead, we find (see Appendix A.5):
\be{PnN}
 \RR^r_+ \, \up{2}{\Phi}\^r_{n}\,\up{1}{\Phi}\^r_{n+N}\,=\,
 \up{1}{\Phi}\^r_{n+N}\, \up{2}{\Phi}\^r_{n}\, R_- \ , \ \ \ \
 \RR^l_- \, \up{1}{\Phi}\^l_{n}\,\up{2}{\Phi}\^l_{n+N}\,=\,
 \up{2}{\Phi}\^l_{n+N}\, \up{1}{\Phi}\^l_{n}\, R_+ \ .
\ee
A similar situation is found for the braid relations (\ref{braid})-%
(\ref{braid'}). It turns out that here we need to apply eqs.
(\ref{ax3}) for the structure data of the vertex operators. Let us
demonstrate their role by investigating the first eq.~in (\ref{braid})
(i.e, the case $n <m$) with $n$ replaced by $n+N$,
$$
 \ar{c}
 \up{1}{\Phi}\^r_{n+N}\,\up{2}{\Phi}\^r_m \,=\,
 \up{1}{\Phi}\^r_{n}\,\up{1}{M}\^r_n \,\up{2}{\Phi}\^r_m  \,=\,
 \up{1}{v}\^{-1}_a \,{\up{1}{D}}_r\,\up{1}{\Phi}\^r_{n} \,
 \up{2}{\Phi}\^r_m  \,=\,
 \up{1}{v}\^{-1}_a \,{\up{1}{D}}_r\, \RR^r_- \,
 \up{2}{\Phi}\^r_m\,\up{1}{\Phi}\^r_{n}  \,= \\ [1.5mm]
 =\, \up{1}{v}\^{-1}_a \, \RR^r_+ \,\up{2}{\s}(D_r)\,
 \up{2}{\Phi}\^r_m\,\up{1}{\Phi}\^r_{n} \,=\,
 \RR^r_+ \, \up{2}{\Phi}\^r_m\, \up{1}{v}\^{-1}_a \,{\up{1}{D}}_r\,
 \up{1}{\Phi}\^r_{n}  \,=\,
 \RR^r_+ \, \up{2}{\Phi}\^r_m\, \up{1}{\Phi}\^r_{n+N} \ . \er
$$
We see that the result coincides with the first eq.~in (\ref{braid'}),
which is natural since $n+N >m$. Proceeding in the same way, one can
show that the braid relations (\ref{braid}) and (\ref{braid'}) hold,
for all $n,m\in {\Bbb Z}$ such that $|n-m|<N$, $n \neq m$. Thus, the 
equations (\ref{ax3}) became an important ingredient for a self-consistent
extension of the lattice theory beyond the interval $0 \leq n < N$.

\subsection{Construction of the local field $g_n$.}
As we have shown above, the local properties of lattice vertex
operators are the same as those we studied in the toy model case.
Therefore, we can repeat the construction of Subsection 3.4 and
introduce the objects
\be{g}
 g_n \,:=\, \S_a(\Phi^l_n) \, \Phi^r_n \  \in\G_a\o\W_N \  .
\ee
To proceed, we need some more information about the representation
theory of lattice current algebras. As we mentioned before, the
algebras $\K_N$ admit a series of irreducible representations
on spaces $W^{IJ}_N$ where $I,J$ run through classes of irreducible
representations of the quantum algebra $\G$. These spaces $W^{IJ}_N$
are of the form
$$    W^{IJ}_N  \ = \ V^I \o  V^J \o \Re^{\o_{N-1}}\ \ \mbox{ where }
   \ \ \ \Re = \bigoplus_{K} V^{\bar K} \o V^K\ \ . $$
Suppose that we describe $\K_N$ in terms of the holonomies $U^\a_n,
n = 1, \dots, N-1$, the monodromies $M^\a_0$ and the local elements
$N_n, n = 0, \dots, N-1$ (notice that the currents can be
reconstructed from holonomies and monodromies). We divide these
generators into two sets, the first containing all $U_n^\a$ and
$N_m$ for $m \neq 0$ while we put $M^\a_0$ and $N_0$ into the
second set. This choice is made so that objects which were not part
of the toy model are separated from objects we met in Section 3
already. In \cite{AFFS}, an action of $\K_N$ on $W^{IJ}_N$ was
constructed for which objects in the first set, i.e., holonomies
$U^\a_n$ and elements $N_m, m \neq 0$, act trivially on the
factor $V^I \o V^J$ in $W^{IJ}_N$ and irreducibly on
$\Re^{\o_{N-1}}$.

It is then straightforward to see that our algebra $\W_N$ of
vertex operators on the lattice possesses only one irreducible
representation on the total space
$$  \cM_N \ = \ \bigoplus_{I,J} \ W^{IJ}_N \ \cong
    \ \cM \o \Re^{\o_{N-1}} \ \  ,$$
where each summand $W^{IJ}_N$ appears with multiplicity one.
By now, the picture resembles very much the situation in the
toy model: we have the model space $\cM_N$ on which $\W_N$
acts irreducibly. Therefore, we may look for operators that
can be restricted to the diagonal subspace
$$   \H_N \ = \ \bigoplus_J \ W^{\bar J J}_N \ \cong \
     \H \o \Re^{\o_{N-1}}\ \subset \cM_N\ \ \ . $$
This is certainly possible for all elements in $\K_N$. But
in addition, we may restrict the field $g_n$ to $\H_N$. As
in Subsection 3.4, the diagonal subspace is characterized by
the constraint $\rf = \S_{lr}(\rf)$ for all $\rf \in \C^r
\subset \K_N$. If we adjust left and right structure
data according to eqs.~(\ref{rltrafo}) \footnote{This can
be done simultaneously for all sites, since the structure data
do not depend on the lattice site $n$ (see Subsection 5.2).}, the
constraint to $\H_N$ is compatible with the construction
of $g_n$, i.e., (\ref{fSf}) holds with $g$ replaced by $g_n,
n = 0, \dots, N-1$. The properties of the restricted field
are spelled out in the following proposition.

\begin{prop} {\em (Properties of $g_n$)} When restricted to the
diagonal subspace $\H_N$, the element $g_n\in\G_n\o\End(\H_N)$ obeys
the following relations:
\label{prgn}
\ba
\up{2}{g}_n \, \up{1}{g}_n \,=\, \D_a (g_n) \ \ \ & ,&\ \ \
 R_\pm \, \up{2}{g}_n \, \up{1}{g}_n \,=\,
 \up{1}{g}_n \, \up{2}{g}_n \, R_\pm \ , \label{gn} \\ [1mm]
 M^l_n\,g_n \,=\, g_n \, M^r_n \ \ \ & , & \ \ \
 \S_a(g_n) \ = \ g_n^{-1} \ , \\ [1mm]
 g_{n+N} \,=\, g_n \ \ \ &,&\ \ \ \up{1}{g}_n \, \up{2}{g}_m \,=\,
 \up{2}{g}_m\, \up{1}{g}_n\ \ \mbox{ for }\ n\neq m\ ,\label{loc}\\[1mm]
  \up{1}{M}\^r_n\up{2}{g}_n\ R_-\, =\, \up{2}{g}_n\ R_+ \up{1}{M}\^r_n
 \ \ \ &,& \ \ \ \up{1}{M}\^l_n \  R_- \up{2}{g}_n\, = \,
 R_+ \up{2}{g}_n\   \up{1}{M}\^l_n \ , \label{gM}\\[1mm]
 \D_n(\xi)\, g_n \,=\, g_n \, \D'_n(\xi) \ \ \ &,& \ \ \
 \up{1}{N}_{n,\pm} \, R_{\pm} \, \up{2}{g}_n \,=\,
 \up{2}{g}_n \, R_{\pm} \, \up{1}{N}_{n,\pm} \label{ign}
\ea
for all $\xi\in\G$ and $g_n$ commutes with all $\iota_m(\xi) \in \G_n
\subset \K_N$  for $ m \neq n$.
\end{prop}

The properties listed above, and in particular the locality and
periodicity relations (\ref{loc}), allow to regard $g_n$ as an
observable in the lattice WZNW-model. It is a discrete analogue
of the group valued field $g(x)$. Some remarks on the proof of
Proposition \ref{prgn} can be found in Appendix A.6. To complete
the description of $g_n$, let us give its exchange relations with
the chiral currents. Using (\ref{PJn}), we obtain
\ba
 \up{2}{g}_n \, \up{1}{J}\^r_n \,=\,
 \up{1}{J}\^r_n \, \up{2}{g}_n \, R_- \ \ \ &,& \ \ \
 \up{2}{g}_n \, \up{1}{J}\^l_n \,=\,
 \up{1}{J}\^l_n \, R_- \, \up{2}{g}_n   \  ,
 \nn \\ [1mm]
 \up{1}{J}\^r_{n+1}\,\up{2}{g}_n  \,=\,
 \up{2}{g}_n \, R_+ \, \up{1}{J}\^r_{n+1} \ \ \ &,&\ \ \
 \up{1}{J}\^l_{n+1}\,\up{2}{g}_n  \,=\,
 R_+ \, \up{2}{g}_n \, \up{1}{J}\^l_{n+1}  \  , \label{gJ} \\ [1mm]
 \up{1}{J}\^\a_m \, \up{2}{g}_n \,=\, \up{2}{g}_n \, \up{1}{J}\^\a_m \, ,
 \ \a=l,r \ \ \ & &\ \ \  {\mbox{ for } } m\neq n,n+1 \ (mod\, N). \nn
\ea

\subsection{Lattice vertex operators for $\Z_q$.}

Let us construct the algebra $\W_N$ in the case of $\G = \Z_q$. To this
end we have to add the chiral vertex operators introduced in Subsection
3.5 to the lattice $U(1)$-current algebra discussed in Subsection 4.4.
As a result we get the algebra generated by components of the following
elements belonging to $\G_a\o\W_N$:
\ba
  \Phi^\a_n \,=\, \sum_{s=0}^{p-1} \, P^s \o (\wh{Q}^\a_n)^s \,=\,
 e^{\wh{p}\,\o\,\wh{\vs}_n^{\,\a} } \ &,&
 N_{n,\pm} \,=\, \sum_{s=0}^{p-1} \, P^s \o h_n^{\pm s} \,=\,
 q^{\pm \wh{p}\,\o\,\wh{p}_n } \ , \nn \\ [1mm]
  J^r_n \,=\, \sum_{s=0}^{p-1} \, q^{\frac 12 s^2} \, P^s\o (\wh{W}_n^r)^s 
 \,=\, k_a^{-\frac 12}\, e^{\wh{p}\,\o\,\wh\vp_n^r} &,&
  J^l_n \,=\, \sum_{s=0}^{p-1} \, q^{-\frac 12 s^2} \, P^s\o (\wh{W}_n^l)^s 
 \,=\, k_a^{\frac 12}\, e^{\wh{p}\,\o\,\wh\vp_n^l} \ , \nn
\ea
where $\a=r,l$ and $n=0,..,N-1$. According to eqs.~(\ref{nchir})
and (\ref{Pn}), not all the generators are independent. Namely,
the following relations are to be fulfilled:
$$
  \wh{W}_n^l \,=\, h_{n-1}^{-1}\,\wh{W}_n^r\,h_n^{-1} \ , \ \ \ \
  \wh{Q}^\a_n \,=\, \wh{Q}^\a_0 \, \wh{W}_1^\a \dots \wh{W}_n^\a \ .
$$
Due to the Campbell-Hausdorff formula these equalities may be 
re-expressed in terms of the generators 
$\wh\vp_n^\a, \wh{\vs}_n^\a$ as follows:
\be{constr}
 \wh\vp_n^l \,=\, \wh\vp_n^r - \ln q\ (\wh{p}_n + \wh{p}_{n-1}) \ , \ \ \ \
 \wh{\vs}_n^{\,\a} \,=\, \wh{\vs}_0^{\,\a} + \sum_{k=1}^n \,
 \wh\vp_k^\a   \ .
\ee
It is easy to see that all the formulae between $\Phi^\a_n$, $J^\a_n$
and $N_{n,\pm}$ spelled out in Subsections 5.1-5.3 are satisfied if we
add to eqs.~(\ref{hW})-(\ref{WW}) or, alternatively, to eqs.~(\ref{zv})
the following relations:
\ba
 h_n \, \wh{Q}^\a_n \,=\, q\, \wh{Q}^\a_n \, h_n \ , &&
 h_m \, \wh{Q}^\a_n \,=\, \wh{Q}^\a_n \, h_m \ \
 {\rm for\ } m\neq n \ , \nn \\ [1mm]
 \wh{W}_{n+1}^r \, \wh{Q}^r_n \,=\,q\, \wh{Q}^r_n \, \wh{W}_{n+1}^r \ , &&
 \wh{W}_n^r \, \wh{Q}^r_n \,=\, q\, \wh{Q}^r_n \, \wh{W}_n^r \ , \nn \\[1mm]
 \wh{W}_{n+1}^l \, \wh{Q}^l_n \,=\, q^{-1}\, \wh{Q}^l_n \, \wh{W}_{n+1}^l
 \ , && \wh{W}_n^l \, \wh{Q}^l_n \,=\, q^{-1}\, \wh{Q}^l_n \, \wh{W}_n^l
 \ , \nn \\ [1mm]
 \wh{W}^\a_m \, \wh{Q}^\a_n \,=\, \wh{Q}^\a_n \, \wh{W}^\a_m \  &&
 {\rm for\ } m\neq n, n+1 \ , \nn
\ea
which can be rewritten as follows:
\be{os}
 [ \wh{p}_m\,,\, \wh{\vs}_n^{\,\a} ] \,=\, \dl_{m,n} \ , \ \ \
 [ \wh\vp_m^r \,,\, \wh{\vs}_n^{\,r} ] \,=\, -
 [ \wh\vp_m^l \,,\, \wh{\vs}_n^{\,l} ] \,=\, \ln q \,
 ( \dl_{m,n+1}+\dl_{m,n} ) \ .
\ee
Since we already discussed properties of the vertex operators at a
fixed site for the $\Z_q$-theory in Subsections 2.6 and 3.5, we shall
concentrate on the aspects of locality and periodicity here.
Actually, the latter simplify in the case of $\Z_q$ due to
the circumstance that all our monodromies $M^\a_n$ of the same
chirality coincide (since all they are given by (\ref{zM})).
This allows to rewrite eqs.~(\ref{braid})-(\ref{braid'}) and
(\ref{PnN}) in the form (recall that in the case of $\Z_q$ we
have $\RR_\pm=R_\pm=R^{\pm 1}$ with $R$ given in Subsections
2.1 and 2.2):
\be{last33}
  \up{1}{\Phi}\^r_n\,\up{2}{\Phi}\^r_m \,=\,
  R^{\beta(n-m)} \, \up{2}{\Phi}\^r_m\,\up{1}{\Phi}\^r_n \ , \ \ \ \
 \up{2}{\Phi}\^l_n\,\up{1}{\Phi}\^l_m \,=\,
 R^{-\beta(n-m)} \, \up{1}{\Phi}\^l_m\,\up{2}{\Phi}\^l_n \
\ee
for all $n,m\in {\Bbb Z}$. Here $\beta(n-m)=1+2[\frac{n-m}{N}]$
($[x]$ stands for the entire part of $x$) for $n\neq m\,(mod\, N)$ and
$\beta(n-m)=1+[\frac{n-m}{N}]$ for $n=m\,(mod\,N)$. In the
derivation  we have also used the following consequence of 
eqs.~(\ref{os}):
$[\sum_{k=1}^N \wh\vp_k^r \,,\, \wh{\vs}_n^{\,r}]=
 -[\sum_{k=1}^N \wh\vp_k^l \,,\, \wh{\vs}_n^{\,l}]=2\ln q$.
Notice that, since $R^p=e\o e$, the above relations are actually
periodic with a period $N'=pN$ for odd $p$.  That is, the theory
lives on a $p$-fold covering of the circle so that vertex operators
for $\Z_q$ are periodic on a lattice of size $N' = pN$.

Now let us introduce the field $g_n$. We repeat the construction
of Subsection 3.5 and define $g_n$ as follows:
$$
 g_n\,=\,\sum_{s=0}^{p-1} \, P^s \o (\wh{Q}^{\,r}_n)^{\!s} 
 \,(\wh{Q}^{\,l}_n)^{-s}
 \,=\,  e^{\wh{p}\,\o\,(\wh{\vs}^{\;r}_{\,n}-\wh{\vs}^{\;l}_{\,n}) } \ .
$$
It obviously admits restriction to the diagonal subspace
$\H_N$ of the model  space  $\cM_N$ (cf. Subsections 3.5 and 5.4).
The locality of $g_n$ is evident from eqs.~(\ref{last33}) and its
periodicity $g_{n+N}=
(M^l)^{-1} g_n M^r = g_n$ is, in fact, reduced to the Weyl-type
relation which we explained in detail at the end of Subsection 3.5.

\section{AUTOMORPHISMS AND DISCRETE DYNAMICS}
\setcounter{equation}{0}

In this section we shall demonstrate that the lattice theory which we
constructed above indeed may be regarded as a discretization of the WZNW
model. For this purpose we investigate the exchange relations of
currents and some automorphisms of our lattice algebra  in the
classical continuum limit and recover the Poisson structure and
the dynamics of the classical WZNW model, respectively.

\subsection{Remarks on the classical continuum limit.}

Let us briefly discuss the classical continuum limit of the algebra of
vertex operators. Following ideology of \cite{AFS}, we rewrite the exchange
relations (\ref{K4}), (\ref{JJ}) and (\ref{DI}), (\ref{II}) for the chiral
currents in a more compact form:
\ba
 R_{n-m,+}^{-1} \, \up{2}{J}\^r_n \, R_{n-m+1,+} \, \up{1}{J}\^r_m &=&
 \up{1}{J}\^r_m \, R_{n-m-1,-}^{-1} \, \up{2}{J}\^r_n \, R_{n-m,-} \ ,
 \label{nm1} \\ [1mm]
 R_{n-m,+}^{-1} \, \up{1}{J}\^l_n \, R_{n-m+1,-} \, \up{2}{J}\^l_m &=&
 \up{2}{J}\^l_m \, R_{n-m-1,+}^{-1} \, \up{1}{J}\^l_n \, R_{n-m,-} \ ,
 \label{nm2}
\ea
where $R_{n,\pm}:=\dl_{n,0} R_\pm + (1-\dl_{n,0})e\o e$ is, as usual,
an element of $\G_a\o\G_a$. Now we consider these relations in the
limit where $a=2\pi/{\small N}\rar 0$ and $\hbar\rar 0$.

Since for our theory $q=exp\{i\ga\hbar\}$ (cf. Introduction), we
can expand the universal $R$-matrix according to $R_\pm = e\o e +
i\ga\hbar\, r_\pm + O(\hbar^2)$ . On the other hand, the lattice
fields approach their continuum counterparts as $a$ becomes small:
\be{lim}
   J^\a_n \ \rar \ e\o \re - a \, j^\a(x) \ , \ \ \ \
       \Phi^\a_n \ \rar \ \Phi^\a(x) \ , \ \ \ \
   g_n \ \rar g(x) \  ,
\ee
where $x=an$. Bearing in mind that $\frac 1a \dl_{n,0}\rar \dl(x)$
when $a\rar 0$, we obtain the following Poisson brackets from
(\ref{nm1})-(\ref{nm2}):
\ba
 \{\up{1}{j}\^r(x),\up{2}{j}\^r(y)\} &=& {\ns\frac{\ga}{2} }\,
 [C,\up{1}{j}\^r(x)-
 \up{2}{j}\^r(y)]\,\dl(x-y) + \ga \, C \dl'(x-y) \ , \nn \\ [0.5mm]
 \{\up{1}{j}\^l(x),\up{2}{j}\^l(y)\} &=&-{\frac{\ga}{2} } \,
 [C,\up{1}{j}\^l(x)-
 \up{2}{j}\^l(y)]\,\dl(x-y) - \ga \, C \dl'(x-y) \  ,\nn
\ea
where \hbox{$C=(r_+ -r_-)\o \re$}. These are the standard brackets
for the chiral WZNW currents \cite{KnZa,Wit}, and the deformation
parameter $\ga$ is identified with the coupling constant.
\footnote{One may prefer to renormalize the currents by $1/\ga$
so that the $\dl'$-term acquires a coefficient $1/\ga$ which,
in the classical theory, coincides with the level $k$ of the
KM algebra. The quantum correction $1/\ga\rar k+\nu$ is explained,
e.g., in \cite{AFS}.}

The exchange relations (\ref{PJn})-(\ref{PJn'}) can be treated
similarly. Namely, we rewrite them as follows
$$
 \up{1}{J}\^r_{n}\,\up{2}{\Phi}\^r_m \, R_{n-m,-} \,=\;
 \up{2}{\Phi}\^r_m \, R_{n-m-1,+} \, \up{1}{J}\^r_{n} \ , \ \ \ \
 \up{2}{J}\^l_{n}\,\up{1}{\Phi}\^l_m \, R_{n-m,+} \, = \;
 \up{1}{\Phi}\^l_m \, R_{n-m-1,-} \, \up{2}{J}\^l_{n} \ ,
$$
and get the following Poisson brackets for vertex operators in
the classical continuum limit: 
$$
 \{ \up{1}{j}\^r(x),\up{2}{\Phi}\^r(y) \} \,=\, \ga \,
 \up{2}{\Phi}\^r(x) \, C \,\dl(x-y)  \ , \ \ \ \
 \{ \up{1}{j}\^l(x),\up{2}{\Phi}\^l(y) \} \,=\, -\ga \,
 \up{2}{\Phi}\^l(x) \, C \,\dl(x-y)  \ .
$$
These relations are classical  counterparts of the commutation
relations known for the chiral primary fields in the continuum
WZNW model \cite{KnZa}.

Substitution of the expansions $\RR^\a_\pm = e\o e\o\re+i\ga\hbar\,
{\bf r}^\a_\pm + O(\hbar^2)$ into eqs.~(\ref{PPn}) and passing to
the classical continuum limit gives
\footnote{In general, the classical $r$-matrices ${\bf r}^\a_\pm$
keep a non-trivial dependence on variables belonging to $\C^\a$.}
$$
 \{\up{1}{\Phi}\^r(x) ,\up{2}{\Phi}\^r(y)\} \,=\, - \chi^r(x-y) \,
 \up{1}{\Phi}\^r(x)\, \up{2}{\Phi}\^r(y) \ , \ \ \ \
 \{\up{1}{\Phi}\^l(x)\,\up{2}{\Phi}\^l(y)\} \,=\, \chi^l(x-y) \,
 \up{1}{\Phi}\^l(x)\,\up{2}{\Phi}\^l(y) \ ,
$$
where $\chi^\a(x-y)=\ve(x-y)\,\ga\,{\bf r}^\a_+ + \ve(y-x)\,\ga\,
{\bf r}^\a_-$, and $\ve(x)=1$ if $x>0$ and $\ve(x)=0$ if $x<0$.
Such brackets were obtained for the classical WZNW model in
\cite{Fa1,BDF,CGO,Fa2}.

The same technique may finally be applied to the relations
(\ref{gJ}) involving the lattice field $g$ and the resulting
formulae for the classical counterpart of eqs.~(\ref{gJ})
coincide with formulae in \cite{KnZa}, namely,
$$
 \{ \up{1}{j}\^r(x),\up{2}{g}(y) \} \,=\,
 \ga \, \up{2}{g}(x) \, C \,\dl(x-y)\ , \ \ \
 \{ \up{1}{j}\^l(x),\up{2}{g}(y) \} \,=\,
 \ga \, C \, \up{2}{g}(x)\, \dl(x-y)\ .
$$
Thus, in the limit $\hbar \rar 0, a \rar 0$, our main exchange relations
for the chiral currents and chiral vertex operators reproduce the Poisson
structure known for the classical WZNW model.

\subsection{Automorphisms induced by the ribbon element.}

The ribbon element, due to its specific properties, allows to obtain
certain inner automorphisms of the algebra $\W_N$. These are the
subject of the present subsection. \vspace*{3mm}

\noindent
{\it Non-local automorphism induced by global ribbon elements.}\
Consider an automorphism of the form:
\be{aut1}
 A \; \mapsto \; \rv_r^{-1} \rv_{l}^{\phantom{-1}} \!\!\!  A \;
 \rv_r \rv_l^{-1} \ , \ \ \
{\mbox{ for all}}\ \ A\in \W_N \ .
\ee
Here $\rv_r\in\C^r$ and $\rv_{l}^{\phantom{-1}}\!\!\!\!\in\C^l$. 
We call the ribbon elements
$\rv_\a$ {\em global} because they are constructed from
the monodromies $M^\a$, which are non-local.

Since the subalgebras $\C^\a$ constitute the center of $\K_N$,
all the elements of the current algebra $\K_N \subset \W_N$ are
invariant under the transformation (\ref{aut1}). For the
vertex operators this transformation is nontrivial and may
be rewritten with the help of (\ref{DD}), (\ref{sFn}) and
(\ref{MDMn}) so that it becomes
$$
 \Phi_{n+kN}^\a \ \mapsto \ \Phi_{n+kN}^\a \, M_n^\a \,=\,
 \Phi_{n+(k+1)N}^\a  \ ,  \ \ \ 0\leq n<N, \ \ k\in {\Bbb Z} \ .
$$
Thus, the automorphism (\ref{aut1}) is non-local, i.e., it corresponds
to a shift $n\mapsto n+N$ or, in other words, it rotates the lattice
by angle $2\pi$. Being restricted on the diagonal subspace, the field
$g_n$ is periodic (see Proposition \ref{prgn}), and hence it is
invariant under the transformation (\ref{aut1}). In this sense, the
automorphism (\ref{aut1}) separates ``physical'' variables living on
the circle from ``non-physical'' ones (like the vertex operators)
living on a covering of the circle. \vspace*{3mm}

\noindent
{\it Local automorphism induced by  local ribbon elements.}\
Recall that the $n^{th}$ site of the lattice is supplied with a copy
$\G_n$ of the symmetry algebra $\G$. Therefore we can use the {\em
local ribbon elements} $v_n\in\C_n\subset\G_n$ to construct the
following transformation:
\be{aut2}
 A \; \mapsto \; v_0 v_1 \dots v_{N-1} \, A \,
 (v_0 v_1 \dots v_{N-1})^{-1} \ , \ \ \ \mbox{ for any }\ \ A\in \W_N \ .
\ee
Here the product is taken over all sites of the lattice.

To obtain more explicit formulae for the automorphism (\ref{aut2}),
we have to use relations (\ref{K1}), (\ref{FI}), (\ref{xPn}), employ
equation (\ref{rib}) and remember that $v_n$ belongs to the center
of $\G_n$. As a result we get
\be{a2'}
 \ar{ll}
 J^r_n \ \mapsto \ N_{n-1} \, J^r_n \, N_{n}^{-1}  \ , &
 J^l_n \ \mapsto \ \wt{N}_{n-1}^{-1} \, J^l_n \, \wt{N}_{n}\ ,\\ [1.5mm]
 \Phi^r_n \ \mapsto \ v_a \, \Phi^r_n \, N_{n}^{-1}  \ ,&
 \Phi^l_n \ \mapsto \ v_a \, \Phi^l_n \, \wt{N}_{n} \ , \er
\ee
where we used notations of Subsection 2.2, i.e., $N=N_+N_-^{-1}$,
$\wt{N}=N_+^{-1}N_-$. The elements of $\G_n$ and $\C^\a$ remain
invariant under (\ref{aut2}), in particular, $N_{n,\pm} \mapsto
N_{n,\pm}$. With the help of these explicit expressions we also
obtain
\be{a2''}
 M^r_n \ \mapsto \ N_{n} \, M^r_n \, N_{n}^{-1}  \ , \ \ \
 M^l_n \ \mapsto \ \wt{N}_{n}^{-1} \, M^l_n \, \wt{N}_{n} \ ,\ \ \
 g_n \ \mapsto \ \wt{N}_{n} \, g_n \, N_{n}^{-1}  \ .
\ee
We know already that these formulae describe an automorphism of
the algebra $\W_N$ because they were obtained by conjugation with
a unitary element, namely the product of local ribbon elements, in
formula (\ref{aut2}). Without this knowledge, it would be a quite
non-trivial task to check the automorphism property directly
for the expressions in (\ref{a2'})-(\ref{a2''}). To do this,
one would need to apply the relations (\ref{JN}) and (\ref{NPn})
many times.  \vspace*{3mm}

\noindent
{\it Local automorphism induced by $\k_n$.}\
To construct one more inner automorphism of $\W_N$ we employ
the square roots of the local ribbon elements,
\be{aut3}
 A \; \mapsto \; \k_0 \k_1 \dots \k_{N-1} \, A \,
 (\k_0 \k_1 \dots \k_{N-1})^{-1} \ ,
  \ \ \  \mbox{ for all }\ A\in \W_N \ .
\ee
Here $\k_n \in \C_n\subset \G_n$, $\k_n^2=v_n$
and the product is taken over all sites.

Computations similar to those performed above (and making use of
(\ref{SS})) allow to rewrite the transformation (\ref{aut3}) in the
following explicit form:
$$
 \ar{rcll}
 J^r_n & \mapsto & N_{n-1,+} (S_{n-1})^{-1} \, J^r_n\, S_{n}
 N_{n,+}^{-1}&=\,
 N_{n-1,+}  \, ((J^r_n)^*)^{-1} \, N_{n,+}^{-1}  \nn \ , \\ [1.5mm]
 J^l_n & \mapsto & N_{n-1,-}^{-1} S_{n-1} \,J^l_n \, (S_{n})^{-1}
 N_{n,-} &=\,
 N_{n-1,-}^{-1} \, ((J^l_n)^*)^{-1} \, N_{n,-}  \nn \ ,\\ [1.5mm]
 \Phi^r_n & \mapsto & \k_a \, \Phi^r_n \, S_n N_{n,+}^{-1}
 &=\, \k_a ((\Phi^r_n)^*)^{-1} \, N_{n,+}^{-1} \nn \ ,\\ [1.5mm]
 \Phi^l_n & \mapsto & \k_a \, \Phi^l_n \, (S_n)^{-1} N_{n,-}
 &=\, \k_a ((\Phi^l_n)^*)^{-1} \, N_{n,-} \nn \ .   \er
$$
$$
 \ar{c} M^r_n \ \mapsto \ N_{n,+}\,((M^r_n)^*)^{-1}
 \, N_{n,+}^{-1} \ , \ \
 \ M^l_n \ \mapsto \ N_{n,-}^{-1} \,
 ((M^l_n)^*)^{-1}\,N_{n,-}\  \er
$$
and $N_{n,\pm} \mapsto N_{n,\pm}$, as before. Here the r.h.s. of all
formulae have been rewritten with the help of the $*$-operation introduced
for elements of $\W_N$ in Sections 4 and 5. Having done so, we see that
the image of all basic objects $X \in\G_a \o \W_N$ under the automorphism
(\ref{aut3}) coincides with $(X^*)^{-1}$ up to a multiplication with
factors $N_{n,\pm}$. We can now accept an inverse logic -- we may say
that the automorphism (\ref{aut3}) together with the rules $(N_\pm)^*\,
=\,N_\mp$ defines the $*$-operation on $\W_N$. This picture reveals the
naturalness of our $*$-operation, which might have appeared somewhat
artificial in the previous sections. It also makes the role of the ribbon
element in our theory even more remarkable.

To conclude this discussion, we would like to mention that for a lattice
of even length, i.e., for $N=0\,(mod\,2)$, one may also consider
automorphisms of $\W_N$ generated by the alternating products of
$v_n^{\pm 1}$ or $\k_n^{\pm 1}$.

\subsection{Discrete dynamics.}

As we saw above, the exchange relations of the algebra $\W_N$
allowed to recover the Poisson structure of the classical WZNW-model
in the classical continuum limit. However, this is certainly not
sufficient for a construction of the lattice WZNW model.
Indeed, the complete description of a classical theory involves
an evolution equation for the dynamical variables in addition
to the specification of the Poisson structure. Similarly, the
formulation of a discrete quantum model requires not only a set
of exchange relations between quantum operators but also some
one parameter family of automorphisms of the algebra generated
by operators in the quantum theory. The parameter is interpreted
as time variable. For a theory on a discrete space it is natural
to discretize the time as well so that the parameter essentially
runs through the set of integers only. In this case the whole
family of automorphisms can be reconstructed from the automorphism
which provides the evolution for an elementary step in time. Such
an automorphism of a lattice model must be local, i.e., the result
of its action on the variables assigned to a given site (or link)
can only involve variables assigned to some neighboring sites (or
links). In the previous subsection we considered three automorphisms
of the algebra $\W_N$. The first of them was non-local and hence did
not correspond to any dynamics.\footnote{However, one can use it to
describe dynamics of the toy model (see \cite{AF2}).} The second
and the third automorphism were local and, in principle, one may
use them in constructing the corresponding classical continuum
models. However, the dynamics of such models do not reproduce the
dynamics of the WZNW theory.

In this subsection we are going to consider local automorphisms
which can be interpreted as dynamics of the discrete WZNW model.
Let us recall that in the continuum WZNW model the equation of
motion for the $G$-valued field $g(x)$ takes the form:
\be{cg}
   \partial_+ \partial_- \, g \,=\,
   (\partial_+ g) \, g^{-1} \, (\partial_- g) \ ,
\ee
where $\partial_\pm= \frac 12 (\partial_0 \pm \partial_x)$. {}From
the field $g(x)$ one may construct the following Lie algebra valued
currents
\be{cur}
 j^r \,=\, g^{-1} \partial_- g \ \ \  , \ \ \
 j^l=(\partial_+ g) g^{-1}\ \ .
\ee
They turn out to be chiral objects in the sense that their equations
of motion are trivialized:
\be{cj}
   \partial_+ \, j^r \,=\, \partial_- \, j^l \,=\, 0 \ .
\ee
In the Hamiltonian approach, the initial data are provided by the
values of $g(0)$, $j^r(x)$ and $j^l(x)$ at time $t=0$. To recover
the dynamics of $g(x)$ one solves the equations
\be{cd}
  \partial_0 \, g \,=\, j^l \, g + g \, j^r \ , \ \ \
  \partial_x \, g \,=\, j^l \, g - g \, j^r \  .
\ee
Eqs.~(\ref{cj}) and (\ref{cd}) can be derived with the help of
the Poisson brackets given in Subsection 6.1 if the Hamiltonian
and the total momentum are chosen as follows:
\be{HP}
 H \,=\, \frac 1{2\ga} \, 
    \int \tr \left[ (j^r(x))^2 + (j^l(x))^2 \right] \, dx \ , \ \ \ \ 
 P \,=\, \frac 1{2\ga} \, 
    \int \tr \left[ (j^r(x))^2 - (j^l(x))^2 \right] \, dx \ ,
\ee
where the integration is taken over the whole circle and $\tr$ 
is the usual trace in the corresponding Lie algebra.

Let us develop an analogue of the given picture in the quantum lattice
theory. More precisely, we shall consider the ``physical'' subalgebra
$\P_N$ of $\W_N$ generated  by components of the chiral currents $J^\a_n$
and the field $g_n$, $n=1,..,N$ which are subject to the relations
spelled out in Sections 4 and 5. As we have indicated in our general
discussion above, it is natural to work with a discrete time with a
minimal time interval $\tau$ (see also \cite{Fa3,FV1}), so that the
evolution of the quantum theory is described by a single automorphism
of $\P_N$. In addition to this, we shall also introduce an automorphism
which is responsible for the shifts by one lattice unit $a=2\pi/N$ in
space.

\begin{lemma} {\em (Shift and evolution automorphisms)} \label{dyn}
Let $\P_N$ denote the algebra generated by components of $J^\a_n$
and $g_n$ (restriction to the diagonal subspace $\H_N$ is understood).
Then the following two transformations
$T_V, T_U$,
\ba
 T_V(J^\a_n) &=&  J^\a_{n+1}\ , \ \ \a=l,r \ , \label{sJ}  \\ [1mm]
 T_V(g_n) &=& (J^l_{n+1})^{-1}\, g_n \,J^r_{n+1}\ \label{sg} \ ,
\ea
and
\ba
 T_U(J^r_n) &=& J^r_{n-1}\ , \ \ \
 T_U(J^l_n) \ = \ J^l_{n+1}\ , \label{eJ} \\ [1mm]
 T_U(g_n) &=& (J^l_{n+1})^{-1}\, g_n \, (J^r_n)^{-1}  \label{eg}
\ea
extend to automorphisms of the algebra $\P_N$. We call $T_V$ the
{\em shift automorphism} and $T_U$ the {\em evolution automorphism}
of the lattice WZNW-model.
\end{lemma}

It is straightforward to verify that the transformations $T_V$
and $T_U$ preserve all the relations for $J^\a_n$ and $g_n$
given above. Notice also that one may extend $T_U, T_V$ to the
whole algebra $\W_N$ by eqs.~(\ref{sJ}), (\ref{eJ}) and, in
addition, the formulae
\ba
   T_V(\Phi^\a_n) & =& \Phi^\a_n J^\a_{n+1} \ \ , \nn \\[1mm]
   T_U(\Phi^r_n)\ = \ \Phi^r_n \ (J^r_n)^{-1} \ \ \ \ & , & \ \ \ \
   T_U(\Phi^l_n)\ =\ \Phi^l_n \ J^l_{n+1} \nn \ \ .
\ea
These automorphisms are actually combined of two chiral automorphisms
(cf. Subsection 6.4).
After restriction to the diagonal subspace and, hence, to the algebra
$\P_N$, we recover eqs.~(\ref{sg}) and (\ref{eg}).
\footnote{This needs the following variant of eq.~(\ref{dP}):
$\S_a(\Phi^l_{n+1})=(J^l_{n+1})^{-1} \S_a(\Phi^l_{n})$; it is evident
if we take relation (\ref{theta}) into account.}

Assume now that Lemma \ref{dyn} describes inner automorphisms of $\P_N$.
That is, suppose that there exist operators ${\rm V}, {\rm U} \in \P_N$
such that
$$
 T_V(A) \,=\, {\rm V} \, A \, {\rm V}^{-1} \ \ \mbox{ and }  \ \ \
 T_U(A) \,=\, {\rm U} \, A \, {\rm U}^{-1}  \ \ \ \mbox{ for any } \
 A\in\P_N \ .
$$
$\rm V$ and $\rm U$ are usually called {\em shift} and {\em evolution}
operators, respectively. In the classical continuum limit $a\rar 0$,
$\tau\rar 0$ they reproduce the momentum and the Hamiltonian (\ref{HP}):
${\rm V}\rar e+\frac i\hbar a\,P$, ${\rm U} \rar e+\frac i\hbar \tau\,H$.

Our interpretation of the transformations in Lemma \ref{dyn}
as discrete shifts in space and time, motivate to introduce
the the objects $J^\a_n(t), g_n(t) \in \G_a \o \P_N$ such that
\be{dJ}
 \hspace*{-5mm} J^r_n(t+\tau) \,=\, \ J^r_{n-1}(t) \ , \ \ \
 J^l_n(t+\tau) \,=\, J^l_{n+1}(t) \ ,
\ee
\be{dg}
 g_n(t+\tau) \,=\, (J^l_{n+1}(t))^{-1} \, g_n(t) \, (J^r_{n}(t))^{-1} \ ,
 \ \ \
 g_{n+1}(t) \,=\, (J^l_{n+1}(t))^{-1} \, g_n(t) \, J^r_{n+1}(t)\
\ee
and $J^\a_n(0), g_n(0)$ coincide with our usual generators $J^\a_n,
g_n$, respectively. These expressions define $J^\a_n(t)$ and
$g_n(t)$ for $t = k \tau$ with $k$ being integer. Below we shall
also need the equations inverse to (\ref{dJ})-(\ref{dg}):
\ba
 J^r_n(t-\tau) \,=\,  J^r_{n+1}(t) \ ,\ &&\
 J^l_n(t-\tau) \,=\,  J^l_{n-1}(t) \ , \label{-dJ} \\ [1mm]
 g_n(t-\tau) \,=\, J^l_{n}(t) \, g_n(t) \, J^r_{n+1}(t) \ , \ && \
 g_{n-1}(t) \,=\, J^l_{n}(t) \, g_n(t) \, (J^r_{n}(t))^{-1} \ .
 \label{-dg}
\ea
We can now use the rules $J^\a_n(t)  \rar e\o\re - a j^\a(x,t)$ and
$g_n(t) \rar g(x,t)$ from Subsection 6.1 (here $x=na=2\pi n/N$ as before)
to establish that in the classical continuum limit eqs.~(\ref{dJ})-%
(\ref{dg}) become precisely the equations (\ref{cj})-(\ref{cd})
\footnote{In the continuum limit the quotient $c:= a/\tau$ 
(speed of light) is supposed to be fixed.
In fact, eq.~(\ref{cg}) implies that we put $c=1$.}. {}Further,
combining (\ref{dg}) and (\ref{-dg}), we obtain the following relations:
\be{CLR}
 (g_{n-1}(t))^{-1} \, g_n(t-\tau) \,=\, J^r_n \, J^r_{n+1} \ ,\ \ \ \
 g_{n-1}(t)\,(g_n(t+\tau))^{-1} \,=\, J^l_{n+1} \, J^l_{n} \ .
\ee
These are lattice analogues of the definitions (\ref{cur}) of the chiral
currents. Notice that we can express only products of currents on
neighboring links through the field $g_n(t)$ (nevertheless, in the
classical continuum limit, eq.~(\ref{cur}) is certainly recovered).
\footnote{Formally, we can split eqs.~(\ref{CLR}) into the following
relations: $J^r_n=(g_{n-1}(t))^{-1} g_{n-\frac 12}(t-\frac 12 \tau)$
and $J^r_{n+1}=(g_{n-\frac 12}(t-\frac 12 \tau))^{-1}\,g_n(t-\tau)$.
However, the variables assigned to half integer sites or times are not
defined in the lattice formalism. To avoid this we could consider these
relations as relations in $\W_N$ and re-express the involved $g$-fields
through vertex operators, while using that vertex operators are chiral
to replace formal variables on half integer space-time points by
true objects of the lattice theory. As a result we would get the
obvious relations $J^\a_n(t)=(\Phi^\a_{n-1}(t))^{-1}\Phi^\a_{n}(t)$.}

Eqs.~(\ref{CLR}) allow to obtain the dynamics of the lattice model
in terms of the field $g_n$ only. Indeed, since their r.h.s. are
manifestly chiral objects, the combinations of $g_n$-variables on
the l.h.s. are to be invariant under the substitutions $t\rar t+\tau$,
$n\rar n+1$ and $t\rar t+\tau$, $n\rar n-1$, respectively. Thus, we
derive a lattice analogue of the equation of motion (\ref{cg}):
\be{LG}
 g_{n+1}(t) \, (g_{n}(t-\tau))^{-1} \,=\,
 g_{n}(t+\tau) \, (g_{n-1}(t))^{-1} \ .
\ee
Being a discrete analogue of an equation of second order in both
variables, this relation involves four different points on the
space-time lattice (see Figure 5). A natural choice of the initial data
for eq.~(\ref{LG}) is provided by the set $g_{n}(t-\tau)$ and
$g_{n}(t)$, \hbox{$n=0,..,N-1$} (here $t$ is fixed). It is interesting
to notice that this set is divided into two subsets (black and white
circles on Figure 5) which have an independent evolution;\footnote{Let
us stress that it is not necessary to impose a continuity condition on
the initial data, i.e., to demand that they possess smooth continuum
limit. Moreover, it seems interesting to study the case when the two
independent subsets of initial data have different continuum limits
(cf. also \cite{Fa3}).}
that is, the solution constructed according to eq.~(\ref{LG}) from one 
of the sets never interacts with that constructed from the other set.
According to eqs.~(\ref{CLR}), the initial data $g_{n}(t-\tau)$,
$g_{n}(t)$ (at fixed time $t$) can be restored if we are given the set
of currents $J^\a_n$, $n=1,..,N$ and two values of the $g$-field taken 
at two arbitrary points of the independent subsets, e.g. $g_0(0)$ and 
$g_0(-\tau)$. This is a lattice analogue of the initial data usually
used in the continuum Hamiltonian approach (see above).

\hbox{\begin{picture}(300,75)
    \put(92,30){\circle*{3}} \put(80,34){$\large x-a$}
    \put(148,30){\circle*{3}} \put(136,34){$\large x+a$}
    \put(120,15){\circle*{3}} \put(123,11){$\large t-\tau$}
    \put(120,45){\circle*{3}} \put(123,45){$\large t+\tau$}
    \put(60,30){\line(1,0){120}} \put(120,5){\line(0,1){50}}
    \put(262,30){\circle*{3}} \put(250,34){$\large x-a$}
    \put(318,30){\circle*{3}} \put(306,34){$\large x+a$}
    \put(374,30){\circle*{3}} \put(360,34){$\large x+2a$}
    \put(290,15){\circle*{3}} \put(293,11){$\large t-\tau$}
    \put(234,15){\circle*{3}}
    \put(346,15){\circle*{3}}
    \put(220,30){\line(1,0){170}} \put(290,5){\line(0,1){50}}
    \put(234,30){\circle{2.5}} \put(262,15){\circle{2.5}}
    \put(290,30){\circle{2.5}} \put(318,15){\circle{2.5}}
    \put(346,30){\circle{2.5}} \put(374,15){\circle{2.5}}
\end{picture} }
\begin{center}
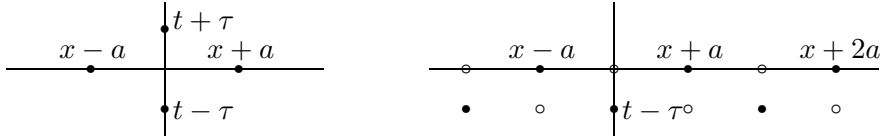
\parbox{13cm}{ \small {\bf Figure 5:} Graphical presentation of the
discrete equation (\ref{LG}) and a possible choice of the initial data.
The two subsets of the initial data have independent evolutions. }
\end{center}
\vspace*{5mm}

To summarize, in this section we have demonstrated that the elements
$J^\a_n$ and $g_n$ which constitute the ``physical'' variables in the
algebra $\W_N$ are indeed quantum lattice analogues of the chiral
currents and the group valued field in the WZNW model.

\subsection{$U(1)$-WZNW model.}

We conclude this section with some comments on the $\Z_q$-case. 
Recall that the $g$-field constructed from lattice vertex operators
in the case of $\G=\Z_q$ is given by (cf. Subsection 5.5) $g_n =
e^{\wh{p}\,\o\,\phi_n} $, where $\phi_{n}=\wh{\vs}^{\;r}_{\,n}-
\wh{\vs}^{\;l}_{\,n}$ is an operator acting on the physical space
$\H_N$ (see Subsections 3.5 and 5.5). In the classical limit
$\phi_{n}(t)$ becomes a lattice variable which, according to
(\ref{LG}), obeys to the following equation of motion:
\be{ff}
 \phi_n(t+\tau) + \phi_n(t-\tau) \,=\, \phi_{n+1}(t) + \phi_{n-1}(t) \ .
\ee
This relation discretizes the equation of motion $\partial_+ \partial_-
\phi(x,t)=0$ of a free field. The latter is known to arise, in particular,
for the continuum $U(1)$-WZNW model. Moreover, in the classical continuum
limit the standard Poisson structure of the abelian WZNW theory is easily
recovered from the exchange relations of our $\Z_q$ lattice model. These
two observations allow to identify the $\Z_q$ lattice theory as a quantized
lattice $U(1)$-WZNW model. In spite of its simplicity, the $U(1)$-theory
has a lot of structure in common with the more complicated nonabelian
models. In fact, the abelian model was used here to illustrate many
elements of our general theory.

It is also worth mentioning that the abelian lattice theory itself has
non-trivial mathematical aspects. In particular, explicit formulae for
shift operators in chiral theories have been worked out in \cite{FV1,
Fa3,AlRe}. We may use these results to present expressions for the
the shift and evolution operators $V$ and $U$.  The latter can be
decomposed into the chiral components: $V=V_l V_r$ and $U=V_l V_r^{-1}$.
When acting on elements of the
algebra $\W_N$, the operators $V_\a \in\W_N^\a$ generate shifts
for the chiral sectors, i.e.,
\be{shifts}
  (e\o V_\a) \, J^\a_n \, (e\o V_\a)^{-1} \,=\, J^\a_{n+1} \ , \ \ \ \
  (e\o V_\a) \, \Phi^\a_n \, (e\o V_\a)^{-1} \,=\, \Phi^\a_{n+1} \ ,
\ee
where $\a=r,l$, $n\in {\Bbb Z}$, and $e\o V_l$, $e\o V_r$ commute with
any  element from  $\G_a\o\W_N^r$ and $\G_a\o\W_N^l$, respectively.

\begin{prop} \label{lapr} {\em (Shifts operators for $\Z_q$)}
Let $\W_N$ be the algebra of lattice vertex operators as defined in 
Subsection 5.5, i.e., it is generated by the elements $\wh{W}_n^\a=
e^{\wh\vp_n^\a}$ (chiral currents) and $\wh{Q}^\a_n=e^{\wh{\vs}_n^{
\,\a}}$ (vertex operators) obeying the relations spelled out in
Subsections 4.4 and 5.5. Let the lattice length $N$ be odd. Then
the chiral shift operators obeying (\ref{shifts}) are given by
\be{cshi}
 V_\a \,=\, Z_\a\, \prod_{k=1}^{N-1} \, \rho_\a({\small N-k}) \ , 
\ee
where $\rho_r(k)=\exp\{-\frac{1}{2\ln q}\, (\wh\vp_{k}^r)^2\}$ and \,
$\rho_l(k)= \exp\{\frac{1}{2\ln q}\, (\wh\vp_{k}^l)^2\}$. The function
$Z_\a$ depends only on the element $C_\a=(\prod_{k=1}^{\frac{N+1}{2}}
\wh{W}_{2k-1}^\a)\,(\prod_{k=1}^{\frac{N-1}{2}} 
\wh{W}_{2k}^\a)^{-1} \in \K_N^\a$.
\end{prop}

To verify that $V_r$ and $V_l$ obey eqs.~(\ref{shifts}) one
proceeds in two steps. The first of them concerns the relations between
$V_\a$ and the chiral currents and it has been performed in
\cite{AlRe}. The computation is based on the following relations:
\footnote{Such relations were used first in
\cite{FV1,Fa3} in the construction of shift operators for $U(1)$-current
algebra for even $N$ (cf. remarks in the text). }
$$
 e^{\wh\vp_n^\a} \, \rho_\a(n+1) \, \rho_\a(n) \,=\,
  \rho_\a(n+1) \, \rho_\a(n) \, e^{\wh\vp_{n+1}^\a} \ , 
  \ \ n=1,..,N-2 \ ,
$$
which hold due to eqs.~(\ref{zv}). The remaining relations for
$n=N-1,N$ are consequences of the others if the function $Z_\a$ is
chosen in a specific way (see \cite{AlRe} for details).

In the second step, one checks the desired properties of $V_\a$ with
respect to chiral vertex operators. To this end we derive the following
relation:
\be{shv}
 \!\! e^{\wh{\vs}_n^{\,\a}} \, \rho_\a(n+1) \,=\,  
 \rho_\a(n+1) \, e^{-\frac 12 \ln q} \, e^{\wh\vp_{n+1}^\a} \,
 e^{\wh{\vs}_n^{\,\a}} \,=\,
 \rho_\a(n+1) \, e^{\wh{\vs}_n^{\,\a} + \wh\vp_{n+1}^\a} \,=\,
 \rho_\a(n+1) \, e^{\wh{\vs}_{n+1}^{\,\a}} \, .
\ee
{}For the first equality we used the commutation relations (\ref{os}).
After this, the Campbell-Hausdorff formula was employed before we
could insert eqs.~(\ref{constr}). Notice that it suffices to prove
eq.~(\ref{shv}) for $n=1$. Due to (\ref{constr}), the relations between
$V_\a$ and the vertex operators assigned to other sites
are consequences of this case (as soon as the relations for $V_\a$
with chiral currents are established). This completes the proof.

Let us comment on the construction of the chiral shift operators for
a lattice of even length $N$ suggested in \cite{FV1,Fa3}. In this case
the shift operators are also
given by (\ref{cshi}) but without the factor $Z_\a$. When checking the
relations between these operators $V_\a$ and chiral currents, Faddeev
and Volkov had to assume that $\prod_{k=1}^{N/2} \wh{W}_{2k-1}^\a \,=\,
\prod_{k=1}^{N/2} \wh{W}_{2k}^\a$. Unfortunately, such a constraint is
incompatible with the exchange relations in the full theory which
includes the objects $N_n$ in addition to chiral currents.
One way to bypass this problem would involve shifts by two lattice
units.

Let us finally mention that the function $\rho$ (which can be
identified as a $\theta$-function, if written in terms of
$\wh{W}^\a_n$) appearing in (\ref{cshi}) admits factorization
into a product of two functions of a q-dilogarithm type
(see \cite{Fa3}). Actually, these objects (the $\theta$-function
and the q-dilogarithm) turn out to be quite universal building
blocks for shift operators. They were employed in the recent
work \cite{FV2} to construct shift operators for the $SU(2)$-%
lattice KM algebra. Since the expressions involving
$\theta$-functions and q-dilogarithms resemble those used
in the abelian theory, one expects that the new operators
of \cite{FV2} serve as shifts not only for the current algebras
but for the whole algebras of vertex operators as well.

\renewcommand{\thesection}{}
\section{CONCLUSION}

In the present paper we have described the construction of lattice
vertex operators for a given modular Hopf algebra. The
investigation of the classical continuum limit reveals a
clear relation between the lattice algebras and the
WZNW-model. Since the latter can be reduced to the affine
Toda model, our technique may be applied to this theory
as well (with certain modifications). Furthermore, there
exist many connections with Chern-Simons theory in $2+1$
dimensions (see \cite{AGS} for lattice constructions
of Chern-Simons observables) which motivate to extend
our framework to two spatial dimensions. 

Let us briefly list some aspects of the presented theory which
have not been developed. As we mentioned before, formulae for
vertex operators are known only for some particular cases. It
would be interesting to work out explicit presentations for
universal vertex operators $\Phi$ of the deformed universal
enveloping algebras $U_q(\sg)$. Alternatively, one may try to
find universal structure data $F$, $\s$ which solve the
discussed set of equations.
A further natural extension is to incorporate infinite dimensional 
structures such as the deformed affine algebras. This might 
allow for a comparison with the approaches in \cite{FrRe,JiMi}
(see also references therein).

Another problem which is to be solved to complete the description
of the quantum lattice WZNW model is an explicit construction of
the shift and evolution operators. By now, exact formulae have been
found for the cases of $U(1)$ and $SU(2)$ \cite{FV1,Fa3,AlRe,FV2}.
These examples, however, hint at some uniform structures
(such as the appearance of q-dilogarithms) that might lead to
new formulae for shift operators in more general theories.

\vspace*{2mm}
\noindent
{\bf Acknowledgments}: We would like to thank A.Yu.Alekseev,
L.D.Faddeev, A.Fring, J.Fr\"ohlich, P.P.Kulish, F.Nill, A.Yu.Volkov 
for useful discussions and E.Jagunova for preparing the pictures. A.B.
is grateful to Prof. R.Schrader for the hospitality
at the Institut f\"ur Theoretische Physik, Freie Universit\"at
Berlin, and to the T\"opfer Stiftung for financial support.

\section{APPENDIX: SOME PROOFS AND FURTHER RELATIONS}
\setcounter{section}{1} \setcounter{subsection}{0}
\setcounter{equation}{0}
\renewcommand{\thesection}{\Alph{section}}
\renewcommand{\thesubsection}{\thesection.\arabic{subsection}}
\renewcommand{\theequation}{\thesection.\arabic{equation}}

\subsection{Proof of Proposition 1.}
It is not difficult to obtain the relations stated in Proposition
\ref{strucprop} from the defining properties of $\Phi$. For instance,
the formula (\ref{Dop}) for $D$  follows from the definition
(\ref{str3}) when $N = RR'$ is re-expressed in terms of the
ribbon element according to eq.~(\ref{rib}). To derive equation
(\ref{ax1}) one needs no more than associativity of the
multiplication in $\V$ along with co-associativity of the
co-product $\D_a$ on $\G_a$,
$$
 \ar{c} \up{3}{\s}(F_{12})\ (\D_a \o id) (F) \ (\D_a \o id)
  \D_a(\Phi) \,=\,  \up{3}{\s}(F_{12}) \ \up{3}{\Phi} \
  \D_a(\Phi)_{12}  \,=\,  \up{3}{\Phi}\ (\up{2}{\Phi}\ \up{1}{\Phi}) \,=
  \\ [2mm]  =\, (\up{3}{\Phi}\ \up{2}{\Phi})\ \up{1}{\Phi} \,=\,
  F_{23}\  \D_a(\Phi)_{23}\ \up{1}{\Phi}
  \,=\,   F_{23} \ (id \o \D_a) (F) \ (id \o \D_a) \D_a(\Phi) \ . \er
$$
The first relation of eqs.~(\ref{ax3}) is a consequence of the
covariance property (\ref{Ncov}) of $\Phi$,
$$
 \ar{c} \up{1}{D} \  \RR_- \up{2}{\Phi} \ \up{1}{\Phi}
   \; = \; \up{1}{\Phi} \  \up{1}{N} \  \up{1}{\Phi}\^{-1}
   \RR_- \up{2}{\Phi} \ \up{1}{\Phi}
   \; = \;  \up{1}{\Phi} \  \up{1}{N} \  \up{2}{\Phi} R_-
   \, = \;  \up{1}{\Phi} \  \up{2}{\Phi} \ R_+  \up{1}{N} \,= \\[2mm]
    =\,   \RR_+ \up{2}{\Phi} \  \up{1}{\Phi} \ \up{1}{N}
   \, = \,  \RR_+ \up{2}{\Phi} \  \up{1}{D} \ \up{1}{\Phi}
   \; = \,  \RR_+ \up{2}{\s}(D) \ \up{2}{\Phi} \ \up{1}{\Phi} \ . \er
$$
We have inserted the definition  (\ref{str3}) twice and
used commutation relations (\ref{PhPh}).

Next, using definitions (\ref{str1})-(\ref{str2}) of the structure data,
we easily check (\ref{ax2}):
$$
 \ar{c} \up{2}{\s}\,\up{1}{\s} (\rf)  \,=\, \up{2}{\Phi}\, \up{1}{\Phi}
 \, (e\o\rf) \,  \up{1}{\Phi}\^{-1} \, \up{2}{\Phi}\^{-1} \,= \\ [2mm]
 =\,  F \, \D_a(\Phi) \, (e\o\rf) \, \D_a(\Phi^{-1}) \, F^{-1} \,=\,
 F \, \D_a(\Phi\,(e\o \rf)\,\Phi^{-1}) \, F^{-1} 
 \,=\, \D_F (\s(\rf)) \ . \er
$$

Let us finally discuss the computation of $F^*$. 
It is based on the second identity in (\ref{propS}) and on the 
relation  $(e \o S) (id \o \D)(S) = (S \o e) (\D \o id)(S)$
which can be checked in straightforward way.
Applying them and the property (\ref{cov}) to the 
definition (\ref{str1}), we derive
\ba 
    F^* & \vspace*{4mm}= & (\Ph{2}{}\Ph{1}{} \ \D_a(\Phi^{-1}))^*
     \, = \, \D'_a(\Phi \ S) \ (S^{-1})_{13} \Ph{1}{-1} (S^{-1})_{23}
            \Ph{2}{-1} \nn  \\[2mm]
      & = & \D'_a(\Phi)\ (\D' \o id)(S) \
           \left[ (e\o S^{-1})\, (id \o \D')( S^{-1}) \right]_{213}
           \Ph{1}{-1} \Ph{2}{-1} \nn\\[2mm]
      & = & \D'_a(\Phi)\ (\D' \o id)(S) \
           \left[ (id \o \D)( S^{-1})\,(e\o S^{-1}) \right]_{213}
           \Ph{1}{-1} \Ph{2}{-1} \nn\\[2mm]
      & = &  \D'_a(\Phi)\ (S\o \re) \ \Ph{1}{-1} \ \Ph{2}{-1} \,=\, 
           S_a \  \D_a(\Phi)\ \Ph{1}{-1} \ \Ph{2}{-1} \,=\
           S_a \ F^{-1}\ \ . \nn
\ea
The index on $[.]_{213}$ refers to a permutation of tensor
factors in the expression enclosed by the brackets.

All other relations in Proposition \ref{strucprop} are either
obvious or they follow directly from the derived equations.
This applies in particular to eq.~(\ref{qYB}).

\subsection{Proof of Proposition 3.}

In this subsection  we want to construct consistent structure
data for vertex operators of the deformed universal enveloping
algebras $\G = U_q(\sg)$ from their Clebsch-Gordan maps and
$6j$-symbols. To fix our notations, let us recall that the
Clebsch-Gordan maps $C[TL|S]: V^T \o  V^L \mapsto V^S$
have the following properties:
\ba
 C [TL|S]\, (\t^T \o \t^L)(\D(\xi)) & = & \t^S(\xi)\, C[TL|S]
       \ \ \  \mbox{ for all } \ \ \ \xi \in \G \ , \nn\\[2mm]
 (\frac{\k_S}{\k_T \k_L})^{\pm 1}\; C [TL|S]\, R_\pm^{TL}\, C[TL|R]^*
       & = &  \dl_{R,S}  \ \ \ \mbox{ with } \ \ \ \k_L \,=\, \t^L(\k)
        \nn \ ,\\[1mm]
 \sum_Q \Br{\pm}{Q}{S}{L}{T}{R}{J} \, C[JQ|R] \, C_{23} [TL|Q]
 & = &  C [JS|R] \, C_{23} [TL|S] \, (R_\pm^{JT} \o e^L)\  ,\nn\\[1mm]
 \sum_Q \Fus{Q}{S}{L}{T}{R}{J} \, C[QL|R] \, C_{12} [JT|Q]
 & = &  C [JS|R] \, C_{23} [TL|S] \ . \label{3jprop}
\ea
In  the order of appearance here these equations describe the
intertwining properties of the Clebsch-Gordan maps, their
normalization and the definition of the braiding and fusion
matrices (or $6j$-symbols). It is well known that certain
(polynomial) relations for the numbers $\Br{\pm}{\,.}{\,.}{.}
{.}{.}{.}$ and $\Fus{\,.}{\,.}{.}{.}{.}{.}$ follow from their
definitions and properties of the quasi-triangular
Hopf algebra $\G$ (see, e.g., \cite{KR1}). In
particular, one has
\ba
\sum_Q \Fus{Q}{S}{L}{T\ }{R}{J} \, \Fus{N}{R}{L}{Q\ }{P}{I} \,
    \Fus{M}{Q}{T}{J\ }{N}{I} & = &
    \Fus{M}{R}{S}{J\ }{P}{I} \, \Fus{N}{S}{L}{T}{P}{M}\ \ ,\nn \\[1mm]
 \sum_N \Fus{N}{S}{L}{T}{P}{M} \, \Fus{N}{R}{L}{T}{P}{M}^* & = &
    \dl_{S,R}  \  , \label{6jprop} \\[1mm]
 \Omega_+ \Omega_-^{-1} \vvert{N}{P}{T} \,
    \Br{-}{N}{N'}{L}{T}{P}{M} & = & \Br{+}{N}{N'}{L}{T}{P}{M} \,
 \Omega_+ \Omega_-^{-1}  \vvert{M}{N}{T} \  .\nn
\ea
Here we used the notation
\be{Om}
  \Omega_{\pm} \vvert{L}{P}{T} = \Br{\pm}{T}{L}{L}{T\ }{P}{0}
   \ \ \mbox{ and } \ \ \Omega_+ \Omega_-^{-1} \vvert{L}{P}{T} =
   \Omega_+ \vvert{L}{P}{T}\, \Omega_-^{-1}  \vvert{L}{P}{T} =
   \frac{v_L v_T}{v_P}  \ .
\ee
To proceed, we parameterize  the labels $N,M,I$ and  $N'$ in terms
of new variables $\lambda,\vth,\iota$ and $\vth'$ so that
$$
   N = P + \lambda \ \ , \ \ M = P + \lambda + \vth \ \ ,  \ \
   I = P + \lambda + \vth + \iota\  \ ,
   \ \ N' = P  + \vth'\  .
$$
Let us introduce  the following matrices $\C \{ TL|S\}(P)$,
$\R_{\pm}^{TL}(P)$ and $ D^T(P)$ ,
\ba
\C \{ TL|S\}(P)_{\vs, \vth \lambda} & = & \Fus{P+\lambda}{S}
   {L}{T}{P}{P+\vth+\lambda} \, \dl_{\vs,\vth+\lambda} \ ,\nn\\[1mm]
\R_{\pm}^{TL} (P)_{\vth' \lambda',\vth \lambda} & = &
   \Br{\pm}{P+\lambda}{P+\vth'}{L}{T}{P}{P+\vth+\lambda} \,
   \dl_{\vth+\lambda, \vth' +\lambda'} \ , \nn\\[1mm]
D^T(P)_{\vth',\vth} & = & \Omega_+ \Omega_-^{-1}
   \vvert{P+\vth}{P}{T} \, \dl_{\vth',\vth}\ \ . \label{DOm}
\ea
We may think of the $P$-dependent matrices $\C\{ TL|S\} (P)$,
$\R_\pm^{TL} (P) $ and $D^T(P)$ as matrices with entries in
the algebra $\C$. Whenever we do  so,  we will neglect to  write
the $P$-dependence explicitly and use the symbols $\C\{ TL|S\},
\R_\pm^{TL}$ and $D^T$. If  we  introduce in addition the matrix
valued map $\sigma^L$ on $\C$ by
$$ \sigma^L (P)_{\lambda', \lambda}   = (P  + \lambda)
     \  \dl_{\lambda',\lambda} \ , $$
then eqs.~(\ref{6jprop}) for the fusion and braiding matrices become
\ba
\sum_Q \Fus{Q}{S}{L}{T}{R}{J}\  \C\{ QL|R\}\ \s^L( \C_{12} \{ JT|Q\})
   & = & \C \{ JS|R\} \ \C_{23} \{ TL|S\} \  , \nn\\[1mm]
\C \{ TL|S\} \ \C\{ TL|R\}^* & = & \dl_{R,S}  \  , \nn\\[3mm]
(D^T\o e^L)\  \R_-^{TL}  & = & \R_+^{TL} \ \s^L(D^T) \ .\label{DR}
\ea
The last of these equations appears already as a close relative of
eq.~(\ref{ax3}). In fact, for semi-simple $\G$ one can
construct universal objects $\R_\pm \in \G \o  \G\o \C$ and
$ D \in \G \o \C$ so that $\R_\pm^{TL} = (\t^T \o \t^L)(\R_\pm)$
and $D^T = \t^T(D)$. Then eq.~(\ref{DR}) turns into
$$     \up{1}{D}\  \R_- = \R_+ \se{2} (D) \ \ , $$
where $\s : \C \mapsto \G \o \C$ is defined so that $\s^L =
(\t^L \o id) \circ \s$. When $\Omega_+ \Omega_-^{-1}$ is expressed
in terms of the ribbon element $v$ as in eq.~(\ref{Om}),
the definition (\ref{DOm}) of $D$  becomes
$$  D = \s(\rv) \rv^{-1} v_a \ \ . $$
To build the universal element $F \in \G \o \G \o \C$, we combine
the matrices $\C \{ ..|.\}$ with the Clebsch-Gordan maps so that $F^{TL}
= (\t^T \o \t^L) (F) $ is given by
$$  F^{TL} \equiv \sum_S \C \{ TL|S \} ^*\  C[TL|S]\ \ . $$
Multiplying the adjoint of the first eq.~in (\ref{6jprop}) with
eq.~(\ref{3jprop}), taking sum over  $R,S$ and using
the intertwining  properties  of the Clebsch-Gordan maps, we
obtain the equation (\ref{ax1}) for $F$. In  the same way,
one may combine the normalizations for the Clebsch-Gordan maps
$C[TL|S]$ and the matrices $\C \{ TL|S\}$ to derive that
$$ 
  F\,\left( ( \D(\k^{-1}) (\k \o \k)\, R_+^{-1} )\o \re \right) \, 
    F^* \,=\, e\o e \o \re \ ,
$$
and hence $F$ has the required property under the $*$-operation.
{}Finally, we use that the matrix $\C \{ TL|S\} (P)$ is proportional
to $\dl_{\vs, \vth+\lambda}$  so that
$$  
  (\rf(P+\vs) )\ \C \{TL|S\}(P)_{\vs, \vth \lambda} \,=\,
  \C\{TL|S\}(P)_{\vs,\vth \lambda}\  (\rf(P +\vth + \lambda)) \ .
$$
Here $\rf(P)$ is  an arbitrary function of  $P$, i.e., $\rf$ may  be
regarded  as an element in $\C$. With our standard conventions, this
can be stated as a matrix equation
$$   \s^S (\rf) \ \C\{TL|S\}\,=\, \C\{TL|S\}\ \s^L \s^T (\rf) \ . $$
Keeping in mind that $\D_F (\xi) = F(\D  (\xi) \o \re) F^{-1}$,
we discover eq.~(\ref{ax2}). All other properties of the structure
data follow easily from the relations we have discussed here.

\subsection{Structure data for left chiral vertex operators.}
We can obtain the relations for the left structure data
$F_l, \s_l, D_l, \RR_\pm^l$ from eqs.~(\ref{Dop})-(\ref{*D})
if we substitute
\def\hs{\hspace*{.5cm}}
$$   \begin{array}{lclclclclcl}
 \D_a  & \rar & \D'_a & \ \ \ & F & \rar & F_l' & \ \ \ &  
    M  & \rar & (M^l)^{-1}  \\[1mm]
 R_\pm & \rar & R'_\pm && \s & \rar & \s_l  &&   
    \RR_\pm & \rar & (\RR^{l}_\pm)' \\[1mm]
  v_a  & \rar & v_a  && D & \rar & D_l^{-1} &&
    \,\Phi & \rar & \ \Phi^l \ \ .
      \end{array} $$
The prime on $F_l, \RR_\pm^l \in \G_a \o \G_a \o \C^l$ denotes
permutation of the first two tensor factors. Once the
validity of these rules has been checked for the defining
relations (\ref{str1})-(\ref{str3}) of structure data, we can apply
them to eqs.~(\ref{Dop})-(\ref{*D}). Within the notations of
Proposition \ref{strucprop}, the result looks as follows:
\ba \label{lDop} 
      D_l \, =\, (v^{-1}_a \rv)\, \s_l({\rv}^{-1})\ \ &,&
      \label{lax2} \ \ \up{1}{\s}\_l\up{2}{\s}\_l\,(\rf_l)  \, = \,
      \D_{F_l} (\s_l(\rf_l)) \ ,    \\[2mm]
      \label{lax1}
     {} \left[F_l \o e \right]_{1243}\; \Bigl( (\D_a \o id )(F_l)\Bigr) &=&
      \up{1}{\s}(F_l) \;
      \,\Bigl( (id \o \D_a)(F_l)\Bigr) \   , \\[2mm]
      \label{lax3}
      \up{2}{D_l}\,{\RR}^l_+ \, =\, \RR^l_- \,\up{1}{\s}\_l(D_l)\ \
      \  &, &  \ \
      \RR^l_+ \up{1}{D_l}\, =\, \up{2}{\s}\_l(D_l)\,\RR_-^l \ ,\\[2mm]
      \RR^l_{\pm,23}\ \up{2}{\s}\_l(\RR^l_{\pm,13}) \
      \RR^l_{\pm,12}  & = &
      \up{3}{\s}\_l( \RR^l_{\pm,12})\ \RR^l_{\pm,13} \
      \up{1}{\s}\_l(\RR^l_{\pm,23}) \ ,   \label{lqYB} \\[2mm]
      \label{l*F}
      F_l^* = S^{-1}_a F_l^{-1} \ \hs , \ \hs
      \label{l*RR}  \label{l*D}
      (\RR^l_\pm)^* &=& (\RR^l_\pm)^{-1} \ \hs
       , \ \hs D_l^* \;=\; D_l^{-1} \ , \\[2mm]
      \label{l*DF} \label{l*s}  \s_l(\rf_l )^* = \s_l(\rf_l^*)
      \hs  &,& \ \hs \D_{F_l} (\xi)^*\;=\;
      \D_{F_l}(\xi^*) \ .
\ea
Here $\D_{F_l} (\xi) \equiv F_l (\D(\xi)\o\re) F_l^{-1}$
$\in\G_a\o\G_a\o\V_l$ analogously to definition (\ref{DelF}).
Using the fact that $(\S_a \o \S_a)(\D_a (\xi )) =
\D_a'(\S_a(\xi))$, it is simple to see that the objects
(\ref{lstruct}), (\ref{rltrafo}) satisfy the relations
(\ref{lDop})-(\ref{l*s})
if  $F,\s,\D,\RR$ solve eqs.~(\ref{Dop})-(\ref{*D}).
The former equations are actually obtained from the latter by 
acting with the maps $\S^{(n)}_{lr}$ defined in eq.~(\ref{Sdef}).

\subsection{Properties of the field $g$.} The consistency of
the object $g$ with the constraint to the diagonal subspace
$\H$ is certainly its most important property. It was formulated
more precisely in (\ref{fSf}) using notations from Subsection
3.4. A formal proof may be given as follows. Suppose that
$ \rf g = \S_{lr}(\rf) g $ holds for all $\rf \in \C^r$. Then
one finds for arbitrary $\rf \in \C^r$
\ba
   g\ \rf & = & \S_a(\Phi^l)\Phi^r \ \rf \ = \
                 \S_a(\Phi^l)\ \s_r(\rf) \ \Phi^r \ = \
               \S_a\left((\S^{-1} \o id)(\s_r(\rf))\Phi^l\right)\ \Phi^r
                 \nn \\[1mm]
  & = & \S_a\left((\S^{-1} \o \S_{lr})(\s_r(\rf))\Phi^l\right)\ \Phi^r
  \ = \ \S_a\left(\s_l(\S_{lr}(\rf))\,\Phi^l\right)\ \Phi^r \nn \\[1mm]
           & = & \S_a(\Phi^l) \ \S_{lr}(\rf) \ \Phi^r \ = \
           \S_a(\Phi^l) \ \Phi^r \ \S_{lr}(\rf) 
           \ = \ g \ \S_{lr}(\rf)\ \ .   \nn
\ea
The computation makes use of the choice (\ref{rltrafo}) to replace
$\s_r$ by $\s_l$.
{}From now on, we think of $g$ as being restricted to $\H$. We begin
our proof of Proposition \ref{gprop} with the operator product
expansion (\ref{qg}) of $g \in \G_a \o \End(\H)$,
\ba
  \up{2}{g}\ \up{1}{g} &=& \S_a(\up{2}{\Phi}\^l) \ \S_a(\up{1}{\Phi}\^l)
               \  \up{2}{\Phi}\^r \, \up{1}{\Phi}\^r
   \,=\, (\S \o \S \o id) (F'_l \D'_a(\Phi^l))\ F_r\  \D_a(\Phi^r) \
              \nn \\[1mm]
   & = & (\S \o \S \o id) \left( (\S^{-1} \o \S^{-1} \o id )
      (F_r^{-1}) \D'_a(\Phi^l)\right) \ F_r\  \D_a(\Phi^r) \ \nn \\[1mm]
   & = & \D_a(\S_a(\Phi^l))\ F_r^{-1} F_r\ \D_a(\Phi^r) \,=\
   \D_a(g) \nn \ .
\ea
Here $\S_a(\up{1}{\Phi})$ and $\S_a(\up{2}{\Phi})$ are shorthands for
$(\S_a\o id)(\up{1}{\Phi})$ and $(id\o\S_a)(\up{2}{\Phi})$, respectively.
The exchange relations for $g$, the formula $g^{-1} = \S_a(g)$ and the
normalization $\e_a (g) = \re$ follow immediately from the functoriality
relation in (\ref{qg}).

The exchange relations (\ref{qs}) are derived from (\ref{rcov}),
(\ref{lcov}) and the explicit construction of $g$ as a product of
$\S_a(\Phi^l)$ and $\Phi^r$. Let us check the first of them:
$$
 \up{1}{M}\^r \up{2}{g} R_-  =\, \up{1}{M}\^r \, \S_a(\up{2}{\Phi}\^l)
 \, \up{2}{\Phi}\^r \, R_-  =\, \S_a(\up{2}{\Phi}\^l)\, \up{1}{M}\^r \,
 \up{2}{\Phi}\^r \, R_-  =\, \S_a(\up{2}{\Phi}\^l)\, \up{2}{\Phi}\^r \,
 R_+ \,\up{1}{M}\^r \,=\, \up{2}{g}\, R_+ \up{1}{M}\^r \, .
$$
The second relation in (\ref{qs}) can be obtained similarly if we take
the covariance properties (\ref{SPcov}) of $\S_a(\Phi^l)$ into
account.

Verification of the relations (\ref{qr}) makes use of the equality
$\rv_l = \rv_r$ which is valid on $\H$ and follows with the help of
$\S(v) = v$, if the constraint $\S_{lr}(\rf_r)=\rf_l$ is evaluated
on the ribbon element. The second relation in (\ref{qr}) then is
obvious, and for the first we check explicitly:
\ba
 g_n\,M^r_n &=& g_n\,(\Phi^r_n)^{-1}\,v_a^{-1}\,D_r\,\Phi^r_n \,=\,
 \S_a(\Phi^l_n)\,v_a^{-1}\,D_r \,\Phi^r_n \nn \\ [1mm]
 &=& (\Phi^l_n)^{-1} \,\theta_l \, \rv_r^{-1}\,\s_r(\rv_r)\,\Phi^r_n \,=\,
 \rv_l^{-1}\,(\Phi^l_n)^{-1} \,\theta_l \,\Phi^r_n \, \rv_l \,=\,
 (\Phi^l_n)^{-1} \,\s_l(\rv_l^{-1})\,\rv_l \,\theta_l\,\Phi^r_n \nn \\ [1mm]
 &=& (\Phi^l_n)^{-1} \,v_a\,D_l\,\Phi^l_n \,(\Phi^l_n)^{-1} \,\theta_l\,
 \Phi^r_n  \,=\,M^l_n \, \S_a(\Phi^l_n)\,\Phi^r_n \,=\ M^l_n \, g_n \ .\nn
\ea
In this computation it was convenient to insert the formula (\ref{theta})
which expresses $\S_a(\Phi^l)$ in terms of $(\Phi^l)^{-1}$.

\subsection{Properties of lattice vertex operators.} Let us prove that
the structure constants $F_\a$, $\RR^\a_\pm$ and $\s_\a$ appearing in
eqs.~(\ref{DFnp})-(\ref{sFn}) coincide with those of the vertex operators
$\Phi^\a_0$. To this end, we exploit the construction of $\Phi^\a_n$
as a product of $\Phi^\a_0$ and the holonomies $U^\a_n$ (see eq.
(\ref{Pn})). Equation (\ref{sFn}) is actually obvious, since $\rf_\a$
commute with all the elements $U^\a_n$. Furthermore, eq.~(\ref{PPn})
is a simple consequence of (\ref{DFnp}). Hence, we need to prove only
eq.~(\ref{DFnp}) which we do for the right chirality ( the left one
works analogously),
$$
\ar{c}
 \D_a(\Phi_n^r) \,=\, \D_a(\Phi_0^r)\,\D_a(U_n^r)\,= \,
  F_r^{-1}\,\up{2}{\Phi}\^r_0 \,\up{1}{\Phi}\^r_0\,
 R_-^{-1}\,\up{2}{U}\^r_n\,\up{1}{U}\^r_n \, =  \\ [1mm]
 =\, F_r^{-1}\,\up{2}{\Phi}\^r_0\,\up{2}{U}\^r_n\,
 \up{1}{\Phi}\^r_0\,\up{1}{U}\^r_n \,=\,
 F_r^{-1}\,\up{2}{\Phi}\^r_n\,\up{1}{\Phi}\^r_n\ .  \er
$$

The exchange relations (\ref{PJn})-(\ref{RL}) are established by
induction. Indeed, for $n=0$ they are part of Definition \ref{W_N}.
Assume now that eqs.~(\ref{PJn})-(\ref{RL}) hold for a certain $n<N$
so that, in particular, $\Phi^\a_n$ has non-trivial exchange
relations with $J^\a_n$ and $J^\a_{n+1}$ only. Then $\Phi^\a_{n+1}\,=\,
\Phi^\a_{n} J^\a_{n+1}$ necessarily commutes
with all currents $J_m^\a$ except from $J^\a_n$, $J^\a_{n+1}$ and
$J^\a_{n+2}$. It is easy to verify that the exchange relations
with $J^\a_n$ become trivial as well. We demonstrate this for
$\a = r$:
$$
 \up{1}{\Phi}\^r_{n+1} \, \up{2}{J}\^r_n \,=\,
 \up{1}{\Phi}\^r_{n} \, \up{1}{J}\^r_{n+1} \, \up{2}{J}\^r_n \,=\,
 \up{1}{\Phi}\^r_{n} \, \up{2}{J}\^r_n \, R_+ \, \up{1}{J}\^r_{n+1} \,=\,
 \up{2}{J}\^r_n \, \up{1}{\Phi}\^r_{n} \, \up{1}{J}\^r_{n+1} \,=\,
 \up{2}{J}\^r_n \, \up{1}{\Phi}\^r_{n+1} \ .
$$
It can be checked similarly that the relations between
$\Phi^\a_{n+1}$ and $J^\a_{n+1}$, $J^\a_{n+2}$ coincide with those
between $\Phi^\a_{n}$ and $J^\a_{n}$, $J^\a_{n+1}$. This
completes the induction.

Now we have to prove eqs.~(\ref{braid})-(\ref{braid'}). For instance,
using (\ref{UU}) and (\ref{PJn}), we derive the first relation in
(\ref{braid}) for $0\leq n<m<N$ :
$$
 \ar{c} \up{1}{\Phi}\^r_n\,\up{2}{\Phi}\^r_m \,=\,
 \up{1}{\Phi}\^r_0\,\up{1}{U}\^r_n \,
 \up{2}{\Phi}\^r_0\,\up{2}{U}\^r_m \,=\,
 \up{1}{\Phi}\^r_0\,\up{2}{\Phi}\^r_0 \,R_+\,
 \up{1}{U}\^r_n\,\up{2}{U}\^r_m \,=  \\ [1.5mm]
 = \RR^r_- \,\,\up{2}{\Phi}\^r_0\,\up{1}{\Phi}\^r_0\,R_-^{-1}\,
 \up{2}{U}\^r_m\,\up{1}{U}\^r_n \,= \RR^r_- \,
 \up{2}{\Phi}\^r_0\,\,\up{2}{U}\^r_m \,\up{1}{\Phi}\^r_0\,\up{1}{U}\^r_n
 \,=  \RR^r_- \, \up{2}{\Phi}\^r_m\,\up{1}{\Phi}\^r_n \ .\er
$$
The relations (\ref{*Pn})-(\ref{sFn}), (\ref{MDMn})-(\ref{RL}) which
involve vertex operators $\Phi^\a_n$ outside of the initial interval
$n=0,..,N-1$, are derived with the help of eqs.~(\ref{DMn})-(\ref{PMn})
for the
monodromies $M^\a_n$. Since the derivation uses the same technique
as above, we prove only the functoriality relation for $\Phi^r_{n}$.
As a first step, we check the following:
$$ \ar{c}
 \D_a(\Phi^r_{n+N})\,=\,\D_a(\Phi^r_n)\,\D_a(M^r_n) \,=\,
 F_r^{-1}\, \up{2}{\Phi}\^r_{n} \, \up{1}{\Phi}\^r_{n} \,
 R_-^{-1}\,\up{2}{M}\^r_n\,R_+\up{1}{M}\^r_n  \,= \\ [1mm]
 = F_r^{-1}\, \up{2}{\Phi}\^r_{n} \,\up{2}{M}\^r_n\,
 \up{1}{\Phi}\^r_{n}\,\up{1}{M}\^r_n\,  \,=\,
 F_r^{-1}\, \up{2}{\Phi}\^r_{n+N} \, \up{1}{\Phi}\^r_{n+N} \ . \er
$$
Then we use an induction and eq.~(\ref{PMn}) to get the same property
for $\Phi^r_{n+kN}$.

{}Finally, we establish relations (\ref{PnN})
directly with the help of eq.~(\ref{PMn}):
$$
 \RR^r_+ \, \up{2}{\Phi}\^r_{n}\,\up{1}{\Phi}\^r_{n+N}\,=\, \RR_+ \,
 \up{2}{\Phi}\^r_{n}\,\up{1}{\Phi}\^r_{n}\,\up{1}{M}\^r_n\,=
 \up{1}{\Phi}\^r_{n}\,\up{2}{\Phi}\^r_{n}\, R_+\,\up{1}{M}\^r_n\,=\,
 \up{1}{\Phi}\^r_{n}\,\up{1}{M}\^r_n\, \up{2}{\Phi}\^r_{n}\, R_- \,=\,
 \up{1}{\Phi}\^r_{n+N}\, \up{2}{\Phi}\^r_{n}\, R_- \ \ .
$$
Detailed computations for the other relations in (\ref{*Pn})-(\ref{sFn}),
(\ref{MDMn})-(\ref{RL}) can be worked out easily.

\subsection{Properties of lattice field $g_n$.}
Let us notice that the equality (\ref{theta}) holds in the lattice case
for all $\Phi^l_n$ with the same $\theta_l\in\G_a\o\C^l$ (as we explained
in Subsection 5.2, vertex operators of the same chirality assigned to
different sites possess the same structure data). Therefore, we can
proceed as in the toy model case and rewrite the expression for $g_n$
as follows:
$$ g_n \,=\, (\Phi^l_n)^{-1} \,\theta_l \, \Phi^r_n \ .$$
This relation allows to express $g_n$ in terms of $g_0$ and the
holonomies $U^\a_n \in \G_a \o \K_N$:
\be{gng0}
 g_n \,=\, (U^l_n)^{-1} \, g_0 \, U^r_n \ .
\ee
Bearing in mind that elements from $\K_N$, and hence, in particular,
components of the holonomies $U^\a_n$, leave the subspaces
$W_N^{\bar K K}$ of the full representation space $\cM_N =
\bigoplus_{IJ} W_N^{IJ}$ invariant (see Subsection 5.4), the
equality (\ref{gng0}) explains why all $g_n$ can be restricted
on the diagonal subspace $\H_N = \bigoplus_K W_N^{\bar K K } $
simultaneously.

Among the properties of the lattice field $g_n$ in Proposition
\ref{prgn}, only the relations (\ref{ign}) and (\ref{loc}) have
not been considered in the toy model case. Eqs.~(\ref{ign}) follow
immediately from the covariance properties (\ref{xPn}) of the vertex
operators and the remark that, due to eq.~(\ref{theta}),  the second
relation in (\ref{xPn}) can be rewritten as follows:
$$
\S_a(\Phi^l_n)\,\iota_n(\xi) \,=\, \D_n(\xi)\,\S_a(\Phi^l_n)\ , \ \ \ \
  \S_a(\Phi^l_n)\,\iota_m(\xi) \,=\, \iota_m(\xi)\,\S_a(\Phi^l_n)\ \
 {\rm for}\ m\neq n \, (mod N)
$$
for all $\xi\in\G$.

The periodicity of $g_n$ is derived with the help of relations
(\ref{DD}), (\ref{MDMn}) and the second equation in (\ref{qr}):
\ba
 g_{n+N} &=& \S_a(\Phi^l_{n+N}) \, \Phi^r_{n+N} \,=\,
 (\Phi^l_{n+N})^{-1} \, \theta_l \, \Phi^r_{n+N} \,=\,
 (\Phi^l_{n}\,M^l_n)^{-1}\, \theta_l \, \Phi^r_{n}\,M^r_n \nn \\ [1mm]
 &=& (\Phi^l_{n})^{-1} \, (v_a\,D_l)^{-1} \, \theta_l \,
 v_a^{-1} \, D_r\,\Phi^r_{n} \,=\, (\Phi^l_{n})^{-1} \,
 \s_l(\rv_l)\,\rv_l^{-1}\,\theta_l\, \s_r(\rv_r)\,\rv_r^{-1} \,\Phi^r_{n}
 \nn \\ [1mm]  &=& \rv_l\, (\Phi^l_{n})^{-1} \,
 \rv_l^{-1}\,\theta_l \,\rv_r^{-1} \, \Phi^r_{n}\,\rv_r  \,=\,
 \rv_r\, (\Phi^l_{n})^{-1} \, \rv_r^{-1}\,\theta_l \,\rv_l^{-1} \,
 \Phi^r_{n}\,\rv_l \,=\, \S_a(\Phi^l_{n}) \, \Phi^r_{n} \ =\,  g_n \ .\nn
\ea
Due to periodicity, it is sufficient to check the locality
of $g_n$ only for $0\leq n,m <N$. Taking, for definiteness, $n<m$,
we derive:
$$
 \ar{c}  \up{1}{g}_n \, \up{2}{g}_m \,=\,
 (\up{1}{U}\^l_n)^{-1} \, \up{1}{g}\_0 \, \up{1}{U}\^r_n \,
 (\up{2}{U}\^l_m)^{-1} \, \up{2}{g}\_0 \, \up{2}{U}\^r_m \,=\,
 (\up{1}{U}\^l_n)^{-1} \, (\up{2}{U}\^l_m)^{-1} \, R_-^{-1} \,
 \up{1}{g}\_0 \, \up{2}{g}\_0 \, R_+ \, \up{1}{U}\^r_n \,
 \up{2}{U}\^r_m \,=  \\ [1mm]  =\,
 (\up{2}{U}\^l_m)^{-1} \, (\up{1}{U}\^l_n)^{-1} \, R_+ \,
 \up{2}{g}\_0 \, \up{1}{g}\_0 \, R_-^{-1} \, \up{2}{U}\^r_m \,
 \up{1}{U}\^r_n \,=\, (\up{2}{U}\^l_m)^{-1} \, \up{2}{g}\_0 \,
 (\up{1}{U}\^l_n)^{-1} \, \up{2}{U}\^r_m \, \up{1}{g}\_0 \,
 \up{1}{U}\^r_n \,=\ \up{2}{g}_m \, \up{1}{g}_n \ . \er
$$
Here we used eq.~(\ref{gng0}) and the commutation relations between
$g_0$ and the holonomies $U^\a_n$ which are obvious consequences of
eqs.~(\ref{gJ}).

\newcommand{\sbibitem}[1]{\vspace*{-1.5ex} \bibitem{#1}}

\end{document}